\documentclass[11pt]{article}
\usepackage[utf8]{inputenc}
\usepackage{geometry}

\usepackage[dvipsnames]{xcolor}

\usepackage{epsfig}
\usepackage{graphicx}
\usepackage{amsbsy}
\usepackage{amsmath}
\usepackage{amsfonts}
\usepackage{amssymb}
\usepackage{textcomp}
\usepackage{hyperref}
\usepackage{aliascnt}
\usepackage{bbm}
\usepackage{tikz-cd}
\usepackage{enumitem}\setlist[enumerate]{label=\textnormal{(\roman*)}}

\usepackage{quiver}

\usepackage[section]{placeins}

\renewcommand{\leq}{\leqslant}
\renewcommand{\geq}{\geqslant}
\renewcommand{\le}{\leqslant}
\renewcommand{\ge}{\geqslant}
\newcommand{\pmo}{^{\pm1}}

\newcommand{\mcm}[3]{\newcommand{#1}[#2]{{\ensuremath{#3}}}} 
\mcm{\tuple}{1}{\langle #1 \rangle}
\mcm{\name}{1}{\ulcorner #1 \urcorner}
\mcm{\Nbb}{0}{\mathbb{N}}
\mcm{\N}{0}{\mathbb{N}}
\mcm{\Zbb}{0}{\mathbb{Z}}
\mcm{\Rbb}{0}{\mathbb{R}}
\mcm{\Cbb}{0}{\mathbb{C}}
\mcm{\Qbb}{0}{\mathbb{Q}}
\mcm{\Sbb}{0}{\mathbb{S}}
\mcm{\Hbb}{0}{\mathbb{H}}
\mcm{\Kbb}{0}{\mathbb{K}}
\mcm{\Acal}{0}{\cal A}
\mcm{\Bcal}{0}{\mathcal B}
\mcm{\Ccal}{0}{\cal C}
\mcm{\Dcal}{0}{\cal D}
\mcm{\Ecal}{0}{\cal E}
\mcm{\Fcal}{0}{\cal F}
\mcm{\Gcal}{0}{\cal G}
\mcm{\Hcal}{0}{\cal H}
\mcm{\Ical}{0}{\cal I}
\mcm{\Jcal}{0}{\cal J}
\mcm{\Kcal}{0}{\cal K}
\mcm{\Lcal}{0}{\cal L}
\mcm{\Mcal}{0}{\cal M}
\mcm{\Ncal}{0}{\cal N}
\mcm{\Ocal}{0}{{\cal O}}
\mcm{\Pcal}{0}{{\cal P}}
\mcm{\Qcal}{0}{{\cal Q}}
\mcm{\Rcal}{0}{{\cal R}}
\mcm{\Scal}{0}{{\cal S}}
\mcm{\Tcal}{0}{{\cal T}}
\mcm{\Ucal}{0}{{\cal U}}
\mcm{\Vcal}{0}{{\cal V}}
\mcm{\Wcal}{0}{{\cal W}}
\mcm{\Xcal}{0}{{\cal X}}
\mcm{\Ycal}{0}{{\cal Y}}
\mcm{\Zcal}{0}{{\cal Z}}
\mcm{\Mfrak}{0}{\mathfrak M}

\usepackage{amsthm}

\newcommand{\theoremize}[2]{\newaliascnt{#1}{thm} \newtheorem{#1}[#1]{#2} \aliascntresetthe{#1}}

\theoremstyle{plain}
\newtheorem{thm}{Theorem}[section]
\theoremize{lem}{Lemma}

\theoremize{skolem}{Skolem}
\theoremize{fact}{Fact}
\theoremize{setting}{Setting}

\newtheorem{claim}{Claim}[thm]
\theoremize{obs}{Observation}
\theoremize{prop}{Proposition}
\theoremize{cor}{Corollary}
\theoremize{que}{Question}
\theoremize{oque}{Open Question}
\theoremize{con}{Conjecture}

\theoremstyle{definition}
\theoremize{dfn}{Definition}
\theoremize{rem}{Remark}
\theoremize{eg}{Example}
\theoremize{exercise}{Exercise}
\theoremstyle{plain}

\hypersetup{
    draft = false,
    bookmarksopen=true,
    colorlinks,
    linkcolor={red!60!black},
    citecolor={green!60!black},
    urlcolor={blue!60!black}
}

\newcommand{\sm}{\setminus}

\makeatletter
\def\thm@space@setup{\thm@preskip=2pt
\thm@postskip=2pt}
\makeatother

\makeatletter
\newenvironment{pf}[1][\proofname]{\par
  \pushQED{\qed}%
  \normalfont \topsep0\p@\relax
  \trivlist
  \item[\hskip\labelsep\itshape
  #1\@addpunct{.}]\ignorespaces
}{%
  \popQED\endtrivlist\@endpefalse
}
\makeatother

\newenvironment{cproof}{\noindent\textit{Proof of Claim.}}{\phantom{!}\!\hfill$\diamondsuit$\vspace{2pt}}

\usepackage{titlesec}

\titleformat{\section}
  {\normalfont\scshape\filcenter}{\thesection}{1em}{}

\titleformat{\subsection}
 {\normalfont\scshape\filcenter}{\thesubsection}{1em}{}

\titlespacing*{\section}
{0pt}{5pt plus 1pt minus 1pt}{5pt plus 1pt minus 1pt}

\titlespacing*{\subsection}
{0pt}{5pt plus 1pt minus 1pt}{5pt plus 1pt minus 1pt}
\usepackage{tikzit}

\tikzstyle{black dot}=[fill=black, draw=none, shape=circle, scale=0.3]
\tikzstyle{backing}=[fill=white, draw=none, shape=circle, scale=0.85]
\tikzstyle{small back}=[fill=white, draw=none, shape=circle, scale=0.95]
\tikzstyle{medium back}=[fill=white, draw=none, shape=circle, scale=1.2]

\tikzstyle{c0}=[-, line width=0.08em]
\tikzstyle{c0 arrow}=[->, line width=0.08em]
\tikzstyle{c0 bold}=[-, line width=0.15em]
\tikzstyle{c1}=[-, draw={rgb,255: red,135; green,206; blue,235}, line width=0.08em]
\tikzstyle{c1 arrow}=[->, draw={rgb,255: red,135; green,206; blue,235}, line width=0.08em]
\tikzstyle{c1 bold}=[-, draw={rgb,255: red,135; green,206; blue,235}, line width=0.15em]
\tikzstyle{c2}=[-, draw={rgb,255: red,255; green,69; blue,0}, line width=0.08em]
\tikzstyle{c2 arrow}=[->, draw={rgb,255: red,255; green,69; blue,0}, line width=0.08em]
\tikzstyle{c3}=[-, draw={rgb,255: red,102; green,176; blue,50}, line width=0.08em]
\tikzstyle{c3 arrow}=[->, draw={rgb,255: red,102; green,176; blue,50}, line width=0.08em]
\tikzstyle{d1}=[-, dashed, line width=0.08em]
\tikzstyle{d1 arrow}=[->, line width=0.08em, dashed]
\tikzstyle{c1 dashed}=[-, draw={rgb,255: red,135; green,206; blue,235}, line width=0.08em, dashed, dash pattern=on 0.5em off 0.5em]
\tikzstyle{c2 dashed}=[-, draw={rgb,255: red,255; green,69; blue,0}, line width=0.08em, dashed, dash pattern=on 0.5em off 0.5em]
\tikzstyle{sideby}=[-, draw={rgb,255: red,255; green,255; blue,255}, line width=0.08em, side by side={blue!50:red!50}]
\tikzstyle{thickwhitedash}=[-, draw={rgb,255: red,255; green,255; blue,255}, line width=0.24em, dashed, dash pattern=on 0.5em off 0.5em]
\tikzstyle{triple}=[-, draw={rgb,255: red,102; green,176; blue,50}, line width=0.08em, side by side by side={blue!50:red!50}]
\tikzstyle{reallythickwhitedash}=[-, draw={rgb,255: red,255; green,255; blue,255}, line width=0.32em, dashed, dash pattern=on 0.5em off 0.5em]

\usetikzlibrary{nfold}
\makeatletter
\tikzset{
  side by side/.style args={#1:#2}{
    postaction={path only,draw=#1,offset=+.8\pgflinewidth},
    postaction={path only,draw=#2,offset=+-.8\pgflinewidth}},
  side by side'/.style={path only,side by side={#1}},
  offset/.code=
    \tikz@addoption{%
      \pgfgetpath\tikz@temp
      \pgfsetpath\pgfutil@empty
      \pgfoffsetpath\tikz@temp{#1}}}
\makeatother

\makeatletter
\tikzset{
  side by side by side/.style args={#1:#2}{
    postaction={path only,draw=#1,offset=+1.5\pgflinewidth},
    postaction={path only,draw=#2,offset=+-1.5\pgflinewidth}},
  side by side by side'/.style={path only,side by side by side={#1}},
  offset/.code=
    \tikz@addoption{%
      \pgfgetpath\tikz@temp
      \pgfsetpath\pgfutil@empty
      \pgfoffsetpath\tikz@temp{#1}}}
\makeatother

\title{\vspace{-1cm} Graphic matroids with connectivity constraints}
\author{Will J.\ Turner\thanks{TU Freiberg, funded by DFG, project number 546892829.}}
\date{}

\begin{document}

\maketitle
\thispagestyle{empty}

\setcounter{page}{0}

\vspace{-0.75cm}
\begin{center}
\begin{minipage}{0.9\textwidth}
\small
\textbf{Abstract.} In 1959, Tutte proved that a matroid can be realised by a graph if and only if it does not have a minor isomorphic to $U_{2,4}$, $F_7$, $F_7^*$, $M^*(K_5)$ or $M^*(K_{3,3})$. In this paper, we extend this theorem to determine when it is possible to find a graphic realisation $G$ of a matroid $M$, subject to a list of additional constraints of the following form: a subset $F$ of the elements of $M$ must induce a connected subgraph of $G$. We prove that a connected matroid with such constraints---called a `constraint matroid'---is realisable as a graph satisfying those constraints if and only if it does not contain as a minor one of the constraint matroids from a given list: that is, the original $5$ matroids of Tutte, a further $9$ sporadic constraint matroids, and two well-behaved infinite families.

\hspace{0.4cm} Matroids realisable as graphs appear famously in Whitney's planarity criterion: Carmesin proved an extension of Whitney's planarity criterion to the embeddability of local simply connected $2$-complexes into $\mathbb{R}^3$. An application of our theorem is an extension of this theorem beyond the simply connected regime.
\end{minipage}
\end{center}
\vspace{1.8cm}
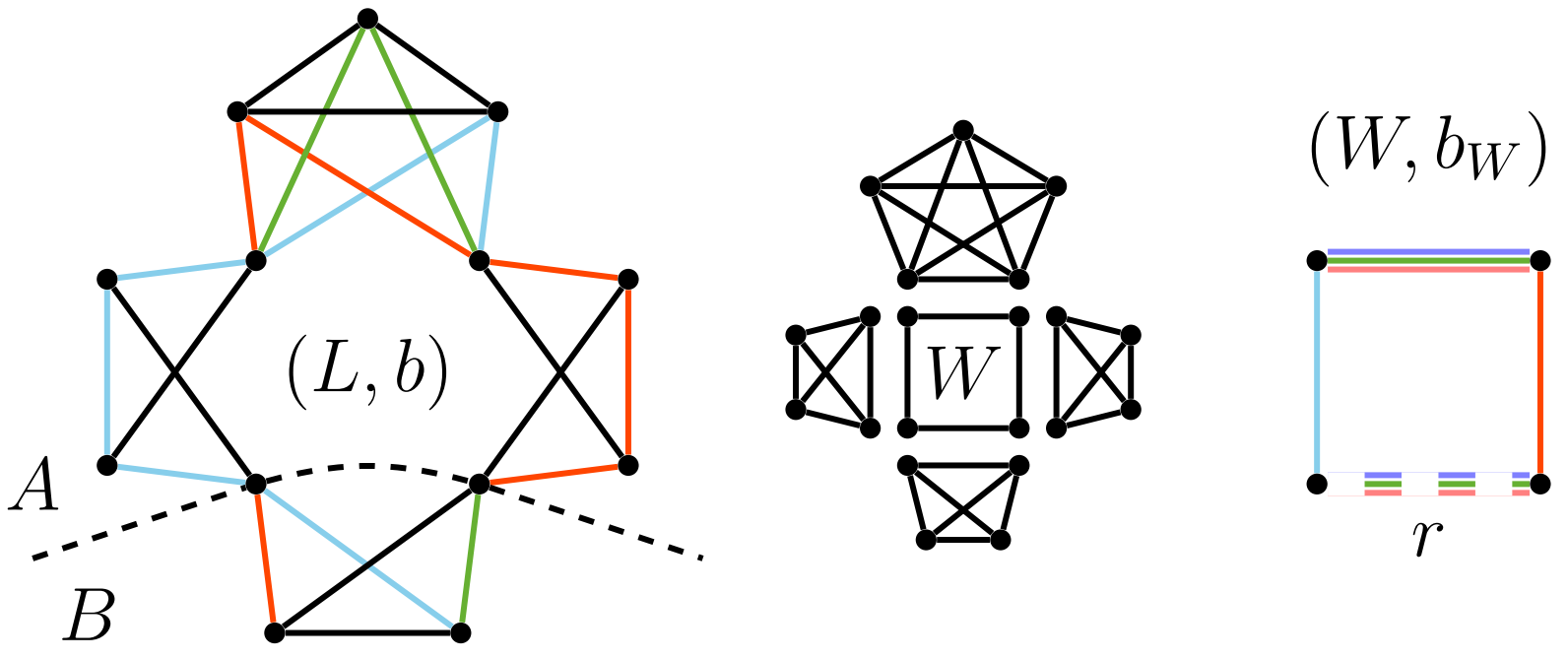
\begin{center}
\begin{minipage}{\textwidth}
        \centering
        \Large
        \scalebox{1}{\begin{tikzpicture}
	\begin{pgfonlayer}{nodelayer}
		\node [style=black dot] (169) at (-24, 6) {};
		\node [style=black dot] (170) at (-24, 3) {};
		\node [style=black dot] (171) at (-27, 3) {};
		\node [style=black dot] (172) at (-27, 6) {};
		\node [style=black dot] (174) at (-23.75, 8) {};
		\node [style=black dot] (175) at (-22, 5.75) {};
		\node [style=black dot] (176) at (-22, 3.25) {};
		\node [style=black dot] (177) at (-24.25, 1) {};
		\node [style=black dot] (178) at (-26.75, 1) {};
		\node [style=black dot] (179) at (-29, 3.25) {};
		\node [style=black dot] (180) at (-29, 5.75) {};
		\node [style=black dot] (197) at (-27.25, 8) {};
		\node [style=black dot] (198) at (-25.5, 9.25) {};
		\node [style=black dot] (199) at (-18.75, 3.75) {};
		\node [style=black dot] (200) at (-18.75, 5.25) {};
		\node [style=black dot] (201) at (-19.75, 4) {};
		\node [style=black dot] (202) at (-19.75, 5) {};
		\node [style=black dot] (203) at (-18.25, 3.25) {};
		\node [style=black dot] (204) at (-18, 2.25) {};
		\node [style=black dot] (205) at (-17, 2.25) {};
		\node [style=black dot] (206) at (-16.75, 3.25) {};
		\node [style=black dot] (207) at (-16.25, 3.75) {};
		\node [style=black dot] (208) at (-16.25, 5.25) {};
		\node [style=black dot] (209) at (-15.25, 5) {};
		\node [style=black dot] (210) at (-15.25, 4) {};
		\node [style=black dot] (211) at (-18.25, 5.75) {};
		\node [style=black dot] (212) at (-16.75, 5.75) {};
		\node [style=black dot] (213) at (-18.75, 7) {};
		\node [style=black dot] (214) at (-16.25, 7) {};
		\node [style=black dot] (215) at (-17.5, 7.75) {};
		\node [style=black dot] (216) at (-18.25, 3.75) {};
		\node [style=black dot] (217) at (-18.25, 5.25) {};
		\node [style=black dot] (218) at (-16.75, 5.25) {};
		\node [style=black dot] (219) at (-16.75, 3.75) {};
		\node [style=black dot] (220) at (-12.75, 6) {};
		\node [style=black dot] (221) at (-9.75, 6) {};
		\node [style=black dot] (222) at (-9.75, 3) {};
		\node [style=black dot] (223) at (-12.75, 3) {};
		\node [style=none] (225) at (-17.5, 4.5) {$W$};
		\node [style=none] (226) at (-11.25, 7.5) {$(W,b_W)$};
		\node [style=none] (230) at (-11.25, 2.25) {$r$};
		\node [style=none] (231) at (-25.5, 4.5) {$(L,b)$};
		\node [style=none] (232) at (-30, 2) {};
		\node [style=none] (233) at (-21, 2) {};
		\node [style=none] (234) at (-30, 3) {$A$};
		\node [style=none] (235) at (-29.25, 1.25) {$B$};
	\end{pgfonlayer}
	\begin{pgfonlayer}{edgelayer}
		\draw [style=c1] (172) to (174);
		\draw [style=c1] (171) to (177);
		\draw [style=c1] (180) to (172);
		\draw [style=c2] (169) to (175);
		\draw [style=c2] (171) to (178);
		\draw [style=c1] (174) to (169);
		\draw [style=c3] (170) to (177);
		\draw [style=c0] (180) to (171);
		\draw [style=c1] (171) to (179);
		\draw [style=c1] (179) to (180);
		\draw [style=c0] (179) to (172);
		\draw [style=c0] (176) to (169);
		\draw [style=c0] (175) to (170);
		\draw [style=c2] (170) to (176);
		\draw [style=c2] (176) to (175);
		\draw [style=c0] (177) to (178);
		\draw [style=c0] (178) to (170);
		\draw [style=c2] (169) to (197);
		\draw [style=c2] (197) to (172);
		\draw [style=c3] (198) to (169);
		\draw [style=c3] (198) to (172);
		\draw [style=c0] (197) to (198);
		\draw [style=c0] (198) to (174);
		\draw [style=c0] (174) to (197);
		\draw [style=c0] (213) to (211);
		\draw [style=c0] (211) to (212);
		\draw [style=c0] (212) to (214);
		\draw [style=c0] (214) to (215);
		\draw [style=c0] (215) to (213);
		\draw [style=c0] (213) to (212);
		\draw [style=c0] (212) to (215);
		\draw [style=c0] (215) to (211);
		\draw [style=c0] (211) to (214);
		\draw [style=c0] (214) to (213);
		\draw [style=c0] (202) to (200);
		\draw [style=c0] (200) to (199);
		\draw [style=c0] (199) to (201);
		\draw [style=c0] (201) to (202);
		\draw [style=c0] (202) to (199);
		\draw [style=c0] (201) to (200);
		\draw [style=c0] (217) to (218);
		\draw [style=c0] (218) to (219);
		\draw [style=c0] (219) to (216);
		\draw [style=c0] (216) to (217);
		\draw [style=c0] (203) to (206);
		\draw [style=c0] (206) to (205);
		\draw [style=c0] (205) to (204);
		\draw [style=c0] (204) to (203);
		\draw [style=c0] (203) to (205);
		\draw [style=c0] (206) to (204);
		\draw [style=c0] (208) to (207);
		\draw [style=c0] (207) to (210);
		\draw [style=c0] (210) to (209);
		\draw [style=c0] (209) to (208);
		\draw [style=c0] (208) to (210);
		\draw [style=c0] (209) to (207);
		\draw [style=triple] (223) to (222);
		\draw [style=c2] (221) to (222);
		\draw [style=c1] (220) to (223);
		\draw [style=triple] (220) to (221);
		\draw [style=reallythickwhitedash] (223) to (222);
		\draw [style=d1] (232.center) to (171);
		\draw [style=d1, bend left=15] (171) to (170);
		\draw [style=d1] (170) to (233.center);
	\end{pgfonlayer}
\end{tikzpicture}}
\end{minipage}
\end{center}
\vspace{2cm}
\begin{center}
\begin{minipage}{0.9\textwidth}
    \small\textbf{Statement on AI.} No LLM tools were used in the development of ideas in this paper, or in its preparation.
\end{minipage}
\end{center}

\pagebreak\section{Introduction}
In 1959, Tutte proved that a matroid can be realised by a graph if and only if it does not have a minor isomorphic to $U_{2,4}$, $F_7$, $F_7^*$, $M^*(K_5)$ or $M^*(K_{3,3})$~\cite{tutte1959matroids}. This was swiftly followed by an efficient algorithm in the case where the input matroid is binary~\cite{tutte1960algorithm}. Since then, many extensions and analogues have been investigated \cite{carmesin2017embeddingIII,geelen2018quasi,maimani2023structure}.
In this paper, we investigate the following.
\begin{dfn}
    A \emph{constraint matroid} is a pair $(M,c)$ where $M$ is a graphic matroid and $c$ is a function assigning each element of the ground set of $M$ to a (possibly empty) set of colours. A \emph{realisation} of a constraint matroid $(M,c)$ is a graphic realisation $G$ of $M$ so that the colour classes defined by $c$ form connected sets of edges in $G$.
\end{dfn}

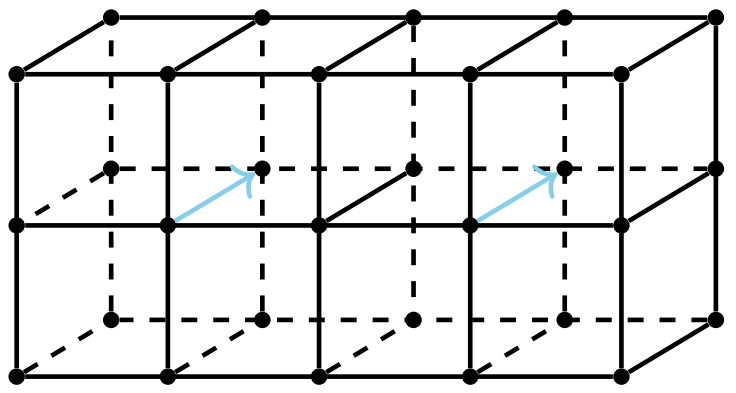
\begin{figure}
    \centering
    \begin{tikzpicture}
	\begin{pgfonlayer}{nodelayer}
		\node [style=black dot] (121) at (-4.75, -2) {};
		\node [style=black dot] (122) at (-4.75, 0) {};
		\node [style=black dot] (123) at (-4.75, 2) {};
		\node [style=black dot] (124) at (-2.75, 0) {};
		\node [style=black dot] (125) at (-2.75, -2) {};
		\node [style=black dot] (126) at (-2.75, 2) {};
		\node [style=black dot] (127) at (-0.75, 2) {};
		\node [style=black dot] (128) at (-0.75, -2) {};
		\node [style=black dot] (129) at (-0.75, 0) {};
		\node [style=black dot] (130) at (1.25, 0) {};
		\node [style=black dot] (131) at (3.25, 0) {};
		\node [style=black dot] (132) at (3.25, 2) {};
		\node [style=black dot] (133) at (1.25, 2) {};
		\node [style=black dot] (134) at (1.25, -2) {};
		\node [style=black dot] (135) at (3.25, -2) {};
		\node [style=black dot] (136) at (-3.5, -1.25) {};
		\node [style=black dot] (137) at (-3.5, 0.75) {};
		\node [style=black dot] (138) at (-3.5, 2.75) {};
		\node [style=black dot] (139) at (-1.5, 0.75) {};
		\node [style=black dot] (140) at (-1.5, -1.25) {};
		\node [style=black dot] (141) at (-1.5, 2.75) {};
		\node [style=black dot] (142) at (0.5, 2.75) {};
		\node [style=black dot] (143) at (0.5, -1.25) {};
		\node [style=black dot] (144) at (0.5, 0.75) {};
		\node [style=black dot] (145) at (2.5, 0.75) {};
		\node [style=black dot] (146) at (4.5, 0.75) {};
		\node [style=black dot] (147) at (4.5, 2.75) {};
		\node [style=black dot] (148) at (2.5, 2.75) {};
		\node [style=black dot] (149) at (2.5, -1.25) {};
		\node [style=black dot] (150) at (4.5, -1.25) {};
	\end{pgfonlayer}
	\begin{pgfonlayer}{edgelayer}
		\draw [style=c0] (123) to (121);
		\draw [style=c0] (121) to (135);
		\draw [style=c0] (135) to (132);
		\draw [style=c0] (132) to (123);
		\draw [style=c0] (122) to (131);
		\draw [style=c0] (133) to (134);
		\draw [style=c0] (128) to (127);
		\draw [style=c0] (126) to (125);
		\draw [style=c0] (138) to (147);
		\draw [style=c0] (147) to (150);
		\draw [style=d1] (150) to (136);
		\draw [style=d1] (136) to (138);
		\draw [style=d1] (140) to (141);
		\draw [style=d1] (142) to (143);
		\draw [style=d1] (149) to (148);
		\draw [style=d1] (137) to (146);
		\draw [style=c0] (123) to (138);
		\draw [style=c0] (126) to (141);
		\draw [style=c0] (127) to (142);
		\draw [style=c0] (133) to (148);
		\draw [style=c0] (147) to (132);
		\draw [style=c0] (146) to (131);
		\draw [style=c0] (144) to (129);
		\draw [style=d1] (137) to (122);
		\draw [style=d1] (121) to (136);
		\draw [style=d1] (125) to (140);
		\draw [style=d1] (128) to (143);
		\draw [style=d1] (134) to (149);
		\draw [style=c0] (135) to (150);
		\draw [style=c1 arrow] (124) to (139);
		\draw [style=c1 arrow] (130) to (145);
	\end{pgfonlayer}
\end{tikzpicture}
    \caption{Illustrated is the $2$-complex $C$, where each $4$-cycle bounds a face.}
    \label{fig:complex}
\end{figure}

Matroids realisable as graphs appear famously in Whitney's planarity criterion, and the problem of realising a constraint matroid arises naturally in the study of the embeddability of $2$-complexes in $\mathbb{R}^3$: for instance, consider the $2$-complex $C$ depicted in \autoref{fig:complex}. We have that~$C$ admits an embedding into $\mathbb{R}^3$. However, it is intuitive that the $2$-complex $C'$ obtained from~$C$ by identifying the blue edges, does not. It is not immediately clear how one might provide an efficient formal certificate for this fact. Using the machinery of~\cite{carmesin2017embeddingIV}, it turns out that to show that $C'$ is not embeddable, it suffices to show that a certain constraint matroid is not realisable.
Roughly speaking, one defines the dual matroid $M^*(C')$ of $C'$ to be the dual of the binary matroid whose circuits are the subsets of faces of $C'$ with null boundary (in the homological sense). An embedding for $C'$ induces a dual graph, which is a graphic realisation of $M^*(C')$. Conversely, a graphic realisation $G$ of $M^*(C')$ corresponds to an embedding of $C'$, as long as it satisfies some additional `dual graph connectivity constraints'. In the case of $C'$ in our example, these require that the edges of $G$ corresponding to the sets of faces of $C'$ incident with each vertex and each edge of $C'$ form connected subgraphs in $G$. Hence, we have translated our embedding problem into one about realising a constraint matroid.

In \cite{carmesin2017embeddingIII}, Carmesin gives a characterisation for the realisability of $3$-connected constraint matroids\footnote{A constraint matroid is $3$-connected if its underlying matroid is $3$-connected.} in terms of finitely many forbidden `constraint minors'.
Given a constraint matroid $(M,c)$, a \emph{constraint minor} $(M',c')$ of $(M,c)$ is a constraint matroid obtained from $(M,c)$ by successive application of the following operations:\vspace{0.1cm}

\begin{tabular}{l p{10cm}}
    (edge contraction) & contracting an element of $M$; \\
    (colour forgetting) & forgetting a colour from the colouring $c$; \\
    (parallel merging) & identifying a pair of parallel elements into a single element which inherits the colours of both original elements.
\end{tabular}

\begin{figure}
\begin{minipage}[t]{0.47\textwidth}
    \begin{tikzpicture}
	\begin{pgfonlayer}{nodelayer}
		\node [style=black dot] (284) at (-7.25, 1.25) {};
		\node [style=black dot] (285) at (-4.75, 1.25) {};
		\node [style=black dot] (286) at (-7.25, -1.25) {};
		\node [style=black dot] (287) at (-4.75, -1.25) {};
		\node [style=black dot] (288) at (-3.5, 1.5) {};
		\node [style=black dot] (289) at (-3.5, -1.5) {};
		\node [style=black dot] (290) at (-0.5, 1.5) {};
		\node [style=black dot] (291) at (-0.5, -1.5) {};
		\node [style=black dot] (292) at (-2, 0) {};
		\node [style=black dot] (293) at (0.5, 3.5) {};
		\node [style=black dot] (294) at (0.5, 0.5) {};
		\node [style=black dot] (295) at (4, 3.5) {};
		\node [style=black dot] (296) at (4, 0.5) {};
		\node [style=black dot] (297) at (1.5, 2) {};
		\node [style=black dot] (298) at (3, 2) {};
		\node [style=black dot] (299) at (5, 3.5) {};
		\node [style=black dot] (300) at (5, 0.5) {};
		\node [style=black dot] (301) at (8.5, 3.5) {};
		\node [style=black dot] (302) at (8.5, 0.5) {};
		\node [style=black dot] (303) at (6, 2) {};
		\node [style=black dot] (304) at (7.5, 2) {};
		\node [style=black dot] (305) at (0.5, -0.5) {};
		\node [style=black dot] (306) at (0.5, -3.5) {};
		\node [style=black dot] (307) at (4, -0.5) {};
		\node [style=black dot] (308) at (4, -3.5) {};
		\node [style=black dot] (309) at (1.5, -2) {};
		\node [style=black dot] (310) at (3, -2) {};
		\node [style=black dot] (311) at (5, -0.5) {};
		\node [style=black dot] (312) at (5, -3.5) {};
		\node [style=black dot] (313) at (8.5, -0.5) {};
		\node [style=black dot] (314) at (8.5, -3.5) {};
		\node [style=black dot] (315) at (6, -2) {};
		\node [style=black dot] (316) at (7.5, -2) {};
	\end{pgfonlayer}
	\begin{pgfonlayer}{edgelayer}
		\draw [style=c0] (284) to (285);
		\draw [style=c0] (284) to (287);
		\draw [style=c0] (287) to (286);
		\draw [style=c0] (286) to (285);
		\draw [style=c0] (292) to (289);
		\draw [style=c0] (288) to (292);
		\draw [style=c0] (292) to (290);
		\draw [style=c0] (290) to (288);
		\draw [style=c0] (292) to (291);
		\draw [style=c0] (291) to (289);
		\draw [style=c1 bold] (284) to (286);
		\draw [style=c1 bold] (285) to (287);
		\draw [style=c1 bold] (288) to (289);
		\draw [style=c1 bold] (290) to (291);
		\draw [style=c0] (295) to (293);
		\draw [style=c0] (296) to (294);
		\draw [style=c1 bold] (293) to (294);
		\draw [style=c1 bold] (295) to (296);
		\draw [style=c1 bold] (293) to (297);
		\draw [style=c1 bold] (297) to (294);
		\draw [style=c1 bold] (298) to (295);
		\draw [style=c1 bold] (298) to (296);
		\draw [style=c0] (297) to (298);
		\draw [style=c0] (301) to (299);
		\draw [style=c0] (302) to (300);
		\draw [style=c1 bold] (299) to (300);
		\draw [style=c1 bold] (301) to (302);
		\draw [style=c1 bold] (299) to (303);
		\draw [style=c1 bold] (303) to (300);
		\draw [style=c1 bold] (304) to (301);
		\draw [style=c0] (303) to (304);
		\draw [style=c0] (307) to (305);
		\draw [style=c0] (308) to (306);
		\draw [style=c1 bold] (305) to (309);
		\draw [style=c1 bold] (309) to (306);
		\draw [style=c1 bold] (310) to (307);
		\draw [style=c1 bold] (310) to (308);
		\draw [style=c0] (309) to (310);
		\draw [style=c0] (313) to (311);
		\draw [style=c0] (314) to (312);
		\draw [style=c1 bold] (313) to (314);
		\draw [style=c1 bold] (311) to (315);
		\draw [style=c1 bold] (315) to (312);
		\draw [style=c1 bold] (316) to (313);
		\draw [style=c0] (315) to (316);
		\draw [style=c0] (304) to (302);
		\draw [style=c0] (305) to (306);
		\draw [style=c0] (307) to (308);
		\draw [style=c0] (311) to (312);
		\draw [style=c0] (316) to (314);
	\end{pgfonlayer}
\end{tikzpicture}
    \caption{The six $3$-connected obstructions. Black edges should be taken to mean having no associated colours. In this figure, and following figures, matroids will routinely be depicted as graphs. One should consider the cycle matroid obtained from the graph depicted.}
    \label{fig:obst3con}
\end{minipage}
\hspace{\fill}
\begin{minipage}[t]{0.47\textwidth}
    \begin{tikzpicture}
	\begin{pgfonlayer}{nodelayer}
		\node [style=black dot] (317) at (-6.25, -2) {};
		\node [style=black dot] (318) at (-3.75, -2) {};
		\node [style=black dot] (319) at (-7.25, 0) {};
		\node [style=black dot] (320) at (-6.25, 2) {};
		\node [style=black dot] (321) at (-3.75, 2) {};
		\node [style=black dot] (322) at (-2.75, 0) {};
		\node [style=black dot] (323) at (-0.5, -2) {};
		\node [style=black dot] (324) at (2, -2) {};
		\node [style=black dot] (325) at (-1.5, 0) {};
		\node [style=black dot] (326) at (-0.5, 2) {};
		\node [style=black dot] (327) at (2, 2) {};
		\node [style=black dot] (328) at (3, 0) {};
		\node [style=black dot] (329) at (5.25, -2) {};
		\node [style=black dot] (330) at (7.75, -2) {};
		\node [style=black dot] (331) at (4.25, 0) {};
		\node [style=black dot] (332) at (5.25, 2) {};
		\node [style=black dot] (333) at (7.75, 2) {};
		\node [style=black dot] (334) at (8.75, 0) {};
		\node [style=none] (335) at (-5, -3.25) {$\text{Hex}$};
		\node [style=none] (336) at (1, -3.25) {$\text{Hex}^+$};
		\node [style=none] (337) at (6.5, -3.25) {$\text{Hex}^{++}$};
	\end{pgfonlayer}
	\begin{pgfonlayer}{edgelayer}
		\draw [style=c0] (319) to (322);
		\draw [style=c1] (319) to (320);
		\draw [style=c2] (320) to (321);
		\draw [style=c2] (317) to (318);
		\draw [style=c3] (321) to (322);
		\draw [style=c3] (318) to (322);
		\draw [style=c1] (319) to (317);
		\draw [style=c1] (325) to (326);
		\draw [style=c2] (326) to (327);
		\draw [style=c2] (323) to (324);
		\draw [style=c3] (327) to (328);
		\draw [style=c3] (324) to (328);
		\draw [style=c1] (325) to (323);
		\draw [style=c1] (331) to (332);
		\draw [style=c2] (332) to (333);
		\draw [style=c2] (329) to (330);
		\draw [style=c3] (333) to (334);
		\draw [style=c3] (330) to (334);
		\draw [style=c1] (331) to (329);
		\draw [style=c1] (325) to (328);
		\draw [style=sideby] (331) to (334);
	\end{pgfonlayer}
\end{tikzpicture}
    \caption{The \emph{unbalanced hexagon obstructions} $\text{Hex},\text{Hex}^+$ and $\text{Hex}^{++}$. Each contains three colours blue, red and green and the chord is (1) colourless; (2) just blue; and (3) both blue and red, respectively.}
    \label{fig:unbalanced_hexagon}
\end{minipage}
\end{figure}

\begin{thm}[Carmesin 2017 \cite{carmesin2017embeddingIII}]
    \label{thm:carmIII}
    A $3$-connected constraint matroid $(M,c)$ is realisable if and only if it does not contain a constraint minor\footnote{Carmesin's list of operations is different to the one in this paper. Indeed, we have replaced the operation of deleting a colourless edge by the operations of parallel merging and colour forgetting. In \autoref{appendix:translation} we justify this change.} isomorphic to one of the constraint matroids depicted in \autoref{fig:obst3con}.
\end{thm}
The first main result of this paper is the following extension of \autoref{thm:carmIII} to all connected constraint matroids. The condition `colour-complete' (see \autoref{dfn:colour-complete}) is a technical condition on the colouring that arises from the behaviour of $2$-separations in the matroid. In \autoref{sec:colour-complete} we will see an example which necessitates this assumption.
\begin{thm}\label{thm:forb}
    Let $(M,c)$ be a connected constraint matroid which is colour-complete. Then $(M,c)$ is realisable if and only if it does not contain as a constraint minor:\vspace{-.2cm}
    \begin{enumerate}
        \item one of the $3$-connected obstructions depicted in \autoref{fig:obst3con}; or\vspace{-.25cm}
        \item an unrealisable constraint diatom (definined \autoref{sec:diatoms}); or\vspace{-.25cm}
        \item an unrealisable constraint radiolarian (defined \autoref{sec:radiolarians}); or\vspace{-.25cm}
        \item one of the unbalanced hexagon obstructions, depicted in \autoref{fig:unbalanced_hexagon}.\vspace{-.2cm}
    \end{enumerate}
\end{thm}
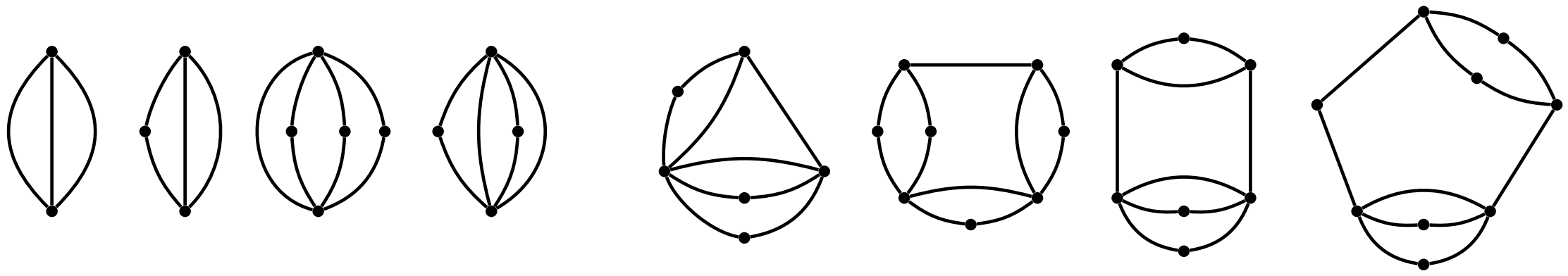
\begin{figure}
    \centering
    \begin{tikzpicture}
	\begin{pgfonlayer}{nodelayer}
		\node [style=black dot] (336) at (-10, 1.5) {};
		\node [style=black dot] (337) at (-10, -1.5) {};
		\node [style=black dot] (338) at (-7.5, 1.5) {};
		\node [style=black dot] (339) at (-7.5, -1.5) {};
		\node [style=black dot] (340) at (-8.25, 0) {};
		\node [style=black dot] (341) at (-5, 1.5) {};
		\node [style=black dot] (342) at (-5, -1.5) {};
		\node [style=black dot] (343) at (-5.5, 0) {};
		\node [style=black dot] (344) at (-4.5, 0) {};
		\node [style=black dot] (345) at (-3.75, 0) {};
		\node [style=black dot] (346) at (-1.75, 1.5) {};
		\node [style=black dot] (347) at (-1.75, -1.5) {};
		\node [style=black dot] (348) at (-2.75, 0) {};
		\node [style=black dot] (349) at (-1.25, 0) {};
		\node [style=black dot] (350) at (1.5, -0.75) {};
		\node [style=black dot] (351) at (4.5, -0.75) {};
		\node [style=black dot] (352) at (3, 1.5) {};
		\node [style=black dot] (353) at (6, -1.25) {};
		\node [style=black dot] (354) at (6, 1.25) {};
		\node [style=black dot] (355) at (8.5, 1.25) {};
		\node [style=black dot] (356) at (8.5, -1.25) {};
		\node [style=black dot] (357) at (10, 1.25) {};
		\node [style=black dot] (358) at (12.5, 1.25) {};
		\node [style=black dot] (359) at (10, -1.25) {};
		\node [style=black dot] (360) at (12.5, -1.25) {};
		\node [style=black dot] (361) at (13.75, 0.5) {};
		\node [style=black dot] (362) at (15.75, 2.25) {};
		\node [style=black dot] (363) at (18.25, 0.5) {};
		\node [style=black dot] (364) at (14.5, -1.5) {};
		\node [style=black dot] (365) at (17, -1.5) {};
		\node [style=black dot] (366) at (1.75, 0.75) {};
		\node [style=black dot] (367) at (3, -1.25) {};
		\node [style=black dot] (368) at (6.5, 0) {};
		\node [style=black dot] (369) at (7.25, -1.75) {};
		\node [style=black dot] (370) at (9, 0) {};
		\node [style=black dot] (371) at (11.25, 1.75) {};
		\node [style=black dot] (372) at (11.25, -1.5) {};
		\node [style=black dot] (373) at (16.75, 1) {};
		\node [style=black dot] (374) at (17.25, 1.75) {};
		\node [style=black dot] (375) at (3, -2) {};
		\node [style=black dot] (376) at (11.25, -2.25) {};
		\node [style=black dot] (377) at (15.75, -1.75) {};
		\node [style=black dot] (378) at (15.75, -2.5) {};
		\node [style=black dot] (379) at (5.5, 0) {};
	\end{pgfonlayer}
	\begin{pgfonlayer}{edgelayer}
		\draw [style=c0, bend right=45, looseness=1.25] (336) to (337);
		\draw [style=c0, bend right=15] (340) to (339);
		\draw [style=c0, bend right=45] (339) to (338);
		\draw [style=c0, bend right=15, looseness=0.75] (338) to (340);
		\draw [style=c0, bend left=15, looseness=0.75] (343) to (341);
		\draw [style=c0, bend left=15] (341) to (344);
		\draw [style=c0, bend left=15] (344) to (342);
		\draw [style=c0, bend left=15] (342) to (343);
		\draw [style=c0, bend left] (341) to (345);
		\draw [style=c0, bend left] (345) to (342);
		\draw [style=c0, bend right=15] (346) to (348);
		\draw [style=c0, bend left=15] (346) to (349);
		\draw [style=c0, bend left=15] (349) to (347);
		\draw [style=c0, bend left=15] (347) to (348);
		\draw [style=c0] (336) to (337);
		\draw [style=c0, bend left=45, looseness=1.25] (336) to (337);
		\draw [style=c0] (338) to (339);
		\draw [style=c0, bend right=75, looseness=1.25] (341) to (342);
		\draw [style=c0, bend right=15] (346) to (347);
		\draw [style=c0, bend left=60, looseness=1.25] (346) to (347);
		\draw [style=c0, bend right=15] (350) to (352);
		\draw [style=c0] (352) to (351);
		\draw [style=c0, bend left=15] (350) to (351);
		\draw [style=c0, bend left=15] (353) to (356);
		\draw [style=c0, bend left] (356) to (355);
		\draw [style=c0] (354) to (355);
		\draw [style=c0] (357) to (359);
		\draw [style=c0, bend left] (359) to (360);
		\draw [style=c0] (360) to (358);
		\draw [style=c0, bend left] (358) to (357);
		\draw [style=c0] (361) to (362);
		\draw [style=c0] (361) to (364);
		\draw [style=c0, bend left] (364) to (365);
		\draw [style=c0] (365) to (363);
		\draw [style=c0, bend right=15, looseness=0.75] (366) to (350);
		\draw [style=c0, bend left=15] (366) to (352);
		\draw [style=c0, bend right=15] (350) to (367);
		\draw [style=c0, bend right=15] (367) to (351);
		\draw [style=c0, bend left] (351) to (375);
		\draw [style=c0, bend left, looseness=0.75] (375) to (350);
		\draw [style=c0, bend right=15] (368) to (354);
		\draw [style=c0, bend left=15] (353) to (379);
		\draw [style=c0, bend left=15] (368) to (353);
		\draw [style=c0, bend left=15] (379) to (354);
		\draw [style=c0, bend right=15] (353) to (369);
		\draw [style=c0, bend right=15] (369) to (356);
		\draw [style=c0, bend left=15] (370) to (356);
		\draw [style=c0, bend left=15] (355) to (370);
		\draw [style=c0, bend left=15] (357) to (371);
		\draw [style=c0, bend left=15] (371) to (358);
		\draw [style=c0, bend right=15] (359) to (372);
		\draw [style=c0, bend right=15] (372) to (360);
		\draw [style=c0, bend left] (360) to (376);
		\draw [style=c0, bend left] (376) to (359);
		\draw [style=c0, bend right=15] (362) to (373);
		\draw [style=c0, bend right=15] (373) to (363);
		\draw [style=c0, bend right=15] (363) to (374);
		\draw [style=c0, bend right=15] (374) to (362);
		\draw [style=c0, bend right=15] (364) to (377);
		\draw [style=c0, bend right=15] (377) to (365);
		\draw [style=c0, bend left] (365) to (378);
		\draw [style=c0, bend left] (378) to (364);
	\end{pgfonlayer}
\end{tikzpicture}
    \caption{An example of diatoms (left) and radiolarians (right). The words `diatom' and `radiolarian' refer to types of phytoplankton and zooplankton, respectively, which both build a hard geometric shell out of glass.}
    \label{fig:sea_critters}
\end{figure}
A `diatom' is a matroid obtained from a cocircuit (bond) (that is, a collection of pairwise parallel edges) by subdividing each edge at most once. A `radiolarian' is a matroid obtained from a circuit by replacing some of its edges by diatoms. Each of these infinite classes does not admit particularly complex behaviour: for instance, checking the realisability of constraint diatoms and constraint radiolarians is relatively straightforward. See Sections~\ref{sec:diatoms} and \ref{sec:radiolarians} and see \autoref{fig:sea_critters} for an example.

In \autoref{appendix:evil_class}, we discuss why it seems unlikely that a neat forbidden-substructure characterisation exists for the realiability of constraint matroids in general. Even so, computing the realisability of constraint matroids in general remains tractable: the second main result of this paper is the following, where we consider the `size' of a constraint matroid to be the sum of its number of elements and number of colours. Importantly, in this theorem we do not make any assumptions about the connectivity of the matroid.
\begin{thm}\label{cor:alg}
    Given an independence oracle, the realisability of a constraint matroid can be decided in polynomial time.
\end{thm}
To prove our two main results in the connected case, we will first build an analogue of the tree decomposition of Cunningham and Edmonds~\cite{cunningham1980combinatorial} for our setting, which splits a connected matroid in a tree-like way along its totally-nested $2$-separations. This will be \autoref{thm:matroids_tutte_analogue}. The whole constraint matroid will then be realisable if and only if each of its pieces are.

\begin{figure}
    \centering
    \begin{tikzpicture}
	\begin{pgfonlayer}{nodelayer}
		\node [style=black dot] (169) at (-7.75, 1.25) {};
		\node [style=black dot] (170) at (-7.75, -1.25) {};
		\node [style=black dot] (171) at (-10.25, -1.25) {};
		\node [style=black dot] (172) at (-10.25, 1.25) {};
		\node [style=black dot] (174) at (-7.5, 3) {};
		\node [style=black dot] (175) at (-6, 1) {};
		\node [style=black dot] (176) at (-6, -1) {};
		\node [style=black dot] (177) at (-8, -3) {};
		\node [style=black dot] (178) at (-10, -3) {};
		\node [style=black dot] (179) at (-12, -1) {};
		\node [style=black dot] (180) at (-12, 1) {};
		\node [style=black dot] (197) at (-10.5, 3) {};
		\node [style=black dot] (198) at (-9, 4) {};
		\node [style=none] (231) at (-9, 0) {$(L,b)$};
		\node [style=black dot] (232) at (7.75, -1.5) {};
		\node [style=black dot] (233) at (10.25, -1.5) {};
		\node [style=black dot] (234) at (7.75, 1.5) {};
		\node [style=black dot] (235) at (10.25, 1.5) {};
		\node [style=black dot] (236) at (5.75, 1.25) {};
		\node [style=black dot] (237) at (12.25, 1.25) {};
		\node [style=black dot] (238) at (5.75, -1.25) {};
		\node [style=black dot] (239) at (12.25, -1.25) {};
		\node [style=none] (240) at (9, 0) {$(N,d)$};
		\node [style=black dot] (241) at (-1.25, 2) {};
		\node [style=black dot] (242) at (-1.25, -2) {};
		\node [style=black dot] (243) at (1.25, 1.25) {};
		\node [style=black dot] (244) at (1.25, -1.25) {};
		\node [style=black dot] (245) at (-2.5, 0) {};
		\node [style=black dot] (246) at (0.5, 3) {};
		\node [style=black dot] (247) at (0.5, -2.75) {};
		\node [style=black dot] (248) at (2.5, 0) {};
		\node [style=none] (249) at (2.25, 2.25) {$(M,c)$};
	\end{pgfonlayer}
	\begin{pgfonlayer}{edgelayer}
		\draw [style=c1] (172) to (174);
		\draw [style=c1] (171) to (177);
		\draw [style=c1] (180) to (172);
		\draw [style=c2] (169) to (175);
		\draw [style=c2] (171) to (178);
		\draw [style=c1] (174) to (169);
		\draw [style=c3] (170) to (177);
		\draw [style=c0] (180) to (171);
		\draw [style=c1] (171) to (179);
		\draw [style=c1] (179) to (180);
		\draw [style=c0] (179) to (172);
		\draw [style=c0] (176) to (169);
		\draw [style=c0] (175) to (170);
		\draw [style=c2] (170) to (176);
		\draw [style=c2] (176) to (175);
		\draw [style=c0] (177) to (178);
		\draw [style=c0] (178) to (170);
		\draw [style=c2] (169) to (197);
		\draw [style=c2] (197) to (172);
		\draw [style=c3] (198) to (169);
		\draw [style=c3] (198) to (172);
		\draw [style=c0] (197) to (198);
		\draw [style=c0] (198) to (174);
		\draw [style=c0] (174) to (197);
		\draw [style=c2] (234) to (235);
		\draw [style=c1] (233) to (232);
		\draw [style=c2] (234) to (236);
		\draw [style=c0] (236) to (238);
		\draw [style=c1] (238) to (232);
		\draw [style=c0] (239) to (233);
		\draw [style=sideby] (235) to (237);
		\draw [style=c0] (237) to (239);
		\draw [style=c0] (236) to (232);
		\draw [style=c0] (238) to (234);
		\draw [style=c0] (235) to (239);
		\draw [style=c0] (233) to (237);
		\draw [style=c0] (246) to (241);
		\draw [style=c0] (242) to (247);
		\draw [style=c0] (243) to (242);
		\draw [style=c0] (241) to (244);
		\draw [style=c0] (244) to (245);
		\draw [style=c0] (245) to (243);
		\draw [style=c0] (241) to (245);
		\draw [style=c0] (245) to (242);
		\draw [style=c1] (246) to (243);
		\draw [style=c1] (247) to (244);
		\draw [style=c1] (243) to (248);
		\draw [style=c0] (248) to (244);
		\draw [style=c0] (242) to (241);
	\end{pgfonlayer}
\end{tikzpicture}
    \caption{So called `unbalanced', `unstable' and `disagreeable' connected constraint matroids $(L,b)$, $(M,c)$ and $(N,d)$, respectively.}
    \label{fig:obs1}
\end{figure}
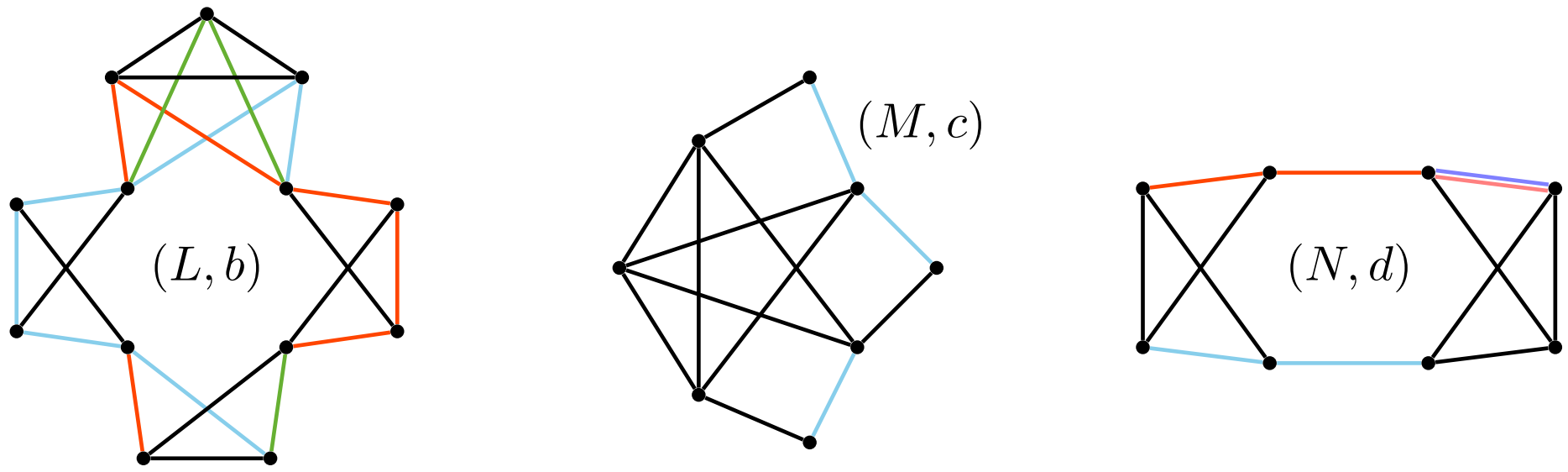

However, in proving \autoref{thm:matroids_tutte_analogue}, we will have to make three further structural assumptions about our connected constraint matroid. In particular, we will need to forbid three types of `non-local' obstruction: We will call these obstructions to decomposability `unbalanced', `unstable' and `diagreeable' and they will be formally defined in Subsections \ref{subsec:unbalanced}, \ref{subsec:unstable} and \ref{subsec:disagreeable}, respectively. In \autoref{fig:obs1}, we depict a prototypical example of each of these types of obstruction.

After proving \autoref{thm:matroids_tutte_analogue}, there are two remaining steps to proving \autoref{thm:forb}, our forbidden substructure result. That is, proving that each of the `non-local' obstructions results in one of our forbidden constraint minors (this will be \autoref{sec:first_steps}) and proving that finding an unrealisable `piece' of the decomposition from \autoref{thm:matroids_tutte_analogue} also results in one of our forbidden constraint minors. This last step will be stratified in terms of the type of piece and split across Sections \ref{sec:circuits_detailed}, \ref{sec:diatoms_detailed} and \ref{sec:done}, completing our proof.

Finally, our algorithmic result will follow from \autoref{thm:new_red_to_con}, which takes an arbitrary constraint matroid and decomposes it into a collection of connected constraint matroids, upon which we can apply our previous theory. This will occur in \autoref{sec:disconnected}. In this section, we will require amongst other things, an analogue of Whiteny's $2$-isomorphism Theorem (\autoref{thm:2-iso}, see~\cite{whitney19332}), which states that the set of graphs that have the same cycle matroid are all $2$-isomorphic. Loosely speaking, that means that they are identical up to `reconfiguring' their maximal $2$-connected pieces and `twisting' maximal $2$-connected components about $2$-separators. Our analogue will be~\autoref{lem:homotwisty}. In \autoref{sec:algorithmics}, we will tie together everything in the algorithmic analysis that completes the proof of \autoref{cor:alg}.

For background on graph theory, we refer the reader to Diestel's textbook \cite{diestel2006}. In this paper, graphs are understood to be undirected multigraphs, and we may treat two graphs as equal if there is a isomorphism between them which preserves edge-labels. For background on matroids, see the textbook of Oxley \cite{oxley2006}. Since almost all of the matroids in this paper are graphic, we will refer to the elements of their ground sets as \emph{edges}. Some particularly pertinent preliminary definitions and results will be included in \autoref{sec:recall}, including the definition of `colour-complete' in \autoref{sec:colour-complete}.

The beginning of this paper is structured as follows: in Sections \ref{sec:circuits}, \ref{sec:diatoms} and \ref{sec:radiolarians} we introduce our two infinite classes of obstruction, the constraint diatoms and constraint radiolarians (of which the `constraint circuits' of \autoref{sec:circuits} are a special case).

\subsection*{Acknowledgments}
The author would like to thank Johannes Carmesin and Tim Planken for their helpful suggestions, as well as the of the Discrete Maths seminar at TU Freiberg who engaged with the development of this work over the past year.

\section{Circuits}\label{sec:circuits}
We refer to a matroid that consists of a single circuit simply as \emph{a circuit}.

\begin{dfn}[Overlay]
    Let $(C,c)$ be a constraint circuit. The \emph{overlay} $\lozenge (C,c)$ of $(C,c)$ is the smallest binary matroid containing $C$ and, for each colour class $d$ of $(C,c)$, an additional element~$e_d$ and an additional circuit $\{e_d\}\cup d$. That is, the circuit space of $\lozenge (C,c)$ is generated (via symmetric difference) from the circuit $C$ and the circuits of the form $\{e_d\}\cup d$.
\end{dfn}
The matroid overlay construction is a specialisation the construction of the same name from~\cite{matroid_overlay}, where a richer theory is developed.

\begin{lem}\label{overlay}
    The realisations of a constraint circuit $(C,c)$ are exactly the restrictions to $E(C)$ of the realisations of its overlay $\lozenge(C,c)$.
\end{lem}
\begin{pf}
    Suppose that $O$ is a realisation of $(C,c)$. Let $O'$ be the graph obtained from $O$ by adding an edge $e_d$ joining the endvertices of the interval in $O$ coloured $d$ for each colour $d$. Then $O'$ is a graphic realisation of $\lozenge (C,c)$. Now let a graphic realisation $G$ of $\lozenge (C,c)$ be given. The graph $G$ contains $C$ as a cycle. Furthermore, for each colour $d$, there is a cycle $C_d$ of $G$ containing exactly the $d$-coloured edges of $C$ and one additional edge $e_d$ which is not in $C$. We conclude that the $d$-coloured edges form an interval in the cyclic order of the edges of $C$. Hence the subgraph induced on $C$ in $G$ is a realisation for the constraint circuit $(C,c)$.
\end{pf}
Recognising whether a given constraint circuit is realisable is thus efficiently computable, for example by using the algorithm of Tutte for recognising graphicness of binary matroids \cite{tutte1960algorithm}.

\section{Diatoms}\label{sec:diatoms}
Recall that a \emph{bond} (also known as a \emph{cocircuit}) is a matroid that consists of a single parallel class. A \emph{diatom} is a matroid obtained from a bond with at least three edges by subdividing each edge at most once.

Let $(D,c)$ be a constraint diatom. Denote its non-singleton series classes by $(\{a_i,b_i\})_{i\in[n]}$. The possible graphic realisations of the underlying matroid $D$ are given by a choice of order for each of these non-singleton series classes. For the purposes of the next few results, we will assume these non-singleton series classes have prespecified `default' order $(a_i,b_i)$.

We say that a constraint diatom is \emph{clean} if no colour colours a whole series class and the \emph{clean version} of any diatom is obtained by forgetting any colours which colour a whole series class.

\begin{obs}
    A constraint diatom is realisable if and only if its clean version is. \qed
\end{obs}
A \emph{signed graph} is a graph whose edges are labeled either `positive' or `negative'.
To characterise when clean constraint diatoms are realisable, we make the following definition.
\begin{dfn}
    Let $(D,c)$ be a clean constraint diatom. The \emph{frustration graph} of $(D,c)$ is the signed graph whose vertices are the non-singleton series classes of $(D,c)$ and where there is: (1) a positive edge between distinct $\{a_i,b_i\}$ and $\{a_j,b_j\}$ when there is a colour that colours either both $a_i$ and $a_j$ or both $b_i$ and $b_j$; and (2) a negative edge between distinct $\{a_i,b_i\}$ and $\{a_j,b_j\}$ when there is a colour that colours either both $a_i$ and $b_j$ or both $b_i$ and $a_j$.
\end{dfn}

A signed graph is \emph{balanced} if it does not have a cycle with an odd number of negative edges.
\begin{lem}\label{lem:quick_check_diatom}
    Let $(D,c)$ be a clean constraint diatom. Then $(D,c)$ is realisable if and only if its frustration graph is balanced. \qed
\end{lem}

\section{Radiolarians}\label{sec:radiolarians}
Let $A$ and $B$ be matroids and $a$ and $b$ edges in $A$ and $B$, respectively. The \emph{$2$-sum of $A$ and $B$ along $(a,b)$} is the matroid $A\oplus_{(a,b)}B$ with edge set $E(A)\cup E(B)\setminus \{a,b\}$, whose circuits are exactly the circuits of $A-a$, the circuits of $B-b$ and the circuits of the form $(C_A\setminus\{a\})\cup(C_B\setminus\{b\})$ where $C_A\ni a$ is a circuit of $A$ and $C_B\ni b$ is a circuit of $B$.

A \emph{radiolarian} $R$ is a connected matroid obtained from a circuit $C$ in the following way: choose a subset $E'\subseteq E(C)$; then for each $e\in E'$, let $D_e$ be a diatom with one of its series classes the singleton~$\{d_e\}$; then set $R$ to be the $2$-sum of $R$ with each $D_e$ along $(e,d_e)$.

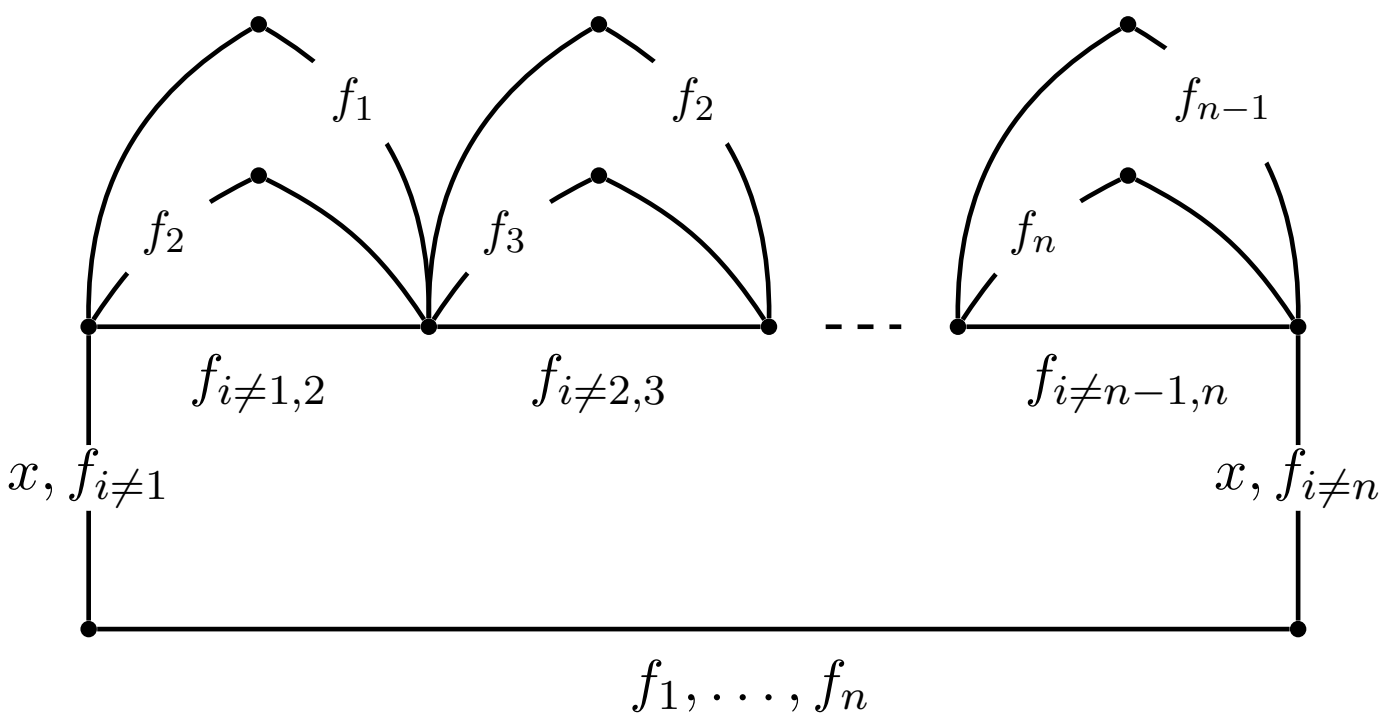
\begin{figure}
\begin{minipage}[t]{0.55\textwidth}
    \begin{tikzpicture}
	\begin{pgfonlayer}{nodelayer}
		\node [style=black dot] (380) at (-6, -2) {};
		\node [style=black dot] (381) at (-6, 2) {};
		\node [style=black dot] (382) at (-1.5, 2) {};
		\node [style=black dot] (383) at (3, 2) {};
		\node [style=black dot] (384) at (5.5, 2) {};
		\node [style=black dot] (385) at (10, 2) {};
		\node [style=none] (386) at (3.75, 2) {};
		\node [style=none] (387) at (4.75, 2) {};
		\node [style=black dot] (388) at (-3.75, 4) {};
		\node [style=black dot] (389) at (7.75, 4) {};
		\node [style=black dot] (390) at (7.75, 6) {};
		\node [style=black dot] (391) at (-3.75, 6) {};
		\node [style=black dot] (392) at (0.75, 4) {};
		\node [style=black dot] (393) at (0.75, 6) {};
		\node [style=black dot] (394) at (10, -2) {};
		\node [style=none] (397) at (2.75, -2.75) {$f_1,\ldots,f_n$};
		\node [style=none] (399) at (-3.75, 1.25) {$f_{i\neq1,2}$};
		\node [style=none] (400) at (0.75, 1.25) {$f_{i\neq2,3}$};
		\node [style=none] (401) at (7.75, 1.25) {$f_{i\neq n-1,n}$};
		\node [style=backing] (402) at (-5, 3.25) {$f_2$};
		\node [style=backing] (403) at (-2.5, 5) {$f_1$};
		\node [style=backing] (404) at (2, 5) {$f_2$};
		\node [style=backing] (405) at (-0.5, 3.25) {$f_3$};
		\node [style=backing] (406) at (9, 5) {$f_{n-1}$};
		\node [style=backing] (407) at (6.5, 3.25) {$f_n$};
		\node [style=medium back] (408) at (-6, 0) {};
		\node [style=none] (409) at (-6, 0) {$x,f_{i\neq1}$};
		\node [style=medium back] (410) at (10, 0) {};
		\node [style=none] (411) at (10, 0) {$x,f_{i\neq n}$};
	\end{pgfonlayer}
	\begin{pgfonlayer}{edgelayer}
		\draw [style=d1] (386.center) to (387.center);
		\draw [style=c0] (380) to (394);
		\draw [style=c0] (394) to (385);
		\draw [style=c0, bend right=15] (385) to (389);
		\draw [style=c0, bend right=15] (389) to (384);
		\draw [style=c0] (384) to (385);
		\draw [style=c0, bend left] (390) to (385);
		\draw [style=c0, bend right] (390) to (384);
		\draw [style=c0, bend right=15] (383) to (392);
		\draw [style=c0, bend right=15] (392) to (382);
		\draw [style=c0, bend left] (382) to (393);
		\draw [style=c0, bend left] (393) to (383);
		\draw [style=c0] (383) to (382);
		\draw [style=c0, bend right=15] (382) to (388);
		\draw [style=c0, bend left=15] (381) to (388);
		\draw [style=c0, bend right] (391) to (381);
		\draw [style=c0, bend left] (391) to (382);
		\draw [style=c0] (382) to (381);
		\draw [style=c0] (381) to (380);
	\end{pgfonlayer}
\end{tikzpicture}
\end{minipage}
\hspace{\fill}%
\begin{minipage}[t]{0.43\textwidth}
\vspace{-3cm}
    \caption{A constraint radiolarian $(R,c)$ with the colours $x,f_1,\ldots,f_n$. It is not realisable. However, any constraint matroid obtained from $R$ by contracting an edge or forgetting a colour is realisable. We do not need to consider parallel merging since $R$ has no parallel edges.}
    \label{fig:radio_gaga}
\end{minipage}
\end{figure}

\begin{eg}
    The constraint radiolarian depicted in \autoref{fig:radio_gaga} is unrealisable and minimal with this property under the constraint minor operations.
\end{eg}

Determining whether a constraint radiolarian $(R,c)$ is realisable requires a mixture of the machinery from \autoref{sec:circuits} and \autoref{sec:diatoms}, along with a few new ideas. The necessary tools will be developed in the course of the proof of \autoref{thm:matroids_tutte_analogue}. Intuitively, one replaces each diatom in $R$ by a series class of two edges, coloured appropriately, and then checks to see if the resulting constraint circuit is realisable using \autoref{overlay}.

\section{Preliminaries: matroids and connectivity}
\label{sec:recall}
In this section, we recall some matroid terminology and introduce `mixed-sum tree decompositions', which slightly generalise the tree-decomposition like the so-called `Tutte decomposition' (see below).

Let $M$ be a matroid. We say that $M$ is \emph{connected} if each pair of edges appear in some shared cycle of $M$. The maximal connected subsets in $E$ are called the \emph{components} of $M$ and are thought of as matroids themselves (by restriction). Let $(A,B)$ be a bipartition of $E$. We say that $(A,B)$ is a \emph{$2$-separation} if both $|A|\geq2$ and $|B|\geq2$ and $r(A)+r(B)=r(M)-1$. Let $(A,B)$ and $(C,D)$ be $2$-separations of $M$. We write $(A,B)\leq(C,D)$ if $A\subseteq C$ and $D\subseteq B$. We say that $(A,B)$ and $(C,D)$ are \emph{nested} with eachother if $(A,B)$ and $(C,D)$ are $\leq$-comparable, possibly after reordering one of the pairs. A $2$-separation $(A,B)$ is \emph{totally-nested} if it is nested with all other $2$-separators in $M$.

\begin{lem}[fokelore; \cite{oxley2006}]\label{concon}
    Let $G$ be a loopless graph without isolated vertices and that $|V(G)|\geq3$. Then the cycle matroid $M(G)$ is a connected matroid if and only if $G$ is a $2$-connected graph.
\end{lem}

With \autoref{concon} in mind, we define a $2$-separations for a connected graph as a $2$-separation of $M(G)$. Given a connected graph $G$ and a $2$-separation $(A,B)$ of $G$, there are exactly two vertices meeting edges from both $A$ and $B$, call this pair the \emph{separator} of $(A,B)$ in $G$.

Given a graph $G$, two edges are in \emph{series} if they are in series in the matroid $M(G)$. A connected set of pairwise-series edges in a graph forms a path whose internal vertices all have degree two.
\begin{dfn}[Mixed-sum for graphs]
    Let $G_1$ and $G_2$ be disjoint graphs. Let $\Ecal_1$ be a connected set of pairwise-series edges in $G_1$ and $\Ecal_2$ a connected set of pairwise-series edges in $G_2$. Then a \emph{mixed-sum of $G_1$ and $G_2$ along $(\Ecal_1,\Ecal_2)$} is a graph obtained by identifying the endvertices of the path formed by $\Ecal_1$ in $G_1$ with the path formed by $\Ecal_2$ in $G_2$ and then deleting the interior of the paths $\Ecal_1$ and $\Ecal_2$. There are two choices of how to make this identification and so two mixed-sums for each choice of data.
\end{dfn}

\begin{dfn}[Mixed-sum for matroids]
    Let $M_1$ and $M_2$ be disjoint graphs. Let $\Ecal_1$ be a set of pairwise-series edges in $G_1$ and $\Ecal_2$ a set of pairwise-series edges in $G_2$. Then the \emph{mixed-sum of $M_1$ and $M_2$ along $(\Ecal_1,\Ecal_2)$} is the matroid obtained after a choice of $e_1\in\Ecal_1$ and $e_2\in\Ecal_2$ as the $2$-sum of $M_1/(\Ecal_1-e_1)$ and $M_2/(\Ecal_2-e_2)$ along the edges $e_1$ and $e_2$. We denote the mixed-sum by $M_1\oplus_{(\Ecal_1,\Ecal_2)}M_2$.
\end{dfn}
The mixed sum has ground set $(E(M_1)\cup E(M_2))\setminus(\Ecal_1\cup\Ecal_2)$ and the following circuits:\vspace{-.25cm}
\begin{enumerate}
        \item for each circuit $o\in M_1$ disjoint from $\Ecal_1$, the set $o$ appears as a circuit in $M_1\oplus_{(\Ecal_1,\Ecal_2)}M_2$; and\vspace{-.25cm}
        \item for each circuit $o\in M_2$ disjoint from $\Ecal_2$, the set $o$ appears as a circuit in $M_1\oplus_{(\Ecal_1,\Ecal_2)}M_2$; and\vspace{-.25cm}
        \item for each pair of circuits $o_1$ from $M_1$ and $o_2$ from $M_2$ with $\Ecal_1\subseteq o_1$ and $\Ecal_2\subseteq o_2$, the set $(o_1\cup o_2)\setminus(\Ecal_1,\Ecal_2)$ appears as a circuit in $M_1\oplus_{(\Ecal_1,\Ecal_2)}M_2$.\vspace{-.25cm}
\end{enumerate}
Unlike mixed-sum for graphs, the mixed-sum of matroids is uniquely defined by the starting data.

\begin{lem}[fokelore]
    Let $M:=M_1\oplus_{(\{e_1\},\{e_2\})}M_2$ a matroid which decomposes as a $2$-sum. Then $M_1$ can be obtained from $M$ by contracting and parallel-merging edges from $E(M_2)\setminus\{e_2\}$, and then relabeling one edge from $E(M_2)\setminus\{e_2\}$ to $e_1$.
\end{lem}
\begin{pf}[Sketch of proof.]
    Let $D$ be a bond (co-circuit) in $M$ that uses edges from both $M_1$ and $M_2$. Let $D_2:=E(M_2)\cap D$. Let $M'$ be the matroid obtained from $M$ by contracting all of the edges in $E(M_2)$ but not in $D_2$. Then $M'$ is the matroid obtained from $M_1$ by replacing $e_1$ by a parallel class, whose edges are exactly those from $D_2$.
\end{pf}

\begin{lem}\label{lem:mixed-sum_graphic}
    Let $G_1$ and $G_2$ be disjoint graphs. Let $\Ecal_1$ be a connected set of pairwise-series edges in $G_1$ and $\Ecal_2$ a connected set of pairwise-series edges in $G_2$. Let $G$ be a mixed-sum of $G_1$ and $G_2$ along $(\Ecal_1,\Ecal_2)$. Then we have that $M(G)=M(G_1)\oplus_{(\Ecal_1,\Ecal_2)}M(G_2)$.
\end{lem}
\begin{pf}
    The proof is a simple exercise.
\end{pf}

\begin{dfn}[Mixed-sum tree decomposition for matroids]
    A \emph{mixed-sum tree decomposition} $(T,\{M_u\},\{\Ecal_u^e\})$ is a tree $T$ along with the following data: a set $M_u$ of matroids indexed by the nodes $u$ of $T$; and for each edge $e=uv$ in $T$ a pair $(\Ecal_u^e,\Ecal_v^e)$ of sets of series edges with $\Ecal_u^e\subseteq E(M_u)$ and $\Ecal_v^e\subseteq E(M_v)$. We require that the sets $\Ecal_u^e$ are all disjoint from eachother for a fixed node $u$.
    
    The matroid \emph{built} by $(T,\{M_u\},\{\Ecal_u^e\})$ is the matroid $M$ constructed via repeated mixed-sums of a matroid corresponding to a leaf of $M$ with its neighbour via the pair corresponding to its connecting edge. In this case, we say that $(T,\{M_u\},\{\Ecal_u^e\})$ is a \emph{mixed-sum tree decomposition of $M$}.
\end{dfn}

A connected matroid $M$ is said to be \emph{$3$-connected} if $|E(M)|\ge 4$ and $M$ does not decompose as a $2$-sum $M=A\oplus_{(\{a\},\{b\})} B$ where both $|E(A)|,|E(B)|\geq 3$ (equivalently, $M$ has no $2$-separation). The following, which is a direct consequence of \autoref{thm:2-iso}, we will use throughout, sometimes without direct reference.

\begin{lem}
    $3$-connected matroids have unique graphic realisations, up to relabeling vertices. \qed
\end{lem}

We say that a mixed-sum tree decomposition is a \emph{$2$-sum tree decomposition} when all the mixed sums involve only pairs of singletons. That is, they are standard $2$-sums.

\begin{thm}[Cunningham and Edmonds 1973 \cite{cunningham1980combinatorial}]\label{thm:matroids_tutte}
    Let $M$ be a connected matroid. Then there exists a unique $2$-sum tree decomposition of $M$ into circuits, bonds and $3$-connected matroids so that no pair of neighbouring nodes in the decomposition tree $T$ both correspond to circuits or both correspond to bonds. Furthermore, the separations induced by the edges of the tree $T$ are exactly the (unordered) totally-nested $2$-separations of $M$.
\end{thm}

Given a connected matroid $M$, we will refer to the tree decomposition obtained from \autoref{thm:matroids_tutte} as the \emph{Tutte decomposition} of $M$ and the individual matroids of this tree decomposition as the \emph{torsos} of $M$. The auxiliary tree $T$ is called the \emph{Tutte decomposition tree}.

\begin{cor}\label{cor:torso_graphic}
    A connected matroid $M$ is graphic if and only if each of its torsos are graphic.
\end{cor}
\begin{pf}
    Combine \autoref{lem:mixed-sum_graphic} and \autoref{thm:matroids_tutte}.
\end{pf}

Let $M$ be a connected matroids. The torsos of $M$ which are circuits will be called \emph{circular}. Let $e$ be an element of a torso $W$ of $M$. Then we say $e$ is a \emph{torso edge} if $e$ is not present in the matroid $M$. That is, it is `used up' in one of the $2$-sum operations between a pair of torsos.

\begin{dfn}[Branch at a torso edge]
Let $T$ be the Tutte decomposition tree for a connected matroid $M$, and $w$ one of its nodes corresponding to a torso $W$. Then each torso edge $e$ of $W$ corresponds to an edge $wu$ of $T$ for some node $u\in V(T)$. The edge $wu$ induces a totally-nested $2$-separation $(A_e,B_e)$ of $M$ where $B_e$ consists of the non-torso edges of the torsos corresponding to nodes in the connected component of $T-wu$ which does not contain $w$. We call $B_e$ the \emph{branch of $M$ at the torso edge $e$}.
\end{dfn}

\subsection{Colour-completion}\label{sec:colour-complete}
In our main result \autoref{thm:forb}, we assumed our constraint matroids were `colour-complete'. We can now define this property and explain why it is necessary. Let $(M,c)$ and $(M,c')$ be constraint matroids with the same underlying matroid. We say the colouring $c'$ is an \emph{extension} of $c$ if we can obtain $c$ from $c'$ by forgetting colours.
\begin{dfn}[Colour-complete]\label{dfn:colour-complete}
    Let $(M,c)$ be a connected constraint matroid. For a colour $f$ and a totally-nested $2$-separation $(A,B)$ of $M$, we say that \emph{$f$ completes $(A,B)$} if $f$ colours none of the edges in $A$ and all of the edges in $B$. We say that $(M,c)$ is \emph{colour-complete} if all of its totally-nested $2$-separations are completed by some colour. A \emph{colour-completion} of $(M,c)$ is a colour-complete constraint matroid $(M,c')$ where $c'$ is an extension of $c$ by adding colours which complete totally-nested $2$-separations of $M$.
\end{dfn}

\begin{obs}
    Let $(M,c)$ be a connected constraint matroid. Then a colour-completion of $(M,c)$ always exsits and uses at most $|E(M)|$ more colours than $c$.
\end{obs}
\begin{pf}
    It is fokeloric that $M$ has at most $|E(M)|$ totally-nested $2$-separations.
\end{pf}

\begin{obs}
    Let $(M,c)$ be a connected constraint matroid and $(M,c')$ a colour-completion thereof. Then $(M,c)$ is realisable if and only if $(M,c')$ is. \qed
\end{obs}

In the following example, we explore an infinite antichain with respect to the constraint minor relation, which shows that colour-complete is a necessary condition in \autoref{thm:forb}.
\begin{figure}
    \centering
    \begin{tikzpicture}
	\begin{pgfonlayer}{nodelayer}
		\node [style=none] (265) at (-8, 2) {$(M,c)$};
		\node [style=black dot] (266) at (-14, 2.5) {};
		\node [style=black dot] (267) at (-14, 0) {};
		\node [style=black dot] (268) at (-10, 0) {};
		\node [style=black dot] (269) at (-6, 0) {};
		\node [style=black dot] (270) at (-2, 0) {};
		\node [style=black dot] (271) at (-2, 2.5) {};
		\node [style=black dot] (272) at (-12.75, -3) {};
		\node [style=black dot] (273) at (-11.25, -3) {};
		\node [style=black dot] (274) at (-8.75, -3) {};
		\node [style=black dot] (275) at (-7.25, -3) {};
		\node [style=black dot] (276) at (-4.75, -3) {};
		\node [style=black dot] (277) at (-3.25, -3) {};
		\node [style=none] (278) at (-14.5, 0) {};
		\node [style=none] (279) at (-9.5, 0) {};
		\node [style=none] (280) at (-10.5, 0) {};
		\node [style=none] (281) at (-6.5, 0) {};
		\node [style=none] (282) at (-5.5, 0) {};
		\node [style=none] (283) at (-1.5, 0) {};
		\node [style=none] (284) at (-0.5, -0.25) {};
		\node [style=none] (285) at (-0.5, 0) {};
		\node [style=none] (286) at (-0.5, 0.25) {};
		\node [style=none] (287) at (0.5, 0) {};
		\node [style=none] (288) at (0.25, 0.25) {};
		\node [style=none] (289) at (0.25, -0.25) {};
		\node [style=none] (290) at (9.5, 0.5) {$(C,c'')$};
		\node [style=black dot] (291) at (2, 1.5) {};
		\node [style=black dot] (292) at (2, -1.5) {};
		\node [style=black dot] (293) at (7, -1.5) {};
		\node [style=black dot] (294) at (12, -1.5) {};
		\node [style=black dot] (295) at (17, -1.5) {};
		\node [style=black dot] (296) at (17, 1.5) {};
		\node [style=black dot] (297) at (4.5, -1.5) {};
		\node [style=black dot] (298) at (9.5, -1.5) {};
		\node [style=black dot] (299) at (14.5, -1.5) {};
		\node [style=backing] (302) at (-13.5, -2) {$f_1$};
		\node [style=backing] (303) at (-10.5, -2) {$f_2$};
		\node [style=backing] (304) at (-9.5, -2) {$f_2$};
		\node [style=backing] (305) at (-6.5, -2) {$f_3$};
		\node [style=backing] (306) at (-5.5, -2) {$f_3$};
		\node [style=backing] (307) at (-2.5, -2) {$f_4$};
		\node [style=backing] (308) at (3.25, -1.5) {$f_1$};
		\node [style=backing] (309) at (5.75, -1.5) {$f_2$};
		\node [style=backing] (310) at (8.25, -1.5) {$f_2$};
		\node [style=backing] (311) at (10.75, -1.5) {$f_3$};
		\node [style=backing] (312) at (13.25, -1.5) {$f_3$};
		\node [style=backing] (313) at (15.75, -1.5) {$f_4$};
		\node [style=medium back] (320) at (-14, 1.25) {};
		\node [style=none] (321) at (-14, 1.25) {$x,f_1$};
		\node [style=medium back] (322) at (-2, 1.25) {};
		\node [style=none] (323) at (-2, 1.25) {$x,f_4$};
		\node [style=medium back] (324) at (2, 0) {};
		\node [style=none] (325) at (2, 0) {$x,f_1$};
		\node [style=medium back] (326) at (17, 0) {};
		\node [style=none] (327) at (17, 0) {$x,f_4$};
		\node [style=none] (328) at (-12, -3.75) {$e_1$};
		\node [style=none] (329) at (-8, -3.75) {$e_2$};
		\node [style=none] (330) at (-4, -3.75) {$e_3$};
	\end{pgfonlayer}
	\begin{pgfonlayer}{edgelayer}
		\draw [style=c0] (267) to (273);
		\draw [style=c0] (272) to (273);
		\draw [style=c0] (272) to (268);
		\draw [style=c0] (268) to (275);
		\draw [style=c0] (275) to (274);
		\draw [style=c0] (274) to (269);
		\draw [style=c0] (269) to (277);
		\draw [style=c0] (277) to (276);
		\draw [style=c0] (276) to (270);
		\draw [style=c0, bend left=15] (266) to (271);
		\draw [style=d1] (279.center)
			 to [bend right=90, looseness=0.50] (278.center)
			 to [bend right=90, looseness=0.50] cycle;
		\draw [style=d1] (280.center)
			 to [bend right=90, looseness=0.50] (282.center)
			 to [bend right=90, looseness=0.50] cycle;
		\draw [style=d1] (281.center)
			 to [bend right=90, looseness=0.50] (283.center)
			 to [bend right=90, looseness=0.50] cycle;
		\draw [style=c0] (284.center) to (286.center);
		\draw [style=c0] (285.center) to (287.center);
		\draw [style=c0] (287.center) to (288.center);
		\draw [style=c0] (287.center) to (289.center);
		\draw [style=c0, bend left=15] (291) to (296);
		\draw [style=c0] (297) to (293);
		\draw [style=c0] (293) to (298);
		\draw [style=c0] (298) to (294);
		\draw [style=c0] (294) to (299);
		\draw [style=c0] (295) to (296);
		\draw [style=c0, bend right=15] (275) to (269);
		\draw [style=c0, bend right=15] (269) to (276);
		\draw [style=c0, bend left=15] (274) to (268);
		\draw [style=c0, bend left=15] (268) to (273);
		\draw [style=c0] (270) to (271);
		\draw [style=c0, bend right=15] (267) to (272);
		\draw [style=c0, bend left=15] (270) to (277);
		\draw [style=c0] (291) to (292);
		\draw [style=c0] (292) to (297);
		\draw [style=c0] (299) to (295);
		\draw [style=c0] (266) to (267);
	\end{pgfonlayer}
\end{tikzpicture}
    \caption{A constraint matroid $(M,c)$ with colours $x,f_1,f_2,f_3,f_4$, three distinguished colourless edges $e_1,e_2,e_3$ and a constraint minor $(C,c'')$ thereof. Marked with dotted lines are the separators corresponding to the totally-nested $2$-separations of $M$. We obtain $(C,c'')$ from $(M,c)$ by contracting $e_1,e_2,e_3$ and parallel merging.}
    \label{fig:irreducible1}
\end{figure}
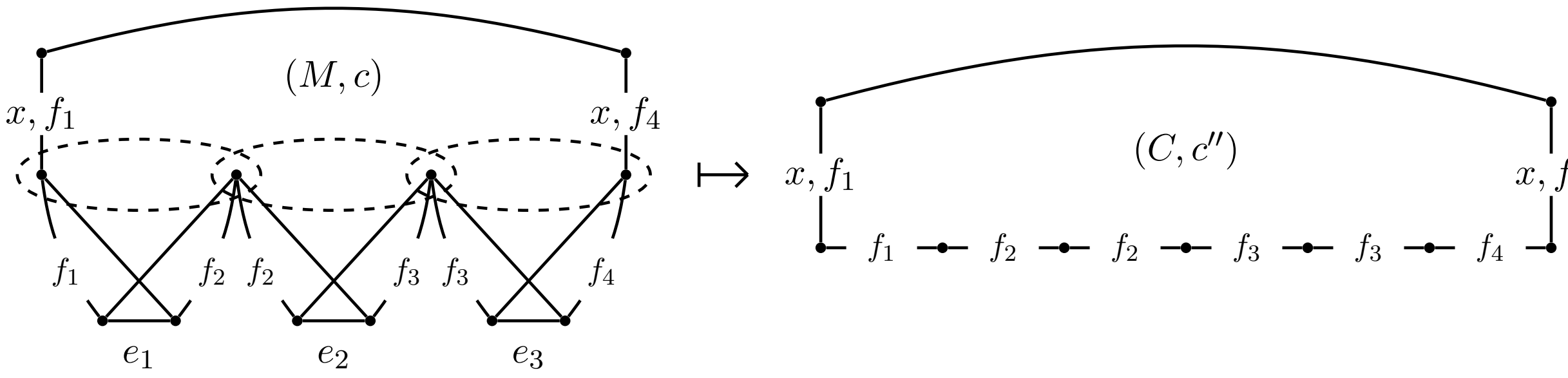
\begin{eg}
    See \autoref{fig:irreducible1}. The constraint matroid $(M,c)$ is obtained from a circuit via $2$-sums with three copies of $M(K_4)$, and it coloured according to the illustration with colours $x,f_1,f_2,f_3,f_4$. We observe that it is not realisable. Further observe that any proper constraint minor of $(M,c)$ is realisable. It is clear that we might extend this example to a minimally-unrealisable constraint matroid with $n\ge 4$ colours, and all such examples are not in our list of obstructions from \autoref{thm:forb}. Hence, they must not be colour-complete.

    Now consider the constraint matroid $(M,c')$ obtained from $(M,c)$ by adding a new colour $f_{(A,B)}$ for every totally-nested $2$-separation $(A,B)$ of $M$ which colours exactly $B$. Let $(C,c''')$ be the constraint minor obtained from $(M,c')$ by contracting $e_1,e_2,e_3$ and parallel merging all new parallel classes. Observe that now $(C,c''')$ is not realisable, and in fact it is a constraint circuit, a special case of a constraint radiolarian, which is in the list from \autoref{thm:forb}.
\end{eg}

\section{Connected constraint matroids: a mixed-sum tree-decomposition}
\label{sec:con}
\subsection{Overview}\label{subsec:overview_1}

The main result of this section is the following, which is a constraint analogue of \autoref{thm:matroids_tutte}.
\begin{thm}\label{thm:matroids_tutte_analogue}
    Let $(M,c)$ be a connected constraint matroid which is neither unbalanced nor unstable nor disagreeable. Then there exists a mixed-sum tree decomposition $(T,\{(M_u,c_u)\},\{\Ecal_u^e\})$ of $M$ into constraint circuits, constraint diatoms and $3$-connected constraint matroids so that the realisations of $(M,c)$ are in one-to-one correspondence with a choice of realisation for each $(M_u,c_u)$.
\end{thm}
In particular, the realisations of $(M,c)$ are obtained by mixed-sum along the adhesions $\Ecal_u^e$ from the realisations of the parts $(M_u,c_u)$.

In the remainder of this section, we will develop the necessary tools to prove the above theorem. In Subsections~\ref{subsec:semicircuits} and \ref{subsec:strong_weak} we introduce some important notation. In Subsections~\ref{subsec:unbalanced}, \ref{subsec:unstable} and \ref{subsec:disagreeable}, we introduce the three `non-local' obstructions to realisability. Finally, in \autoref{subsec:proof1}, we give a proof of \autoref{thm:matroids_tutte_analogue}.

\subsection{Semicircuits}\label{subsec:semicircuits}
Let $M$ be a connected matroid and $(A,B)$ a totally-nested $2$-separation of $M$. An \emph{semicircuit spanning $(A,B)$ in $A$} is the intersection of $A$ with a cycle $o$ that has nonempty intersection with both $A$ and $B$. In the context of a colouring $c$ for $M$ and a particular colour $f$, we say that a semicircuit is an $f$-semicircuit if all of its edges receive colour $f$.

\begin{dfn}
    Let $(M,c)$ be a connected constraint matroid and $(A,B)$ a totally-nested $2$-separation of $M$. For some colour $f$, we say that $(A,B)$ is \emph{$f$-spanned in $A$} if there is an $f$-semicircuit spanning $(A,B)$ in $A$. We say that $(A,B)$ is \emph{$f$-spanned} if it is $f$-spanned in either $A$ or $B$.
\end{dfn}

\begin{lem}\label{lem:f-spanned_equiv}
     Let $(M,c)$ be a connected constraint matroid and $(A,B)$ a totally-nested $2$-separation of $M$. Then the following are equivalent:\vspace{-.25cm}
     \begin{enumerate}
         \item the totally-nested $2$-separation $(A,B)$ is $f$-spanned; and\vspace{-.25cm}
         \item in every realisation $G$ of $M$, there is a path in $G$ between the separating vertices of $(A,B)$, all whose edges receive colour $f$; and\vspace{-.25cm}
         \item in some realisation $G$ of $M$, there is a path in $G$ between the separating vertices of $(A,B)$, all whose edges receive colour $f$.\vspace{-.25cm}
     \end{enumerate}
\end{lem}
\begin{pf}
    \textbf{(i)$\Rightarrow$(ii).} Let $G$ be a realisation of $M$ and $\{x,y\}$ denote the separator of $(A,B)$ in $G$. Let $o$ be a circuit in $M$ so that $A\cap o=S$. Then $o$, as a cycle in $G$, contains both $x$ and $y$ and the two segments of $o$ between $x$ and $y$ correspond exactly to the edges in $A$ and the edges in $B$. The segment in $A$ is comprised of the edges of $S$ and gives us our required path.

    \textbf{(ii)$\Rightarrow$(iii).} This follows trivially.

    \textbf{(iii)$\Rightarrow$(i).} Let $P$ be a path in $G$ between the vertices in the separator of $(A,B)$ whose edges all receive colour $f$. Then $P$ is contained completely in either $A$ or $B$, say $A$. Let $Q$ be a path in $G$ between the vertices in the separator of $(A,B)$ whose edges all lie in $B$. Then $P$ and $Q$ together form our required circuit in the matroid $M$.
\end{pf}

\begin{cor}\label{cor:f-span}
    Let $(M,c)$ be a connected constraint matroid and $(A,B)$ a totally-nested $2$-separation of $M$. Let $G$ be a realisation of $(M,c)$. If each of the vertices in the separator of $(A,B)$ are incident with an edge of colour $p$, then $(A,B)$ is $f$-spanned.
\end{cor}
\begin{pf}
    Since $G$ is a realisation of $(M,c)$ the edges receiving colour $f$ form a connected subgraph $F$ of $G$. Hence there is a path in $F$ between the vertices in the separator of $(A,B)$, all whose edges receive colour $f$. Our result follows by applying \autoref{lem:f-spanned_equiv}.
\end{pf}

\subsection{Strong and weak colouring}\label{subsec:strong_weak}
Let $M$ be a connected constraint matroid and $W$ a torso of the Tutte decomposition of $M$. Let $(A,B)$ be a totally-nested $2$-separator of $M$. We say that $(A,B)$ is \emph{left-incident} with $W$ if the edge corresponding to $(A,B)$ in the Tutte decomposition tree is incident with the node correspondng to $W$, and $B$ is one of the branches at a torso edge of $W$. Roughly speaking, that is when $W$ is the `nearest' torso to $B$ that is in some sense still `contained' in $A$. We say $(A,B)$ is \emph{right-incident} with $W$ when $(B,A)$ is left-incident with $W$.

Let $W$ be a torso of a connected matroid $M$. Let $c$ be a colouring of $M$ and define the \emph{torso colouring} $c_W$ of $W$ as follows: for each non-torso edge of $W$, we colour it the same as in $c$; and for each torso edge $e$ of $W$, corresponding to a totally-nested $2$-separation $(A,B)$ left-incident with $W$, we colour $e$ using all of the colours present in $B$.
Additionally, we introduce the following distinction in the way in which an edge in $W$ can receive a colour: we say that an edge $e$ in $W$ that is coloured by $f$ is \emph{strongly coloured by $f$} if either: (1) it is a non-torso edge; or (2) it is a torso edge of $W$ and the totally-nested $2$-separator $(A,B)$ left-incident with $W$ is $f$-spanned in $B$. Otherwise, we say that $e$ is \emph{weakly coloured by $f$}. Notice that edges of $W$ that are weakly coloured by some colour are necessarily torso edges. We call the constraint matroid\footnote{Formally speaking, this is no longer just a constraint matroid, since we also carry the information about which colourings are weak and which are strong. We hope that it is clear from context when this extra information is necessary.} $(W,c_W)$ a \emph{coloured torso} of~$(M,c)$.

Let $P$ be a path in a graph $G$. The \emph{end-edges} of $P$ are the first and last edges of $P$ (in the path order from $P$). All other edges are called \emph{internal}.
\begin{dfn}[$f$-strongly connected]\label{dfn:strongly_connected}
    Let $(M,c)$ be a connected constraint matroid and $(W,c_W)$ a coloured torso thereof. Let a realisation $H$ of $W$ be given. In this context, we say that a path $P$ in $H$ is an \emph{internally strong $f$-path} if all of its edges receive colour $f$ and no internal edges are weakly-coloured by $f$. Suppose that $H$ is not a triangle having only edges weakly-coloured in $f$. Then we say that $H$ is \emph{$f$-strongly connected} if every pair of edges that receive colour $f$ are the end-edges from an internally strong $f$-path in $H$.
\end{dfn}

\begin{figure}
    \centering
    \begin{tikzpicture}
	\begin{pgfonlayer}{nodelayer}
		\node [style=black dot] (98) at (-8, -1) {};
		\node [style=black dot] (99) at (-4.5, -1) {};
		\node [style=black dot] (100) at (-6.25, 2) {};
		\node [style=black dot] (103) at (-2.75, 0.25) {};
		\node [style=black dot] (104) at (-7.75, -3) {};
		\node [style=black dot] (105) at (-9.75, 0.25) {};
		\node [style=black dot] (106) at (-8.25, 3) {};
		\node [style=black dot] (109) at (-4.75, -3) {};
		\node [style=black dot] (110) at (-4.25, 3) {};
		\node [style=none] (120) at (-10, -2) {$(M,c)$};
		\node [style=black dot] (121) at (4.25, -1) {};
		\node [style=black dot] (122) at (7.75, -1) {};
		\node [style=black dot] (123) at (6, 2) {};
		\node [style=none] (124) at (6, -0.25) {$(W_0,c_{W_0})$};
		\node [style=black dot] (134) at (3.75, -0.75) {};
		\node [style=black dot] (135) at (2, 0.5) {};
		\node [style=black dot] (136) at (3.5, 3.25) {};
		\node [style=black dot] (137) at (5.5, 2.25) {};
		\node [style=black dot] (138) at (6.5, 2.25) {};
		\node [style=black dot] (139) at (8.25, -0.75) {};
		\node [style=black dot] (140) at (7.75, -1.5) {};
		\node [style=black dot] (141) at (4.25, -1.5) {};
		\node [style=black dot] (142) at (4.5, -3.5) {};
		\node [style=black dot] (143) at (7.5, -3.5) {};
		\node [style=black dot] (144) at (8.5, 3.25) {};
		\node [style=black dot] (145) at (10, 0.5) {};
		\node [style=none] (146) at (0.75, 2.25) {$(W_1,c_{W_1})$};
		\node [style=none] (147) at (11.25, 2.25) {$(W_2,c_{W_2})$};
		\node [style=none] (148) at (9.75, -2.5) {$(W_3,c_{W_3})$};
	\end{pgfonlayer}
	\begin{pgfonlayer}{edgelayer}
		\draw [style=c2] (98) to (105);
		\draw [style=c0] (105) to (100);
		\draw [style=c0] (100) to (106);
		\draw [style=c0] (106) to (105);
		\draw [style=c0] (106) to (98);
		\draw [style=c0] (98) to (104);
		\draw [style=c0] (104) to (99);
		\draw [style=c2] (99) to (103);
		\draw [style=c0] (104) to (109);
		\draw [style=c0] (109) to (98);
		\draw [style=c2] (109) to (99);
		\draw [style=c0] (103) to (110);
		\draw [style=c0] (110) to (100);
		\draw [style=c0] (103) to (100);
		\draw [style=c0] (99) to (110);
		\draw [style=c2 dashed] (121) to (123);
		\draw [style=c2 dashed] (123) to (122);
		\draw [style=c2 dashed] (122) to (121);
		\draw [style=c0] (135) to (136);
		\draw [style=c0] (136) to (137);
		\draw [style=c0] (137) to (135);
		\draw [style=c0] (136) to (134);
		\draw [style=c0] (138) to (144);
		\draw [style=c0] (144) to (145);
		\draw [style=c0] (145) to (138);
		\draw [style=c0] (139) to (144);
		\draw [style=c0] (141) to (142);
		\draw [style=c0] (142) to (143);
		\draw [style=c0] (143) to (141);
		\draw [style=c0] (142) to (140);
		\draw [style=c2 dashed] (141) to (140);
		\draw [style=c2 dashed] (139) to (138);
		\draw [style=c2 dashed] (137) to (134);
		\draw [style=c2] (140) to (143);
		\draw [style=c2] (139) to (145);
		\draw [style=c2] (134) to (135);
	\end{pgfonlayer}
\end{tikzpicture}
    \caption{A connected constraint matroid $(M,c)$ its coloured torsos. A dashed coloured line signifies that this torso edge is weakly-coloured using that colour.}
    \label{fig:weak_triangle}
\end{figure}
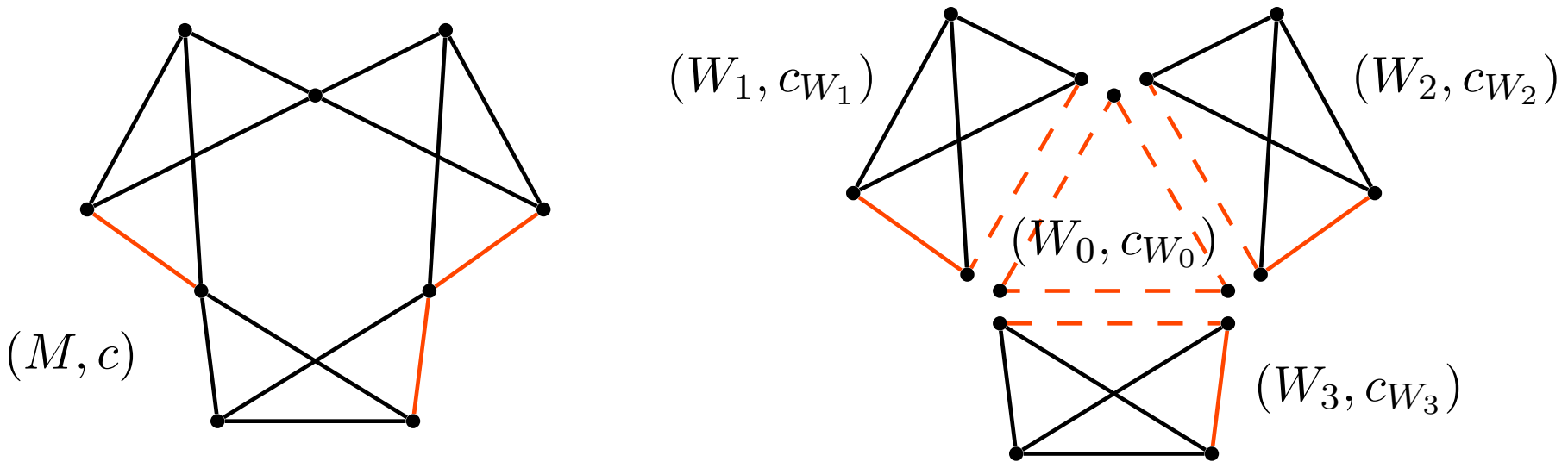
\begin{eg}
    In \autoref{fig:weak_triangle} we illustrate a connected constraint matroid $(M,c)$ along with a coloured $(W_0,c_{W_0})$ thereof. Observe that each edge of the coloured torso $(W_0,c_{W_0})$ is weakly-coloured by red. Each pair of edges forms the end-edges of a strong red-path, however $(M,c)$ admits no realisation. In fact, no connected constraint matroid with such a torso is realisable and this means our ommision of this case in \autoref{dfn:strongly_connected} ($f$-strongly-connected) will not cause us problems. The reason we make this definition is so that later definitions are smoother.
\end{eg}

\begin{lem}\label{lem:strong_torsos}
    Let $(M,c)$ be a connected constraint matroid and $(W,c_W)$ a coloured torso thereof. Every realisation $G$ of $(M,c)$ induces a realisation $H$ of $(W,c_W)$ that is $f$-strongly connected for each colour $f$.
\end{lem}
\begin{pf}
    Let $e_1$ and $e_2$ be edges in $H$ receiving colour $f$. Then since $G$ is a realisation of $(M,c)$, there is a path $P$ in $G$ whose end-edges are $e_1$ and $e_2$ and all whose edges receive colour $f$. Observe that $P$ is formed by concatenating non-torso edges of $W$ with paths in branches corresponding to torso edges of $W$. Let $P'$ be the path in $H$ obtained by replacing each of these subpaths in a branch by their corresponding torso edge. Observe that if one of these subpaths in a branch does not contain an end-edge of $P$, then it spans the totally-nested $2$-separator at this branch and hence by \autoref{lem:f-spanned_equiv} the corresponding torso edge is strongly-coloured by $f$. Since all non-torso edges in $P$ are by definition not strongly-coloured by $f$, we are done.
\end{pf}

\subsection{Circular coloured torsos: unbalanced}
\label{subsec:unbalanced}

\begin{dfn}[$r$-related]
    Let $(M,c)$ be a connected constraint matroid and $(W,c_W)$ a circular coloured torso thereof. Let $r$ be an edge of $W$.

    Let $\Rcal$ be the set of colours which weakly colour $r$ in $c_W$. For each colour $f\in \Rcal$, let $F$ denote the set of edges in $W-r$ which receive colour $f$. Let $F_S\subseteq F$ denote the edges in $W-r$ coloured strongly by $f$.
    For every pair $p,q\in\Rcal$, we say that $p$ is \emph{$r$-relatively contained in $q$} and write $p\le_rq$ if the following conditions hold:\vspace{-.25cm}
    \begin{enumerate}
        \item we have $P\subseteq Q$; and\vspace{-.25cm}
        \item if $P\neq E(W-r)$ and $Q=E(W-r)$, then $P\subseteq Q_S$.\vspace{-.25cm}
        \item if $P=Q= E(W-r)$, then $P_S=Q_S$.\vspace{-.25cm}
    \end{enumerate}
    For every pair $p,q\in\Rcal$, we say that $p$ and $q$ are \emph{$r$-related} and write $p\sim_r q$ if either $p\le_r q$ or~$q\le_r p$.
\end{dfn}

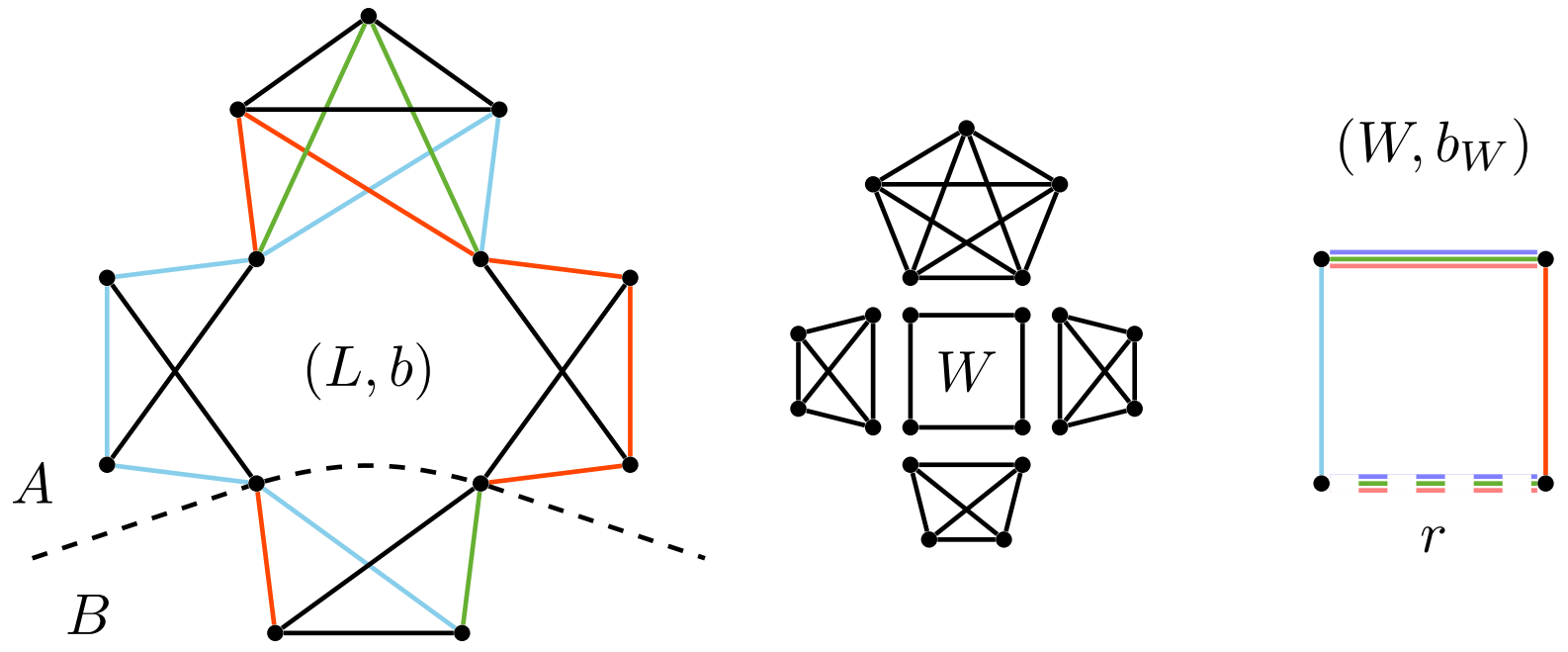
\begin{figure}
    \centering
    \begin{tikzpicture}
	\begin{pgfonlayer}{nodelayer}
		\node [style=black dot] (169) at (-24, 6) {};
		\node [style=black dot] (170) at (-24, 3) {};
		\node [style=black dot] (171) at (-27, 3) {};
		\node [style=black dot] (172) at (-27, 6) {};
		\node [style=black dot] (174) at (-23.75, 8) {};
		\node [style=black dot] (175) at (-22, 5.75) {};
		\node [style=black dot] (176) at (-22, 3.25) {};
		\node [style=black dot] (177) at (-24.25, 1) {};
		\node [style=black dot] (178) at (-26.75, 1) {};
		\node [style=black dot] (179) at (-29, 3.25) {};
		\node [style=black dot] (180) at (-29, 5.75) {};
		\node [style=black dot] (197) at (-27.25, 8) {};
		\node [style=black dot] (198) at (-25.5, 9.25) {};
		\node [style=black dot] (199) at (-18.75, 3.75) {};
		\node [style=black dot] (200) at (-18.75, 5.25) {};
		\node [style=black dot] (201) at (-19.75, 4) {};
		\node [style=black dot] (202) at (-19.75, 5) {};
		\node [style=black dot] (203) at (-18.25, 3.25) {};
		\node [style=black dot] (204) at (-18, 2.25) {};
		\node [style=black dot] (205) at (-17, 2.25) {};
		\node [style=black dot] (206) at (-16.75, 3.25) {};
		\node [style=black dot] (207) at (-16.25, 3.75) {};
		\node [style=black dot] (208) at (-16.25, 5.25) {};
		\node [style=black dot] (209) at (-15.25, 5) {};
		\node [style=black dot] (210) at (-15.25, 4) {};
		\node [style=black dot] (211) at (-18.25, 5.75) {};
		\node [style=black dot] (212) at (-16.75, 5.75) {};
		\node [style=black dot] (213) at (-18.75, 7) {};
		\node [style=black dot] (214) at (-16.25, 7) {};
		\node [style=black dot] (215) at (-17.5, 7.75) {};
		\node [style=black dot] (216) at (-18.25, 3.75) {};
		\node [style=black dot] (217) at (-18.25, 5.25) {};
		\node [style=black dot] (218) at (-16.75, 5.25) {};
		\node [style=black dot] (219) at (-16.75, 3.75) {};
		\node [style=black dot] (220) at (-12.75, 6) {};
		\node [style=black dot] (221) at (-9.75, 6) {};
		\node [style=black dot] (222) at (-9.75, 3) {};
		\node [style=black dot] (223) at (-12.75, 3) {};
		\node [style=none] (225) at (-17.5, 4.5) {$W$};
		\node [style=none] (226) at (-11.25, 7.5) {$(W,b_W)$};
		\node [style=none] (230) at (-11.25, 2.25) {$r$};
		\node [style=none] (231) at (-25.5, 4.5) {$(L,b)$};
		\node [style=none] (232) at (-30, 2) {};
		\node [style=none] (233) at (-21, 2) {};
		\node [style=none] (234) at (-30, 3) {$A$};
		\node [style=none] (235) at (-29.25, 1.25) {$B$};
	\end{pgfonlayer}
	\begin{pgfonlayer}{edgelayer}
		\draw [style=c1] (172) to (174);
		\draw [style=c1] (171) to (177);
		\draw [style=c1] (180) to (172);
		\draw [style=c2] (169) to (175);
		\draw [style=c2] (171) to (178);
		\draw [style=c1] (174) to (169);
		\draw [style=c3] (170) to (177);
		\draw [style=c0] (180) to (171);
		\draw [style=c1] (171) to (179);
		\draw [style=c1] (179) to (180);
		\draw [style=c0] (179) to (172);
		\draw [style=c0] (176) to (169);
		\draw [style=c0] (175) to (170);
		\draw [style=c2] (170) to (176);
		\draw [style=c2] (176) to (175);
		\draw [style=c0] (177) to (178);
		\draw [style=c0] (178) to (170);
		\draw [style=c2] (169) to (197);
		\draw [style=c2] (197) to (172);
		\draw [style=c3] (198) to (169);
		\draw [style=c3] (198) to (172);
		\draw [style=c0] (197) to (198);
		\draw [style=c0] (198) to (174);
		\draw [style=c0] (174) to (197);
		\draw [style=c0] (213) to (211);
		\draw [style=c0] (211) to (212);
		\draw [style=c0] (212) to (214);
		\draw [style=c0] (214) to (215);
		\draw [style=c0] (215) to (213);
		\draw [style=c0] (213) to (212);
		\draw [style=c0] (212) to (215);
		\draw [style=c0] (215) to (211);
		\draw [style=c0] (211) to (214);
		\draw [style=c0] (214) to (213);
		\draw [style=c0] (202) to (200);
		\draw [style=c0] (200) to (199);
		\draw [style=c0] (199) to (201);
		\draw [style=c0] (201) to (202);
		\draw [style=c0] (202) to (199);
		\draw [style=c0] (201) to (200);
		\draw [style=c0] (217) to (218);
		\draw [style=c0] (218) to (219);
		\draw [style=c0] (219) to (216);
		\draw [style=c0] (216) to (217);
		\draw [style=c0] (203) to (206);
		\draw [style=c0] (206) to (205);
		\draw [style=c0] (205) to (204);
		\draw [style=c0] (204) to (203);
		\draw [style=c0] (203) to (205);
		\draw [style=c0] (206) to (204);
		\draw [style=c0] (208) to (207);
		\draw [style=c0] (207) to (210);
		\draw [style=c0] (210) to (209);
		\draw [style=c0] (209) to (208);
		\draw [style=c0] (208) to (210);
		\draw [style=c0] (209) to (207);
		\draw [style=triple] (223) to (222);
		\draw [style=c2] (221) to (222);
		\draw [style=c1] (220) to (223);
		\draw [style=triple] (220) to (221);
		\draw [style=reallythickwhitedash] (223) to (222);
		\draw [style=d1] (232.center) to (171);
		\draw [style=d1, bend left=15] (171) to (170);
		\draw [style=d1] (170) to (233.center);
	\end{pgfonlayer}
\end{tikzpicture}
    \caption{The constraint matroid $(L,b)$ from \autoref{fig:obs1}, a totally-nested $2$-separation $(A,B)$, the Tutte decomposition of $M$ (centre), and the coloured torso $(W,b_W)$ of its circular torso $W$. A dashed coloured line in the coloured torso signifies that this torso edge is weakly-coloured using that colour. The torso edge $r$ corresponds to the separation $(A,B)$.}
    \label{fig:intransitive_nightmare}
\end{figure}

\begin{obs}\label{obs:this_one_is_actually_trans}
    The relation `$r$-relatively contained' is transitive, but `$r$-related' may not be.  \qed
\end{obs}

\begin{obs}
    The relation `$r$-related' is always reflexive and symmetric. \qed
\end{obs}

\begin{eg}
    In \autoref{fig:intransitive_nightmare}, we depict a connected constraint matroid $(M,c)$, along with a circular coloured torso $(W,c_W)$ having an edge $r$ so that the $r$-related relation is not transitive. Indeed, let $R,G,B$ denote the sets of red, green and blue edges in $W-r$, respectively, we get that neither $R,G$ nor $B$ is the whole of $E(W-r)$ but that $G\subseteq R$ and $G\subseteq B$ but neither $R\subseteq B$ nor $B\subseteq R$.
\end{eg}

\begin{lem}\label{lem:relation_when_rep}
    Let $(M,c)$ be a connected constraint matroid and $(W,c_W)$ a circular coloured torso thereof. Let $r$ be a torso edge of $W$ and $\Rcal$ denote the set of colours which weakly colour $r$ in the torso colouring $c_W$.
    If $(M,c)$ is realisable, then $\sim_r$ is an equivalence relation on $\Rcal$.
    
    Furthermore, in every realisation of $(M,c)$, the subsets of colours of $\Rcal$ of the edges incident with each of the vertices in the totally-nested $2$-separator corresponding to $r$ are the equivalence classes of $\sim_r$. In particular, $\sim_r$ has at most two equivalence classes.
\end{lem}
\begin{pf}
    Let $(A,B)$ be the totally-nested $2$-separation of $M$ corresponding to the torso edge $r$. Let $G$ be a realisation of $(M,c)$ and $O$ the cycle induced by $G$ which realises $(W,c_W)$. 
    We make the following claims.
    
    \begin{claim}\label{claim:chromodynamics}
        For each colour $f$ in $\Rcal$, the following holds:\vspace{-.25cm}
        \begin{enumerate}
            \item the cycle $O$ is $f$-strongly connected; and\vspace{-.25cm}
            \item there are at most two edges in $O$ weakly-coloured by $f$, one of which is $r$; and\vspace{-.25cm}
            \item not all of the edges of $W-r$ are strongly-coloured by $f$; and\vspace{-.25cm}
            \item the edges coloured by $f$ form an initial or final subpath of the path $O-r$; and\vspace{-.25cm}
            \item the edges strongly-coloured by $f$ form an initial or final subpath of the path $O-r$.\vspace{-.25cm}
        \end{enumerate}
    \end{claim}
    \begin{cproof}
        By \autoref{lem:strong_torsos}, the cycle $O$ is $f$-strongly connected.
        The fact that there are at most two edges in $O$ weakly-coloured by $f$ follows from the fact that $O$ is a cycle and is $f$-strongly connected. Observe that $r$ is weakly-coloured by $f$ because $(A,B)$ is not $f$-spanned. Since $O$ is a realisation for $(W,c_W)$, the edges coloured by $f$ form a segment in the linear order of the path $O-r$. We observe that the strongly-coloured edges have to form a subsegment because $O$ is $f$-strongly connected. Now observe that if all the edges of $O-r$ are strongly-coloured by $f$, then both endvertices of $r$ meet edges coloured by $f$ in $G$ and by \autoref{cor:f-span} this contradicts the fact that $(A,B)$ is not $f$-spanned.
    \end{cproof}\
    
    Let $x$ and $y$ denote the vertices in the separator of $(A,B)$ in $G$. We will also use $x$ and $y$ to denote the corresponding vertices in $O$ which are the endvertices of the torso edge $r$. Let $\Xcal$ and~$\Ycal$ denote the colours in $\Rcal$ which colour the edges meeting $x$ and $y$, respectively, in the side $A$ in the realisation $G$ of $(M,c)$. We will prove that $\Xcal$ and $\Ycal$ are exactly the equivalence classes of $\sim_r$. Certainly $\Xcal\cup \Ycal=\Rcal$, since at least one of $x$ or $y$ must be incident in $A$ with an edge in each colour appearing in both $A$ and $B$.
    \begin{claim}\label{claim:yes}
        Let $p,q\in X$. Then $p\sim_r q$.
    \end{claim}
    \begin{cproof}
        By \autoref{claim:chromodynamics} (iv) we get that either $P\subseteq Q$ or $Q\subseteq P$. Since otherwise we would be done, assume that either $P$ or $Q$ consist of all the edges of $W-r$. Without loss of generality, assume $P$ consists of all the edges of $W-r$. So all the edges of $W$ receive colour $p$. Let $e_x$ be the edge in $O-r$ that is incident with $x$ and $e_y$ the edge incident with $y$. Since $x$ is incident in $A$ with and edge coloured by $p$ and an edge coloured by $q$, and by \autoref{cor:f-span}, $y$ is not, we get that $e_x$ is coloured by both $p$ and $q$ and is strongly-coloured by $p$ and, in the case that $q$ coloured more than one edge in $O-r$, strongly-coloured by $q$. By \autoref{claim:chromodynamics} (v) we get that $Q_S\subseteq P_S$. Since there is at most one edge of $W-r$ coloured by $p$, respectively $q$, we get that either $Q_S=P_S$ and so $Q=P$, or $Q\subseteq P_S$. This completes the proof of our claim.
    \end{cproof}

    The above claim applies analogously to pairs $p,q\in \Ycal$.
    \begin{claim}\label{claim:no}
        Let $p\in \Xcal$ and $q\in \Ycal$. Then $p\nsim_r q$.
    \end{claim}
    \begin{cproof}
        Similarly to in the proof of \autoref{claim:yes} we observe that the edge $e_x$ incident with $x$ in $O-r$ is either strongly-coloured by $p$ or the only edge of $O-r$ coloured by $p$. Similarly, we have the edge $e_y$ incident with $y$ in $O-r$ is either strongly-coloured by $q$ or the only edge of $O-r$ coloured by~$q$. Condition (i) in the definition of $\sim_r$ is not satisfied since by \autoref{claim:chromodynamics} (iv) we have that each colour class forms a segment in the path $O-r$. Condition (ii), (iii) and (iv) are not satisfied for similar reasons.
    \end{cproof}

    Combining \autoref{claim:yes} and \autoref{claim:no} completes our proof.
\end{pf}

\begin{dfn}[Unbalanced]
    Let $(M,c)$ be a connected constraint matroid. We say a circular coloured torso $(W,c_W)$ of $(M,c)$ is \emph{unbalanced} if there is an edge $r$ of $W$ so that either:\vspace{-.25cm}
    \begin{enumerate}
        \item the relation $\sim_r$ (on the colours $\Rcal$ weakly colouring $r$) is not transitive; or\vspace{-.25cm}
        \item the relation $\sim_r$ is transitive, but it admits more than two equivalence classes.
    \end{enumerate}
\end{dfn}

\begin{cor}\label{cor:not_unbalanced}
    Let $(M,c)$ be a connected constraint matroid. If $(M,c)$ is realisable, then it contains no unbalanced circular coloured torsos.
\end{cor}
\begin{pf}
    We apply \autoref{lem:relation_when_rep} and the definition of unbalanced.
\end{pf}

\subsection{$3$-connected coloured torsos: unstable}
\label{subsec:unstable}

\begin{dfn}[Unstable]\label{dfn:stable}
    Let $(M,c)$ be a connected constraint matroid and $(W,c_W)$ a $3$-connected coloured torso thereof. Let $r$ be an edge of $W$ and $\Rcal$ denote the set of colours which weakly colour~$r$.
    
    We say that the coloured torso $(W,c_W)$ is \emph{stable with respect to $r$} if the following holds: Since $W$ is $3$-connected it has a unique graphic realisation, say $H$. Let $x$ and $y$ be the endvertices of $r$ in $H$. Then we have for each colour $f\in \Rcal$ that there is a vertex $z\in\{x,y\}$ so that either:\vspace{-.25cm}
    \begin{enumerate}
        \item only $z\in\{x,y\}$ sees $f$ in $H-r$; or\vspace{-.25cm}
        \item both $x$ and $y$ see $f$ in $H-r$ but exactly one $z\in\{x,y\}$ is incident with an edge strongly-coloured by $f$ in $H-r$.\vspace{-.25cm}
    \end{enumerate}
    If $(W,c_W)$ is not stable with respect to every edge in $W$, then we call it \emph{unstable}.
\end{dfn}

\begin{lem}\label{lem:not_unstable}
    Let $(M,c)$ be a connected constraint matroid and $(W,c_W)$ a $3$-connected coloured torso thereof. If $(M,c)$ is realisable, then $(W,c_W)$ is stable.
\end{lem}
\begin{pf}
    Let $r$ be an edge of $W$. Since otherwise it is trivial, suppose that $r$ is a torso edge. Let $(A,B)$ be the totally-nested $2$-separation of $M$ which corresponds to $r$. Let $G$ by a realisation of $(M,c)$ and $H$ the unique graph which realises $W$. By \autoref{lem:strong_torsos}, $H$ is a realisation of $(W,c_W)$. Let $x$ and $y$ denote both the vertices of the separator of $(A,B)$ in $G$ and the corresponding endvertices of $r$ in $H$. Observe that one of $x$ or $y$ must be incident with an edge coloured $f$ in $H-r$, since otherwise $H$ would not be a realisation of $(W,c_W)$.
    
    First suppose that both $x$ and $y$ are incident with edges receiving colour $f$ but that all such edges are weakly-coloured by $f$. Let $e_1$ and $e_2$ be such edges incident with $x$ and $y$ in $H-r$, respectively. Since $r$ is weakly-coloured, we have that there is no strong $f$-path in $H$ with $r$ and $e_1$ as its end-edges. But by \autoref{lem:strong_torsos}, the graph $H$ is strongly $f$-connected, a contradiction.
    
    Secondly suppose that both $x$ and $y$ are incident with edges receiving colour $f$ in $H-r$ and each one is incident with an edge strongly coloured by $f$ in $H-r$. Then both $x$ and $y$ are incident with edges coloured $f$ in $G$. But this is a contradiction because $(A,B)$ is not spanned by $f$.
\end{pf}

\subsection{Disagreeable separators}
\label{subsec:disagreeable}
In this subsection, we define what it means for a totally-nested $2$-separation to be `disagreeable'.

We say that a colour $f$ \emph{crosses} a totally-nested $2$-separation $(A,B)$ if both $A$ and $B$ contain edges which receive the colour $f$. We further say that $f$ \emph{weakly crosses} $(A,B)$ if it crosses $(A,B)$ but $(A,B)$ is not $f$-spanned.

Let $G$ be a graph, each of whose edges have been assigned some nonnegative number of colours, and let $z$ be a vertex of $G$. We say that $z$ \emph{sees} a set $\Scal$ of colours if for each $s\in \Scal$ there is an edge incident with $z$ in $G$ which is coloured using $s$.

\begin{figure}
    \centering
    \begin{tikzpicture}
	\begin{pgfonlayer}{nodelayer}
		\node [style=black dot] (149) at (-26.25, -1.25) {};
		\node [style=black dot] (150) at (-23.75, -1.25) {};
		\node [style=black dot] (151) at (-26.25, 1.75) {};
		\node [style=black dot] (152) at (-23.75, 1.75) {};
		\node [style=black dot] (153) at (-28.25, 1.5) {};
		\node [style=black dot] (154) at (-21.75, 1.5) {};
		\node [style=black dot] (155) at (-28.25, -1) {};
		\node [style=black dot] (156) at (-21.75, -1) {};
		\node [style=black dot] (157) at (-15.75, -1.25) {};
		\node [style=black dot] (158) at (-11.25, -1.25) {};
		\node [style=black dot] (159) at (-15.75, 1.75) {};
		\node [style=black dot] (160) at (-11.25, 1.75) {};
		\node [style=black dot] (161) at (-17.75, 1.5) {};
		\node [style=black dot] (162) at (-9.25, 1.5) {};
		\node [style=black dot] (163) at (-17.75, -1) {};
		\node [style=black dot] (164) at (-9.25, -1) {};
		\node [style=black dot] (165) at (-14.75, -1.25) {};
		\node [style=black dot] (166) at (-12.25, -1.25) {};
		\node [style=black dot] (167) at (-14.75, 1.75) {};
		\node [style=black dot] (168) at (-12.25, 1.75) {};
		\node [style=none] (169) at (-25.25, -2.25) {$(N,d)$};
		\node [style=none] (170) at (-23.75, 3.5) {};
		\node [style=none] (171) at (-23.75, -2.5) {};
		\node [style=none] (172) at (-24.5, 2.75) {$A$};
		\node [style=none] (173) at (-23, 2.75) {$B$};
	\end{pgfonlayer}
	\begin{pgfonlayer}{edgelayer}
		\draw [style=c2] (151) to (152);
		\draw [style=c1] (150) to (149);
		\draw [style=c2] (151) to (153);
		\draw [style=c0] (153) to (155);
		\draw [style=c1] (155) to (149);
		\draw [style=c0] (156) to (150);
		\draw [style=sideby] (152) to (154);
		\draw [style=c0] (154) to (156);
		\draw [style=c0] (153) to (149);
		\draw [style=c0] (155) to (151);
		\draw [style=c0] (152) to (156);
		\draw [style=c0] (150) to (154);
		\draw [style=c2] (159) to (161);
		\draw [style=c0] (161) to (163);
		\draw [style=c1] (163) to (157);
		\draw [style=c0] (164) to (158);
		\draw [style=sideby] (160) to (162);
		\draw [style=c0] (162) to (164);
		\draw [style=c0] (161) to (157);
		\draw [style=c0] (163) to (159);
		\draw [style=c0] (160) to (164);
		\draw [style=c0] (158) to (162);
		\draw [style=c1] (166) to (165);
		\draw [style=c2] (167) to (168);
		\draw [style=sideby] (159) to (157);
		\draw [style=sideby] (167) to (165);
		\draw [style=sideby] (168) to (166);
		\draw [style=sideby] (160) to (158);
		\draw [style=thickwhitedash] (159) to (157);
		\draw [style=thickwhitedash] (167) to (165);
		\draw [style=thickwhitedash] (168) to (166);
		\draw [style=thickwhitedash] (160) to (158);
		\draw [style=d1] (171.center) to (150);
		\draw [style=d1] (150) to (152);
		\draw [style=d1] (152) to (170.center);
	\end{pgfonlayer}
\end{tikzpicture}
    \caption{The connected constraint matroid $(N,d)$ from \autoref{fig:obs1}, along with its coloured torsos and a totally-nested $2$-separator $(A,B)$. A dashed coloured line signifies that this edge is weakly-coloured using that colour.}
    \label{fig:well-coloured_however}
\end{figure}
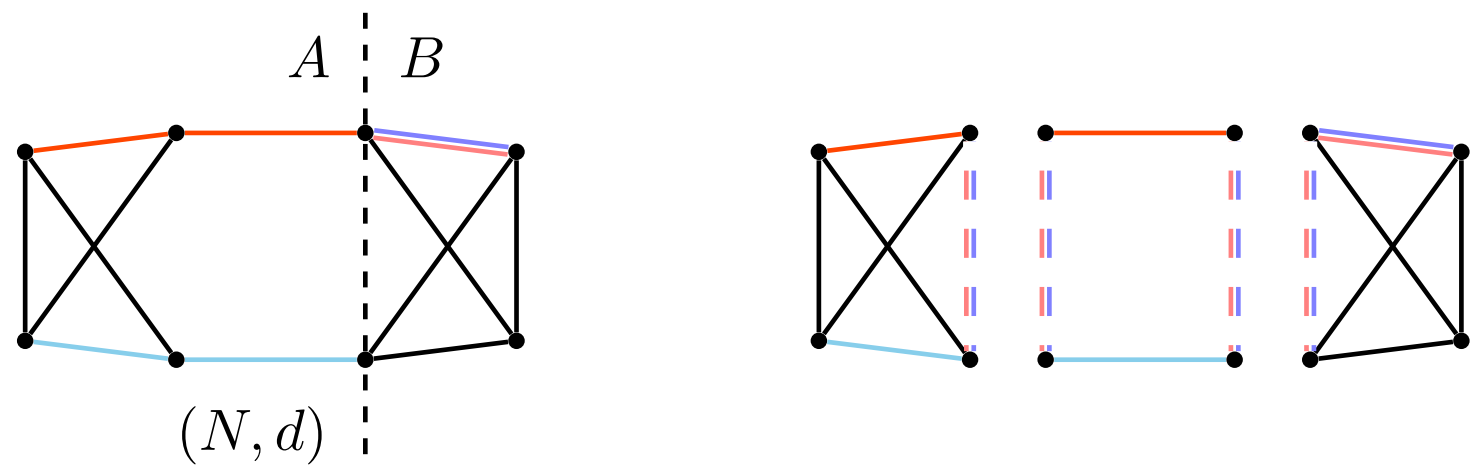

We say that a connected constraint matroid $(M,c)$ is \emph{balanced} or \emph{stable} if none of its coloured torsos are unbalanced or unstable, respectively.

\begin{dfn}
    Let $(M,c)$ be a balanced and stable connected constraint matroid and $(A,B)$ a totally-nested $2$-separation thereof left-incident with a torso $W$ which is not a bond. Let $r$ denote the torso edge of $W$ corresponding to $(A,B)$. Let $\Rcal$ denote the set of colours $f$ which weakly cross $(A,B)$. Then $\Rcal$ is exactly the set of colours that weakly colour $r$ in the torso colouring $c_W$.
    
    We define \emph{a left quasicolouring} of $(A,B)$ to be a bipartition $\Xcal\cup \Ycal$ of $\Rcal$ so that the following cases hold:
    
    \textbf{Case 1.} The torso $W$ is circular. Since it is not unbalanced with respect to $(A,B)$, the $r$-related relation $\sim_r$ is an equivalence relation on $\Rcal$ with at most two equivalence classes. If there is at most one class, set $\Xcal:=\Rcal$ and $\Ycal:=\varnothing$. Otherwise, set $\Xcal$ and $\Ycal$ to be distinct equivalence classes.
    
    \textbf{Case 2.} The torso $W$ is $3$-connected. Let $H$ denote its unique graphic realisation and $x$ and $y$ the endvertices of $r$ in $H$. We assign a colour $f\in \Rcal$ to $\Xcal$ if $x$ is the vertex satifying (i) or (ii) in the definition of stable (\ref{dfn:stable}). That is, either: (i) only $x\in\{x,y\}$ sees $f$ in $H-r$; or (ii) both $x$ and $y$ see $f$ in $H-r$ but only $x\in\{x,y\}$ is incident with an edge strongly-coloured by $f$ in $H-r$. All other colours $f\in \Rcal$ are assigned to $\Ycal$.
\end{dfn}
We have that Case 2 is unchanged when we swap the roles of the symbols $x$ and $\Xcal$ with the roles of and $y$ and $\Ycal$, respectively. With this fact, we observe that left quasicolourings are unique up to swapping the roles of $\Xcal$ and $\Ycal$. Hence we may refer to `a' left quasicolouring as `the' left quasicolouring.

Let $G$ be a graph, each of whose edges have been assigned some nonnegative number of colours, and let $z$ be a vertex of $G$. Let $\Scal$ and $\Tcal$ be a sets of colours with $\Scal\subseteq\Tcal$. We say that $z$ sees \emph{exactly $\Scal$ from $\Tcal$} if there is no edge incident with $z$ that is coloured using some colour in~$\Tcal\setminus \Scal$.

For a graph $G$ and a subset $A\subseteq E(G)$ of edges, write $G[[A]]$ for the smallest subgraph of $G$ which contains these edges. We say that $G[[A]]$ is the subgraph of $G$ \emph{induced} on the edge set $A$.
\begin{lem}\label{lem:quasi3con}
    Let $(M,c)$ be a balanced and stable connected constraint matroid and $(A,B)$ a totally-nested $2$-separation thereof which is left-incident with a $3$-connected torso $W$. Let $\Rcal$ denote the set of colours which weakly cross $(A,B)$. Then the left quasicolouring $\Xcal\cup \Ycal=\Rcal$ of $(A,B)$ has the following property: In every realisation $G$ of $(M,c)$, one vertex of the separator of $(A,B)$ sees exactly $\Xcal$ from $\Rcal$ and the other sees exactly $\Ycal$ from $\Rcal$.
\end{lem}
\begin{pf}
    Since no colour in $\Rcal$ spans $(A,B)$, it suffices to prove our hypothesis in the restriction $G[[A]]$. That is, in $G[[A]]$, one vertex of the separator of $(A,B)$ sees exactly $\Xcal$ from $\Rcal$ and the other sees exactly $\Ycal$ from $\Rcal$.
    
    To that end, let $H$ denote the unique graphic representation of $W$ used in the definition of the left quasicolouring and $r$ the torso edge of $W$ corresponding to the separation $(A,B)$. Let $x$ and $y$ denote the vertices in the separator of $(A,B)$ in $G$ and their corresponding endvertices of $r$ in $H$.

    For each $f\in \Rcal$, since it crosses $(A,B)$ and $G$ is a realisation of $(M,c)$, we have that exactly one vertex of the separator of $(A,B)$ meets an edge of colour $f$ in $G[[A]]$.
    
    If $z\in\{x,y\}$ is incident with an edge of colour $f$ in $G[[A]]$, then clearly it is incident with an edge of colour $f$ in $H-r$. Hence by the definition of the left quasicolouring, it remains to prove that if both $x$ and $y$ are incident with an edge of colour $f$ in $H$, and $z\in\{x,y\}$ is incident with an edge strongly-coloured by $f$, then $z$ is the vertex in $\{x,y\}$ which is incident with an edge coloured $f$ in $G[[A]]$. Let $e$ denote an edge incident with $z$ in $H$, strongly coloured by $f$. If $e$ is not a torso edge of $W$, then it is an edge in $G[[A]]$. If it is a torso edge, then its endvertices $z$ and $z'$ are from the separator of a totally-nested $2$-separation of $M$ (and $G$) spanned by $f$. By \autoref{lem:f-spanned_equiv}, there is a path $P$ in $G$ between $z$ and $z'$, all whose vertices are coloured by $f$. Since $z'$ lies in $G[[A]]$ but not $G[[B]]$, we get that $z$ is incident with an edge of colour $f$ in $G[[A]]$.
\end{pf}

\begin{cor}\label{lem:quasicolouring_main}
    Let $(M,c)$ be a balanced and stable connected constraint matroid and $(A,B)$ a totally-nested $2$-separation thereof which is left-incident with a torso $W$ which is not a bond. Let $\Rcal$ denote the set of colours which weakly cross $(A,B)$. Then the left quasicolouring $\Xcal\cup \Ycal=\Rcal$ of $(A,B)$ has the following property: In every realisation $G$ of $(M,c)$, one vertex of the separator of $(A,B)$ sees exactly $\Xcal$ from $\Rcal$ and the other sees exactly $\Ycal$ from $\Rcal$.
\end{cor}
\begin{pf}
    We simply combine \autoref{lem:relation_when_rep} and \autoref{lem:quasi3con}.
\end{pf}

A left quasicolouring for a $2$-separation $(A,B)$ is also known as a \emph{right quasicolouring} for $(B,A)$.
\begin{dfn}[Disagreeable]
    Let $(M,c)$ be a balanced and stable connected constraint matroid and $(A,B)$ a totally-nested $2$-separator thereof. Suppose that neither torso incident with $(A,B)$ is a bond. We say that $(M,c)$ is \emph{disagreeable at $(A,B)$} if the left and right quasicolourings of $(A,B)$ are not the same. We say that $(M,c)$ is \emph{disagreeable} if it is disagreeable at one of its totally-nested $2$-separators.
\end{dfn}

\begin{cor}\label{cor:not_disagreeable}
     Let $(M,c)$ be a balanced and stable connected constraint matroid and $(A,B)$ a totally-nested $2$-separator thereof. If $(M,c)$ is realisable, then $(M,c)$ is agreeable.
\end{cor}
\begin{pf}
    Use \autoref{lem:quasicolouring_main}.
\end{pf}

\begin{eg}
    In \autoref{fig:well-coloured_however} we depict a connected constraint matroid $(M,c)$ along with its coloured torsos. The constraint matroid $(M,c)$ is neither unstable nor unbalanced, however it does not admit a realisation. A certificate for this can be found by observing that the totally-nested separator $(A,B)$ is disagreeable. Indeed, the left quasicolouring of $(A,B)$ has blue and red in different parts and the right quasicolouring of $(A,B)$ groups both red and blue into the same part.
\end{eg}

In the case where our connected constraint matroid is balanced, stable and agreeable, we use the term \emph{quasiolouring} of a totally-nested $2$-separator to describe the left quasicolouring of one of its orientations. Since every totally-nested $2$-separator has an orientation which is left-incident with some circular or $3$-connected torso, we have that there is a well-defined quasicolouring for each totally-nested $2$-separator.

\subsection{The proof of \autoref{thm:matroids_tutte_analogue}}\label{subsec:proof1}
In the following, we define a way to split a stable, balanced and agreeable connected constraint matroid along a totally-nested $2$-separation into two smaller constraint matroids. We will then prove that this proceedure and its reverse preserve realisablility.

Let $(M,c)$ be a stable, balanced and agreeable connected constraint matroid and $(A,B)$ a totally-nested $2$-separation thereof. Let $A^*$ and $B^*$ be the matroids with edges sets $A\sqcup\{e_\Xcal,e_\Ycal\}$ and $B\sqcup\{e_\Xcal,e_\Ycal\}$, respectively, so that $M=A^*\oplus_{(\{e_\Xcal,e_\Ycal\},\{e_\Xcal,e_\Ycal\})} B^*$. Let $\Xcal\sqcup \Ycal=\Rcal$ be the quasicolouring of $(A,B)$ and let $f_{(A,B)}$ denote a new colour which does not appear in $c$. Let $c_A$ be the colouring of~$A^*$ obtained from the restriction of $c$ on to the edges $A$ but where we make the following additional colourations:\vspace{-.25cm}
\begin{enumerate}
    \item the edge $e_\Xcal$ attains each colour $f\in \Xcal$; and\vspace{-.25cm}
    \item the edge $e_Y$ attains each colour $f\in \Ycal$; and\vspace{-.25cm}
    \item the edges $e_\Xcal$ and $e_\Ycal$ both attain each colour $f$ which crosses $(A,B)$ so that $(A,B)$ is $f$-spanned (i.e. crossing colours not in $\Rcal$); and\vspace{-.25cm}
    \item the edges $e_\Xcal$ and $e_\Ycal$ both attain the new colour $f_{(A,B)}$.\vspace{-.25cm}
\end{enumerate}
Let $c_B$ be the colouring defined analoguously for $B^*$. That is, defined to be the restriction of $c$ to the egdes $B$ along with $e_X$ and $e_Y$ which are coloured identically to in $c_A$ except with a new colour $f_{(B,A)}$ in place of $f_{(A,B)}$.

\begin{lem}
    Let $(M,c)$ be a stable, balanced and agreeable connected constraint matroid and $(A,B)$ a totally-nested $2$-separation thereof. Then from every realisation $G$ of $(M,c)$, we can obtain a realisation $G_A$ of $(A^*,c_A)$ by deleting $B$, adding one new vertex, and adding two new edges with labels $e_X$ and $e_Y$.
\end{lem}
\begin{pf}
    Let $x$ and $y$ denote the vertices of $G$ forming the separator of the separation $(A,B)$. By \autoref{lem:quasicolouring_main}, one of the vertices $x$ or $y$ sees $\Xcal$ from $\Rcal$ and the other sees $\Ycal$ from $\Rcal$. Without loss of generality, assume $x$ sees $\Xcal$ and $y$ sees $\Ycal$, or that $\Rcal$ is empty. Furthermore, since $G$ is a realisation of $(M,c)$, for each colour crossing $(A,B)$, there is at least one of $x$ or $y$ incident with an edge in that colour. Let $G_A$ be the graph obtained from $G[[A]]$ by adding a new vertex $v$ and the edge $e_\Xcal$ from $x$ to $v$ and the edge $e_\Ycal$ from $v$ to $y$. One may easily check that $G_A$ is indeed a realisation of $(A^*,c_A)$.
\end{pf}
We call the graph $G_A$ obtained in the above lemma, a \emph{realisation induced on $(A^*,c_A)$ by~$G$}. We remark that $G_A$ might not be unique, even up to relabeling vertices. Indeed, in the case where~$\Rcal$ is empty, our new edges attach depending on how we choose $x$ and $y$.

\begin{lem}\label{lem:replace_by_path}
    Let $(M,c)$ be a stable, balanced and agreeable connected constraint matroid and $(A,B)$ a totally-nested $2$-separation thereof. Then for every pair $G_A$ and $G_B$ of realisations for $(A^*,c_A)$ and $(B^*,c_B)$, respectively, there is a realisation of $(M,c)$, and furthermore it induces both $G_A$ on $(A^*,c_A)$ and $G_B$ on $(B^*,c_B)$.
\end{lem}
\begin{pf}
    Let $G_A$ and $G_B$ be realisations for $(A^*,c_A)$ and $(B^*,c_B)$, respectively. Since the colour $f_{(A,B)}$ colours only the series edges $e_\Xcal$ and $e_\Ycal$ in $G_A$ and $G_B$, the pair corresponds to the edges of a path $P$ in both $G_A$ and $G_B$ whose internal vertex is not incident with any other edges in $G_A$ or $G_B$, respectively.

    Let $G$ be the graph obtained from $G_A$ and $G_B$ by identifying them along the path $P$ and then deleting the internal vertex of this path and its incident edges $e_\Xcal$ and $e_\Ycal$. We claim that $G$ is a realisation of $(M,c)$. For a colour $f$, let $A_f$ and $B_f$ denote the subgraphs of $G[[A]]$ and $G[[B]]$, respectively, induced on the $f$-coloured edges. For any colour $f$ that does not cross $(A,B)$, either $A_f$ or $B_f$ is empty and the nonempty one is connected. So let $f$ by a colour which crosses $(A,B)$. Let $x$ and $y$ denote the vertices of $G$ that form the separation of $(A,B)$. Now suppose $(A,B)$ is $f$-spanned. Then by \autoref{lem:f-spanned_equiv}, there is an $f$-coloured path $P$ in $G$ between $x$ and $y$, contained in one of $A_f$ or $B_f$. Furthermore, both $A_f$ and $B_f$ consist either of a single connected component, meeting one or more of the vertices from $\{x,y\}$, or two connected components, each of which meets a different one of the vertices $x$ and $y$. Hence the $f$-coloured edges form a connected subgraph of~$G$.
    Finally suppose that $(A,B)$ is not $f$-spanned. That is, $f\in \Rcal=\Xcal\sqcup \Ycal$. Then the subgraphs $A_f$ and $B_f$ are connected and each meet one of the vertices $x$ or $y$. Since one of the edges $e_\Xcal$ or $e_\Ycal$ (not both) receive the colour $f$ in both $G_A$ and $G_B$, we have that both $A_f$ and $B_f$ meet the same vertex from $\{x,y\}$, and hence the subgraph of $f$-coloured edges is connected in $G$. This completes the proof.
\end{pf}

\begin{lem}\label{lem:goodness_remains}
    If $(M,c)$ is a stable, balanced and agreeable connected constraint matroid with a totally-nested $2$-separation $(A,B)$, then $(A^*,c_A)$ and $(B^*,c_B)$ are stable, balanced and agreeable.
\end{lem}
\begin{pf}
    We will prove that $(A^*,c_A)$ is stable, balanced and agreeable. The proof for $(B^*,c_B)$ is analogous. To prove these properties, we will first explicate the Tutte decomposition of $A^*$.

    Let $\Tcal=(T,\{M_u\},\{\Ecal_u^e\})$ denote the Tutte decomposition of $G$. Let $e$ be the edge of $T$ corresponding to $(A,B)$ and let $t$ be the endnode of $e$ whose corresponding torso is left-incident with $(A,B)$. First suppose that $M_t$ is a circular torso. Let $T'$ be the connected component of $T-e$ containing $t$ and let $M_t'$ be obtained from $M_t$ by replacing the torso edge corresponding to $(A,B)$ by the series pair $e_\Xcal,e_\Ycal$. Then the Tutte decomposition of $A^*$ is obtained from the restriction $(T',\{M_u\},\{\Ecal_u^e\})$ of $\Tcal$ by replacing $M_t$ by $M_t'$. In this case, it is an elementary exercise to check that $(A^*,c_A)$ is stable, balanced and agreeable.

    Now suppose that $M_t$ is a bond or is $3$-connected and let $e_t$ denote the torso edge therein corresponding to $(A,B)$. Let $T''$ be the tree obtained from $T'$ by appending a new leaf $s$ at $t$ and let $M_s$ be the $3$-circuit consisting of $e_\Xcal$, $e_\Ycal$ and a new edge $e_s$. Then $(A,\{e_\Xcal,e_\Ycal\})$ is a totally-nested $2$-separation of $A^*$ and the Tutte decomposition of $A^*$ is the tree decomposition $(T'',\{M_u\},\{\Ecal_u^e\})$ where the $\Ecal_t^{st}=\{e_t\}$ and $\Ecal_s^{et}=\{e_s\}$.
    In this case, we observe that the coloured torsos of $A^*$ are all stable and balanced and that the left quasicolourings of $(A,\{e_\Xcal,e_\Ycal\})$ and $(\{e_\Xcal,e_\Ycal\},A)$ are the same, namely $\{\Xcal,\Ycal\}$. This completes the proof.
\end{pf}

Let $(M,c)$ be a connected constraint matroid and $W$ a torso of $M$. Let $c$ be a colouring of $M$ and define the \emph{strict torso colouring} of $W$ as follows: for each non-torso edge of $W$, we colour it the same as in $c$; and for each torso edge $e$ of $W$, corresponding to a totally-nested $2$-separation $(A,B)$ left-incident with $W$, we colour $e$ using all of the colours $f$ present in $B$ so that $(A,B)$ is $f$-spanned in $B$. That is, the strict torso colouring is the torso colouring where we forget colour assignments which are weak.

We say that a tree decomposition is \emph{star-shaped} if its auxillary tree is a star. The \emph{central torso} of a star-shaped tree decomposition is the torso corresponding to the central node of the star.
Let $(M,c)$ be a connected constraint matroid and suppose that its Tutte decomposition is star-shaped and that its central torso is $3$-connected. Suppose that each torso of the Tutte decomposition corresponding to a leaf of the star is a circuit with $3$ edges. That is, we obtain $M$ from the central torso by subdividing the torso edges once. If furthermore $(M,c)$ is stable, balanced and agreeable, then we call it an \emph{extended $3$-block}.

\begin{lem}\label{lem:3_block_rectify}
    Let $(M,c)$ be an extended $3$-block. Then $(M,c)$ is realisable if and only if its central torso is realisable with respect to its strict torso colouring.
\end{lem}
\begin{pf}
    Let $W$ denote the central $3$-connected torso of $M$ and $c'$ its strict torso colouring.
    
    Firstly, suppose that $G$ is a realisation of $(M,c)$ and $H$ is a minor of $G$ realising to $W$ (perhaps after relabelling its torso edges). Let~$f$ be a colour that appears in $c'$. Let $p$ and $q$ denote edges in $H$ coloured by $f$. Suppose first that $p$ and $q$ are not torso edges. Then there is an $f$-coloured path $P$ in $G$ whose first and last edges are $p$ and $q$. One may replace all of the subpaths of $P$ by that are not in $H$ by their corresponding torso edges to obtain an $f$-coloured path joining $p$ and $q$ in $H$. Now suppose $p$ is a torso edge and $q$ is not. In this case, we look for an $f$-coloured path $P$ in $G$ whose first vertex is one of the separator vertices from the totally-nested $2$-separation corresponding to $p$ and whose last edge is $q$. By replacing all of the subpaths of $P$ that are not in $H$ by their corresponding torso edges, we are again done. The two remaining cases are (1) that $p$ is not a torso edge and $q$ is a torso edge; and (2) that both $p$ and $q$ are torso edges. In both of these cases, the reasoning is analogous and we get that the $f$-coloured edges form a connected subgraph of $H$, as required.

    Suppose conversely that $H$ is a realisation of $(W,c')$.  Let $c_W$ denote the torso colouring of~$H$ (not the strict torso colouring). Since $(W,c_W)$ is stable, we have for each totally-nested $2$-separation $(A,B)$ of $G$ which is left-incident with $W$ that there is a well-defined left quasicolouring $\Xcal\cup \Ycal=\Rcal$ of the colours $\Rcal$ which weakly cross $(A,B)$. Additionally, since $(M,c)$ is agreeable, this quasicolouring agrees with the left quasicolouring of $(B,A)$. Let $r$ denote the torso edge of $W$ corresponding to $(A,B)$. If $\Rcal$ is nonempty, then there is the following correspondence between the vertices $x,y$ of the separator of $(A,B)$ and the quasicolouring classes~$\Xcal$ and $\Ycal$:\vspace{-.22cm}
    \begin{enumerate}
        \item for each $f\in \Xcal$, $x$ is the vertex satisfying (i) or (ii) in the definition of stable (\ref{dfn:stable}). That is, either: (i) only $x\in\{x,y\}$ is incident with an edge of colour $f$ in $H-r$ under the torso colouring~$c_W$; or (ii) both $x$ and $y$ are incident with edges of colour $f$ in $H-r$ but only $x\in\{x,y\}$ is incident with an edge strongly-coloured by $f$ in $H-r$ under the torso colouring~$c_W$.\vspace{-.25cm}
        \item for each $f\in \Ycal$, $y$ satisfies the above criterion with the symbols `$x$' and `$y$' swapped.\vspace{-.25cm}
    \end{enumerate}
    Let $G$ be the graph obtained from $H$ by replacing each torso edge of $H$ by a path of length two. We will now describe how to assign each colour used by $c$ to the edges of $G$, so that $G$ is a realisation of $(M,c)$.
    Let $C_1,\ldots,C_n$ be triangle graphs realising the non-central torsos of $(M,c)$. Let $G$ be the graph obtained from $H$ by taking the $2$-sum of $H$ with each $C_i$ along their corresponding torso edge, so that the following holds:
    Let $p_i$ and $q_i$ denote the non-torso edges of $C_i$. Let $\Xcal\sqcup \Ycal=\Rcal$ denote the quasicolouring corresponding to the $2$-separator formed from the $2$-sum of $C_i$. Then without loss of generality $p_i$ receives all the colours from $\Xcal$ under $c$ and $q_i$ receives all of the colours of $\Ycal$ under $c$.
    
    We take the $2$-sum so that $p_i$ is incident with $x$ and $q_i$ is incident with $y$, where $x$ and $y$ are the vertices corresponding to $\Xcal$ and $\Ycal$, respectively, as outlined above.
    Observe that for each $f$ colouring $p_i$, either:\vspace{-.25cm}
    \begin{enumerate}
        \item $p_i$ is incident with an edge coloured $f$ in the strict torso colouring of $H$, which form a connected set in $H$, since $H$ is a realisation of $(W,c')$; or\vspace{-.25cm}
        \item $p_i$ is incident with no edge coloured $f$ in the strict torso colouring $c'$ of $H$. Instead, $x$ is the unique vertex incident with edges weakly coloured by $f$ in the torso colouring $c_W$ of $H$.\vspace{-.25cm}
    \end{enumerate}
    We conclude that $G$ is a realisation of $(M,c)$, as required.
\end{pf}

The following is the final ingredient in our proof of \autoref{thm:matroids_tutte_analogue}. \begin{lem}\label{lem:pre_tutte}
    Let $(M,c)$ be a connected constraint matroid which is neither unbalanced nor unstable nor disagreeable. Then there is a mixed-sum tree decomposition $(T,\{(M_u,c_u)\},\{\Ecal_u^e\})$ of $M$ into constraint circuits, constraint diatoms and extended $3$-blocks, so that the following holds: $(M,c)$ is realisable if and only if each $(M_u,c_u)$ is realisable.
\end{lem}
\begin{pf}
    We will define a mixed-sum tree decomposition $(T,\{M_u,c_u\},\{\Ecal_u^e\})$ recursively. At each step, we will make sure all of our torsos are stable, balanced and agreeable. First we take our mixed-sum tree decomposition to consist of a single bag containing the whole of $(M,c)$. While there is a bag of $(T,\{M_u,c_u\},\{\Ecal_u^e\})$ which contains a totally-nested $2$-separation whose sides both have at least three edges, we replace $(T,\{M_u,c_u\},\{\Ecal_u^e\})$ by $(T',\{M_u,c_u\},\{\Ecal_u^e\})$, which is defined as follows: The tree $T'$ is obtained from $T$ by subdividing a node $x$ corresponding to a bag with such a $2$-separation $(A_x,B_x)$ into two nodes $a$ and $b$. Here we have to specify which edges are then incident with the new vertices $a$ and $b$. We use the following claim.
    \begin{claim}
        Every pair of torso edges in $M_x$ corresponding to the same edge in $T$ lie both in the same side of $(A_x,B_x)$.
    \end{claim}
    \begin{cproof}
        Since $e_1$ and $e_2$ are in series, the partition $\{\{e_1,e_2\},E(M_x)\setminus\{e_1,e_2\}\}$ defines a $2$-separation of $M_x$. Since $(A_x,B_x)$ is totally-nested, we get that either $\{e_1,e_2\}\subseteq A_x$ or $\{e_1,e_2\}\subseteq B_x$, as required.
    \end{cproof}
    
    We define $T'$ so that the edges incident with $a$ in $T'$ are those corresponding to torso edges in $A_x$ and the edges incident with $b$ in $T'$ are those corresponding to torso edges in $B_x$.

    Now we have to specify the matroids corresponding to the new nodes $a$ and $b$ or $T$. Let $(M_a,c_a)=(A^*,c_A)$ and $(M_b,c_b)=(B^*,c_B)$ from \autoref{lem:replace_by_path} and $\Ecal_a^{ab}=\Ecal_b^{ab}:=\{e_x,e_y\}$. For each $v\neq b$ a neighbour of $a$ in $T'$, let $\Ecal_a^{av}=\Ecal_v^{av}:=\Ecal_x^{xv}$, and for each $v\neq a$ a neighbour of $b$ in $T'$, let $\Ecal_b^{bv}=\Ecal_v^{bv}:=\Ecal_x^{xv}$. By \autoref{lem:goodness_remains}, all of our torsos are still stable, balanced and agreeable.
    By \autoref{lem:replace_by_path}, the tree decomposition $(T',\{M_u,c_u\},\{\Ecal_u^e\})$ is a tree decomposition of $(M,c)$.
    
    Our process terminates when all of the torsos have no totally-nested $2$-separator whose sides both contain at least three edges. Hence they are circuits, diatoms or extended $3$-blocks.
\end{pf}

\begin{pf}[Proof of \autoref{thm:matroids_tutte_analogue}]
    Let $(M,c)$ be a connected constraint matroid which is neither unbalanced nor unstable nor disagreeable. Then we are required to prove that there exists a mixed-sum tree decomposition $(T,\{(M_u,c_u)\},\{\Ecal_u^e\})$ of $M$ into constraint circuits, constraint diatoms and $3$-connected constraint matroids so that the following holds: $(M,c)$ is realisable if and only if each $(M_u,c_u)$ is realisable.
    
    Let $(T,\{(N_u,d_u)\},\{\Fcal_u^e\})$, denote the tree decomposition obtained via \autoref{lem:pre_tutte}. For each node $u$ of $T$, we let $(M_u,c_u)$ be the matroid $(N_u,d_u)$ if $N_u$ is a circuit or a diatom. Otherwise we set $(M_u,c_u)$ equal to the central torso of the extended $3$-block $M_u$ along with its strict torso colouring. For each $\Fcal_u^{uv}$ where $M_u$ is not an extended $3$-block, let $\Ecal_u^{uv}=\Fcal_u^{uv}$. If $M_u$ is an extended $3$-block, let $\Ecal_u^{uv}:=\{e_v\}$ be the singleton containing the torso edge of the central torso of $M_u$ corresponding to the totally-nested $2$-separation one of whose sides is~$\Fcal_u^{uv}$.
    
    The edges incident with $u$ in $T$ correspond to the torso edges of the central torso that correspond to the original series classes. Our result then follows by appyling \autoref{lem:3_block_rectify}.
\end{pf}

\section{First steps in finding forbidden substructure}\label{sec:first_steps}
\subsection{Overview}
In what follows, we prove that the non-local obstructions actually result in forbidden substructure. In doing so, we will also develop the necessary tools to let us turn the presence of unrealisable bags in the decomposition from \autoref{thm:matroids_tutte_analogue} into forbidden minors.

\subsection{Toolbox}

\begin{obs}\label{lem:carry_it_all}
    Let $(A,B)$ be a totally-nested $2$-separation in a constraint matroid $(M,c)$. Let $f$ be a colour used by $c$ that spans $(A,B)$ in $B$. Let $(M',c')$ be a constraint minor of $(M,c)$ obtained by contracting and parallel-merging all but a subset of edges in $B$, say $B'\subseteq B$, so that $(A,B')$ is a (perhaps improper) $2$-separation of $(M',c')$. Then $f$ spans $(A,B')$ in $B'$, in the colouring $c'$. \qed
\end{obs}

\begin{dfn}[Semibond splitting $(A,B)$ in $B$]
    Let $(A,B)$ be a totally-nested $2$-separation of a connected matroid $M$. A set $D\subseteq B$ is said to be a \emph{semibond splitting $(A,B)$ in $B$} if it is not a bond of $M$ but exends to a bond of $M$ after adding some edges from $A$. Informally, if after minoring $A$ down to a single edge~$t$, the set $D\cup\{t\}$ is a bond. Observe that a semibond is exactly a semicircuit in the dual.
\end{dfn}

\begin{lem}[Colour-avoiding bond]
    Let $(A,B)$ be a totally-nested $2$-separation in a constraint matroid $(M,c)$. Let $f$ be a colour used by $c$ that does not span $(A,B)$ in $B$. Then there is a semibond splitting $(A,B)$ in $B$ that does not contain an edge coloured $f$. 
\end{lem}
\begin{pf}
    Let $G$ be a graphic realisation of the matroid $M$ and $x$ and $y$ the vertices of the separator of $(A,B)$ in $G$. Let $G[[B]]$ be the subgraph of $G$ induced on $B$. Let $H$ be a maximal connected subgraph of $G[[B]]$ containing $x$, all whose edges are coloured $f$ (which might just consist of the vertex $c$). Observe that $y\notin V(H)$ since otherwise there would be an $f$-coloured path from $x$ to $y$ in $G[[B]]$. Let $D$ be the set of edges of $G[[B]]$ with exactly one endvertex in $H$ and one in the connected component of $G[[B]]-E(H)$ containing $y$. The set $D$ forms our required semibond splitting $(A,B)$ in~$B$.
\end{pf}

\begin{cor}\label{lem:cut_control}
    Let $(A,B)$ be a totally-nested $2$-separation in a constraint matroid $(M,c)$. Let $f$ be a colour used by $c$ that does not span $(A,B)$ in $B$. Then there is a constraint minor $(M',c')$ of $(M,c)$ obtained by contracting and parallel-merging all but one edge in $B$, which is not coloured $f$ in $c'$. \qed
\end{cor}

\begin{lem}[Colour-grabbing bond]
    Let $(A,B)$ be a totally-nested $2$-separation in a constraint matroid $(M,c)$. Let $f$ be a colour used by $c$ that appears in $B$. Then there is a semibond splitting $(A,B)$ in $B$ that contains at least one edge coloured $f$. 
\end{lem}
\begin{pf}
    Let $e$ be some edge in $B$ which is coloured by $f$. Since $M$ is connected, it is contained in some bond $D$ of $M$ which also contains an element from $A$. The part of this bond in $B$ is our required semibond.
\end{pf}

\begin{cor}\label{lem:cut_grab}
    Let $(A,B)$ be a totally-nested $2$-separation in a constraint matroid $(M,c)$. Let $f$ be a colour used by $c$ that appears in $B$. Then there is a constraint minor $(M',c')$ of $(M,c)$ obtained by contracting and parallel-merging all but one edge in $B$, which is coloured $f$ in $c'$. \qed
\end{cor}

\begin{lem}\label{lem:unstable_never_happens}
    Let $(M,c)$ be a connected constraint matroid which is colour-complete. If $(M,c)$ is unstable (hence unrealisable), then it contains as a constraint minor a strictly smaller connected constraint matroid which is colour-complete and also unrealisable.
\end{lem}
\begin{pf}
Let $(W,c_W)$ be an unstable $3$-connected coloured torso of $(M,c)$. Let $r$ be a torso edge of $W$ witnessing this instability and let $\Rcal$ denote the set of colours which weakly colour~$r$ in $c_W$.
Let $H$ denote the unique graphic realisation of $W$ and $x$ and $y$ be the endvertices of $r$ in $H$. Using the definition of unstable, we have that there exists a colour $f\in \Rcal$ so that either\vspace{-.25cm}
\begin{enumerate}
    \item neither $x$ nor $y$ sees $f$ in $H$; or\vspace{-.25cm}
    \item both $x$ and $y$ see $f$ in $H-r$ but neither is incident in $H-r$ to an edge coloured strongly by~$f$.\vspace{-.25cm}
\end{enumerate}
In the first case, we can use \autoref{lem:cut_grab} to replace $r$ by a non-torso edge coloured $f$. What we are left with is still unrealisable, contradicting minimality.

Now consider case (ii). Let $E_1$ and $E_2$ be the sets of edges in $H-r$ adjacent to $x$ and $y$, respectively. Let $(M',c')$ be obtained from $(M,c)$ by contracting all of the edges in $(M,c)$ that are not: (1) in the branch at $r$ in $M$; (2) in $E_1$ or $E_2$; or (3) in the branches at the torso edges in $E_1$ and $E_2$. Then $(M',c')$ has a triangular torso in place of $W$, all whose edges are weakly coloured by~$f$. This is strictly smaller than $M$ and still unrealisable, contradicting minimality.
\end{pf}

\begin{dfn}[Diatomic]
    Let $(M,c)$ be a connected constraint matroid and $(A,B)$ a totally-nested $2$-separation of $M$. Then we say that $B$ is \emph{diatomic} if the $B$ side in the $2$-sum that decomposes $M$ along $(A,B)$ is a diatom and the sum is along a series class containing one edge. That is, $B$ consists in $M$ of series classes of size at most two which are all joined in parallel with~$(A,B)$.
\end{dfn}

\begin{eg}\label{eg:aaaaa}
    The result \autoref{lem:unstable_never_happens} will be employed so that we can use the following results, which all require stability. This is because of the following possibility for $B$: pick vertices $x$ and $y$ of $K_4$ and delete the edge between them. Now colour one of the edges incident with $x$ and one of the edges incident with $y$ blue, but in such a way that there is no blue path between $x$ and $y$. This graph can't be minored down to a path of length two, one of whose edges is blue. It might be that the $B$ side of a totally-nested $2$-separation looks like this example, with $x,y$ the separator. But then we can deal with it via \autoref{lem:unstable_never_happens} since it would be unstable.
\end{eg}

\begin{lem}[Crude torsoing with one colour]\label{lem:crude_one-torsoing}
    Let $(M,c)$ be a connected constraint matroid which is stable and $(A,B)$ a totally-nested $2$-separation of $M$ crossed but not spanned by some colour $f$. Then there is a constraint minor $(M',c')$ obtained from $(M,c)$ by contracting and parallel-merging some edges in $B$, except say $B'\subseteq B$, so that; (1) the pair $(A,B')$ is a $2$-separation of $(M',c')$ with $B'$ diatomic; and (2) the colour $f$ colours some edges in $B'$ but does not span $(A,B')$.
\end{lem}
\begin{pf}
    We split into cases based on the type of coloured torso $(W,c_W)$ right-incident with $(A,B)$ in $M$. We begin with special cases, Case * and Case 0.

    \textbf{Case *.} The torso $W$ is a circuit that contains an edge $e$ different from $r$ that is weakly coloured by $f$. Then we may contract all the edges of $W-r$ except $e$ (or their corresponding branches in the case where they're torso edges). Since $e$ corresponds to some totally-nested $2$-separation $(\hat{A},\hat{B})$ of $M$ with $\hat{B}\subseteq B$ and $|\hat{B}|<|B|$, we can apply induction to obtain our required result.

    \textbf{Case 0.} The torso $W$ is a triangle $re_1e_2$ where $e_1$ is coloured strongly by $f$ and $e_2$ is not coloured by $f$.
    If $e_1$ is a torso edge, we replace the branch at $e_1$ by a single edge coloured $f$ using \autoref{lem:cut_grab} and, if $e_2$ is a torso edge, we replace the branch at $e_2$ by a single edge not coloured $f$ using \autoref{lem:cut_grab}. We are done.

    \textbf{Case 1.} The torso $W$ is a circuit with more than three edges, and contains no edges weakly coloured by $f$ except $r$.
    Observe that then $W$ contains two edges aside from $r$, say $e_1$ and $e_2$, so that exactly one of $e_1,e_2$ is coloured strongly by $f$. Now consider the constraint matroid obtained from $(M,c)$ by contracting all the edges in $B$ outside of $e_1,e_2$ or their corresponding branches. We find ourselves in Case 1.

    \textbf{Case 2.} The torso $W$ is $3$-connected. Let $E_1$ and $E_2$ be the sets of edges adjacent to the endvertices of $r$ in the unique graphic realisation of $W$. Since we assumed that $(W,c_W)$ is stable, at most one of $E_1,E_2$ contains an edges strongly coloured by $f$. Consider the matroid obtained from $(M,c)$ by contracting all the edges in $B$ outside of $E_1,E_2$ or the corresponding branches. We find ourselves in Case * or Case 1.

    \textbf{Case 3.} The torso $W$ is a bond. For each torso edge $e$ of $W-r$, let $(A_e,B_e)$ be its corresponding totally-nested $2$-separation so that $W$ is left-incident. By induction, we can contract and parallel-merge edges in $B_e$ for each $e$ so to obtain a diatomic set $B_e'$. Observe that the union of all of these $B_e'$ in our matroid is a diatomic $B'$, with our desired property.
\end{pf}

\begin{lem}[Three weak edges]\label{lem:no_three_weak}
    Let $(M,c)$ be a colour-complete connected constraint matroid that is stable. Let $(W,c_W)$ be a circular coloured torso of $(M,c)$. If there is a colour $f$ so that $(W,c_W)$ has at least three edges coloured weakly by $f$, then $(M,c)$ contains an unrealisable radiolarian as a constraint minor.
\end{lem}
\begin{pf}
    Let $e_1,e_2$ and $e_3$ denote three torso edges of $W$ coloured weakly by $f$ in $c_W$ and $(A_i,B_i)$ their corresponding $2$-separations in $M$, with $W$ left-incident. Let $X$ denote the set of all edges of $M$ outside of $B_1\cup B_2\cup B_3$. By applying \autoref{lem:crude_one-torsoing} at each $(A_i,B_i)$, and then contracting $X$, we obtain a constraint radiolarian $(R,c)$. Observe that since $(M,c)$ is colour complete, each set of edges $B_i'$ produced from application of \autoref{lem:crude_one-torsoing} comes with a unique colour colouring all of $B_i'$. Now observe that $(R,c)$ is not realisable.
\end{pf}

\begin{lem}[Crude torsoing with two colours]\label{lem:crude_two-torsoing}
    Let $(M,c)$ be a connected constraint matroid which is stable and $(A,B)$ a totally-nested $2$-separation of $M$ crossed but not spanned by distinct colours $f_1$ and $f_2$. Then either $(M,c)$ contains an unrealisable radiolarian as a constraint minor, or there is a constraint minor $(M',c')$ obtained from $(M,c)$ by contracting and parallel-merging some edges in $B$, except say $B'\subseteq B$, so that; (1) the pair $(A,B')$ is a $2$-separation of $(M',c')$ with $B'$ diatomic; and (2) the colours $f_1$ and $f_2$ colour some edges in $B'$ but do not span $(A,B')$.
\end{lem}
\begin{pf}
    We proceed by induction on the size of $B$. For the base case, suppose that $B$ is already diatomic. In this case, our matroid is already of the required form. Now suppose that $B$ is not diatomic. Let $W$ be the torso right-incident with $(A,B)$ and $r$ the corresponding torso edge in $W$. We first deal with a special case, and then reduce other possibilities to this case.

    \textbf{Case 0.} The torso $W$ is a triangle $re_1e_2$ with the following properties: at most one of each $e_1,e_2$ is coloured strongly by each $f_1,f_2$ and at most one of each $e_1,e_2$ is coloured weakly by each $f_1,f_2$. Without loss of generality, there are five subcases.

    \textit{Subcase a.} The edge $e_1$ is coloured strongly by $f_1$ the edge $e_2$ is coloured strongly by $f_2$. In this case, since neither $f_1$ nor $f_2$ span $B$, if $e_1$ is coloured by $f_2$, then it is coloured weakly and if $e_2$ is coloured $f_1$, then it is coloured weakly. If either $e_1$ or $e_2$ is a torso edge, we use \autoref{lem:cut_control} at the totally-nested $2$-separations corresponding to it, to replace it by an edge avoiding $f_2$ or $f_1$, respectively. We get that $B'$ consists of two edges, one coloured $f_1$ and the other coloured $f_2$.

    \textit{Subcase b.} The edge $e_1$ is coloured strongly by $f_1$ and weakly by $f_2$ and $e_2$ is not coloured by $f_2$. If $e_2$ is coloured $f_1$, it is coloured weakly by it, and we use \autoref{lem:cut_control} to replace $e_2$ by a colourless edge. Then we use \autoref{lem:cut_grab} to replace $e_1$ by an edge coloured $f_2$. By \autoref{lem:carry_it_all}, this edge also has colour $f_1$.

    \textit{Subcase c.} The edge $e_1$ is coloured strongly by $f_1$ and $f_2$ and $e_2$ is coloured $f_1$ weakly and not coloured $f_2$. Then we use $\autoref{lem:cut_control}$ at at the totally-nested $2$-separations corresponding to $e_2$ to replace it by a colourless edge and we replace $e_1$ by an edge coloured both $f_1$ and $f_2$.

    \textit{Subcase d.} The edge $e_1$ is coloured strongly by $f_1$ and $f_2$ and $e_2$ is coloured weakly by both $f_1$ and $f_2$. Let $(\hat{A},\hat{B})$ be the totally-nested $2$-separation corresponding to $e_2$. Observe that $|\hat{B}|<|B|$ and so by induction we can replace $\hat{B}$ by a diatomic $\hat{B}'$ containing both $f_1$ and $f_2$ but spanned by neither. Now we contract $e_1$ (or in the case that it is a torso edge, its branch) and are done.

    \textit{Subcase e.} The edge $e_1$ is coloured weakly by $f_1$ and weakly by $f_2$ and $e_2$ is coloured by neither $f_1$ nor $f_2$. In this case, let $(\hat{A},\hat{B})$ be the totally-nested $2$-separation corresponding to $e_1$. Observe that $|\hat{B}|<|B|$ and so by induction we can replace $\hat{B}$ by a diatomic $\hat{B}'$ containing both $f_1$ and $f_2$ but spanned by neither. Now we contract $e_2$ (or in the case that it is a torso edge, its branch) and are done.

    \textbf{Case 1.} The torso $W$ is a circuit with more than three edges. In this case, we claim that there exist $e_1,e_2$ distinct edges in $W-r$ so that: (1) at most one of $e_1,e_2$ is coloured strongly by each $f_1,f_2$; (2) at most one of each $e_1,e_2$ is coloured weakly by each $f_1,f_2$; but (3) both colours are present. Indeed, it is clear that we can find a pair $e_1,e_2$ so that (1) and (3) hold. If (2) also holds, which happens for example when the torso $W$ is a balanced circuit, we can contract everything outside of these two branches and proceed as in Case 0.
    Otherwise, there are at least two edges coloured weakly by $f$ in $W-r$, and hence at least three in $W$ and we find a radiolarian obstruction by \autoref{lem:no_three_weak}.

    \textbf{Case 2.} The torso $W$ is $3$-connected. By assumption it is not unstable. Let $H$ be the unique graphic realisation of $W$ and $r$ the torso edge corresponding to $(A,B)$. Denote by $x$ and $y$ the endvertices of $r$ in $H$ and denote by $H'$ the graph obtained from $H$ by contracting all edges not incident with $x$ or $y$. Observe that $H'$ has three vertices, two of them being $x,y$. Hence we obtained a matroid whose torso in the place of $W$ is a triangle, so we continue as in Case 0.

    \textbf{Case 3.} The torso $W$ is a bond. In this case, we treat each non-$r$ torso edge in $W$ using the inductive hypothesis. What we end up with is replacing $W$ by a $2$-sum of a collection of diatoms, which is again a diatom.
\end{pf}

In the above lemma, we don't get control over whether $f_1$ and $f_2$ colour different or the same edges from the series classes of $B'$. In the following (\autoref{lem:careful_two-torsoing}), we get to keep the quasicolouring in some sense, but have to make some sacrifices: in particular, the right-incident torso is not allowed to be a diatom, and we have to assume our constraint matroid is colour-complete.

\begin{dfn}[Diatomic for quasicolourings]
    Let $(M,c)$ be a connected constraint matroid and $(A,B)$ a totally-nested $2$-separation thereof. Let $f_1$ and $f_2$ be colours which cross $(A,B)$ but do not span it. Then we say that $B$ is \emph{diatomic for $\{\{f_1\},\{f_2\}\}$} if in every graphic realisation $G$ of $M$, so that the edges coloured $f_1$ and $f_2$, respectively, form connected sets in $G[[B]]$, then the separator vertices of $(A,B)$ each see one of $f_1$ and $f_2$ in $G[[B]]$. We say that $B$ is \emph{diatomic for $\{\{f_1,f_2\},\varnothing\}$} if the same condition holds except there is exactly one vertex of $(A,B)$ that sees both $f_1$ and $f_2$ in~$G[[B]]$.
\end{dfn}

\begin{lem}[Careful torsoing with two colours]\label{lem:careful_two-torsoing}
    Let $(M,c)$ be a colour-complete connected constraint matroid which is stable and $(A,B)$ a totally-nested $2$-separation of $M$. Suppose that $(A,B)$ is crossed by exactly two colours which don't span it, say $f_1$ and $f_2$. Suppose additionally that the torso right-incident with $(A,B)$ is $3$-connected or a circuit which is balanced and let $\{\Xcal,\Ycal\}$ be the right-quasicolouring of $(A,B)$. Then either $(M,c)$ contains an unrealisable radiolarian as a constraint minor, or there is a constraint minor $(M',c')$ obtained from $(M,c)$ by contracting and parallel-merging some edges in $B$, except say $B'\subseteq B$, so that the pair $(A,B')$ is a $2$-separation of $(M',c')$ with $B'$ diatomic for $\{\Xcal,\Ycal\}$.
\end{lem}
\begin{pf}
    The proof proceeds very similarly to the proof of \autoref{lem:crude_two-torsoing}. Instead of using an induction hypothesis when we replace edges weakly coloured by both $f_1$ and $f_2$, we simply use \autoref{lem:crude_two-torsoing}. In what follows, we will outline how this proof differs from the proof of \autoref{lem:crude_two-torsoing}. We will first work with the subcases of Case 0, since Case 1 and Case 2 reduce to that. Case 3 does not occur since we assumed that $(A,B)$ is not right-incident with a bond torso.

    Subcases a, b and c are the same as in \autoref{lem:crude_two-torsoing}, however we observe that our outcome is always diatomic for $\{\Xcal,\Ycal\}$.

    In Subcase d: We don't contract $e_1$, instead we replace it by an edge coloured both $f_1$ and $f_2$. Then we use \autoref{lem:crude_two-torsoing} on $r$, with the roles of $A$ and $B$, reversed to obtain a diatomic $A'$. We are left with a radiolarian formed from $\hat{B}'$, $A'$ and $e_1$. Since $(M,c)$ was colour-complete, we get that there is a colour unique to the sets $\hat{B}'$ and $A'$, which causes them to be connected in any realisation of our new constraint matroid. Now observe that our radiolarian has no realisation.

    In Subcase e: We don't contract $e_2$, instead we replace it with an edge coloured neither $f_1$ nor $f_2$. Similarly to in Subcase d, we use \autoref{lem:crude_two-torsoing} on $r$, with the roles of $A$ and $B$ reversed, to obtain a diatomic $A'$. We are left with an unrealisable radiolarian formed from $\hat{B}'$, $A'$ and $e_1$.
\end{pf}

Finally, we need to do something similar to \autoref{lem:crude_two-torsoing} but for three crossing colours.

\begin{dfn}[Almost diatomic]
    Let $(M,c)$ be a connected constraint matroid and $(A,B)$ a totally-nested $2$-separation of $M$. Then we say that $B$ is \emph{almost diatomic} if
    
    \begin{tabular}{l p{10cm}}
    (diatomic) & we have $B$ is diatomic; or \\
    ($3$-series) & we have $B$ is $3$ edges in series or obtained from that by adding an edge spanning $(A,B)$; or \\
    ($P$-shape) & we have $B$ is obtained from the $3$-series case by adding an edge spanning two of the edges in series.
\end{tabular}
\end{dfn}

\begin{lem}[Crude torsoing with three colours]\label{lem:crude_three-torsoing}
    Let $(M,c)$ be a connected constraint matroid which is stable and $(A,B)$ a totally-nested $2$-separation of $M$ crossed but not spanned by three distinct colours $f_1,f_2,f_3$. Then either $(M,c)$ contains a radiolarian obstruction, a $M(K_4)$-shaped obstruction, or there is a constraint minor $(M^*,c^*)$ obtained from $(M,c)$ by contracting and parallel-merging some edges in $B$, except say $B^*\subseteq B$, so that; (1) the pair $(A,B^*\subseteq B)$ is a $2$-separation of $(M^*,c^*)$; (2) the side $B^*$ contains all three colours $f_1,f_2,f_3$ and furthermore $(A,B^*)$ is not spanned by $f_1$ or $f_3$; and (3) $B^*$ is almost diatomic.
\end{lem}
\begin{pf}
    First apply \autoref{lem:crude_two-torsoing} at $(A,B)$ with the coloures $f_1$ and $f_3$. That is, we ignore $f_2$. If $f_2$ is present in $B'$, then set $B^*=B'$ we are automatically done. So assume that this is not the case.
    
    Observe that there is a constraint minor $(M^*,c^*)$ of $(M,c)$ obtained by contracting and parallel merging edges in $B$, except say $B'$, so that: (1) $(M',c')$ can be obtained from $(M^*,c^*)$ by contracting one edge $e$ and then parallel merging; and (2) the edge $e$ receives colour $f_2$. The remainder of the proof is split into cases based on how the edge $e$ appears in $B'$. Let $N$ be the matroid obtained from $M$ by replacing $A$ in $(A,B)$ by a single spanning edge $r$. We will not deal with the case where $B'$ is diatomic, since in this case we are immediately done.

    \textbf{Case 1.} The edge $e$ and $r$ form a $K_4$ in $N$ with $e$ and $r$ not adjacent in any graphic realisation of $N$. Since $A$ contains some edge coloured $f_2$, we use \autoref{lem:cut_grab} to find a $K_4$ constraint-minor of $M$ so that is unrealisable because its two $f_2$-coloured edges are not adjacent.

    \textbf{Case 2.} The edge $e$ forms a series class of size $3$. Since otherwise we could contract the edge not coloured $f_1,f_2$ or $f_3$, we know each edge in the series class is exactly $f_1,f_2$ or $f_3$, respectively. In this case, it remains to show that we can reduce all the other series classes in $B'$ to a single edge. This is possible except when there is a series class with one edge coloured $f_1$ and one coloured $f_3$. But in this case, we can contract the $f_1$ or $f_3$ edge in our series class of length $3$ and obtain a diatomic side with our required properties.

    \textbf{Case 3.} The edge $e$ forms a the triangle in a $P$-shaped part, perhaps alongside other series classes in $B'$. Similarly to in Case 2, we get that we can reduce all these other series classes to one edge, or find something diatomic which satisfies our criteria.
\end{pf}

Now we have the necessary tools to fully deal with unbalanced circular coloured torsos.

\begin{lem}\label{lem:unbalanced_resolved}
    Let $(M,c)$ be a colour-complete connected constraint matroid which is stable. Let $(W,c_W)$ be a circular coloured torso of $(M,c)$. If $(W,c_W)$ is unbalanced, then $(M,c)$ contains either a radiolarian obstruction, a $M(K_4)$-shaped obstruction or an unbalanced hexagon as a constraint minor.
\end{lem}
\begin{pf}
    By \autoref{lem:no_three_weak}, we may assume that for no colour $f$ does $(W,c_W)$ have at least three edges coloured weakly by $f$. Since $(W,c_W)$ is unbalanced, there is an edge $r\in E(W)$ so that the relation $\sim_r$ is either intransitive, or is an equivalence relation but has more than two equivalence classes. The remainder of the proof is in two sections based on this distinction.

    \textbf{Intransitivity.} The relation $\sim_r$ is intransitive. Then there exists crossing colours $a,b,c$ so that $a\sim_r b\sim_r c$ but $a\nsim_r c$. Let $A$, $B$ and $C$ (likewise, $A_S$, $B_S$, $C_S$) denote the sets of edges in $E(W-r)$ receiving colours $a$ $b$ and $c$ (likewise, strongly), respectively. Without loss of generality, we may assume that either: (1) $a\le_r b$ and $c\le_r b$; or (2) $b\le_r a$ and $b\le_r c$. Indeed, we cannot have $a\le_r b\le_r c$, since we already observed that $\le_r$ is transitive~(\autoref{obs:this_one_is_actually_trans}).

    \textbf{Case 1.} We have that $a\le_r b$ and $c\le_r b$ but $a\nsim_rc$. Observe that neither $A$ nor $C$ is equal to all of $E(W-r)$. Indeed, if this were the case, then we would have either $A_S=B_S$ and $C\subseteq A_S$ or $C_S=B_S$ and $A\subseteq C_S$, implying either $c\le_r a$ or $a\le_r c$, respectively.
    
    \textit{Subcase 1a.} We have $B=E(W-r)$. Then $A,C\subseteq B_S$. Thus there exists edges $e_a,e_b,e_c$ in $W-r$ satisfying: $e_a\in B_S\cap(A\sm C)$, $e_c\in B_S\cap(C\sm A)$ and $e_b\in B\sm B_S$. Observe that in particular $e_b\notin A\cup C$. We use \autoref{lem:crude_two-torsoing} with $f_1=a,f_2=c$ to replace the branch at $r$ in $M$ by a diatomic side $B_r'$ (here we don't care at all what happens to the colour $b$ in $B_r'$). For each $f\in\{a,b,c\}$ where $e_f$ is a torso edge (and this is always the case for $e_b$), we use \autoref{lem:crude_one-torsoing} with $f$ to replace $e_f$ with a set $B_f'$. If $e_a$ or $e_c$ is not a torso edge, set $B_a':=\{e_a\}$ or $B_c':=\{e_c\}$, respectively. Finally we contract all edges of $M$ not in $B_r'\cup B_a'\cup B_b'\cup B_c'$. We observe that each set we created has a colour unique to it by colour-completion. We observe no series class in $B_b'$ is spanned by $b$. Hence we obtain a radiolarian obstruction.

    \textit{Subcase 1b.} We have $B\neq E(W-r)$. Then there exists edges $e_a,e_b,e_c$ satisfying: $e_a\in B\cap(A\sm C)$, $e_c\in B\cap(C\sm A)$ and $e_0\in E(W-r)\setminus (A\cup B\cup C)$. If they are torso edges, we use our tools to replace the branches at $e_a$ and $e_c$ by diatomic sides $B_a'$ and $B_c'$ and the branch at $e_0$ by a colourless edge, call it also $e_0$.
    Then we use \autoref{lem:crude_three-torsoing} to replace the branch at $r$ in $M$ by either: (1) a diatomic side $B_r'$ containing $a,b,c$ but not spanned by $b$; or (2) three edges in series, each coloured exactly one of $a,b,c$; or (3) something $P$-shaped. In the first outcome, we obtain our unrealisable radiolarian. In the second outcome, it is possible that we instead obtain one of the unbalanced hexagon obstructions. In that last outcome, observe that we can always contract and parallel merge and obtain an unrealisable radiolarian.

    \textbf{Case 2.} We have $b\le_r a$ and $b\le_r c$. Then $B\subsetneq A$ and $B\subsetneq C$. Pick $e_a\in A\setminus (C\cup B)$, $e_c\in C\setminus (A\cup B)$ and $e_b\in A\cap B\cap C$. In this case, we use \autoref{lem:crude_three-torsoing} to replace the branch at $r$ by something almost diatomic, that is not spanned by either $a$ or $c$. Then our reduction proceeds similarly to in the first case.
    
    \textbf{More than two classes.} The relation $\sim_r$ is an equivalence relation but it has more than two equivalence classes. Then there are crossing colours $a,b,c$ pairwise incomparable under $\sim_r$. Let $A,B,C$ (likewise $A_S,B_S,C_S$) denote the sets of edges in $E(W-r)$ receiving colours $a,b,c$ (likewise, strongly), respectively. Without loss of generality, there are three subcases.

    \textbf{Case 1.} We have $B=E(W-r)$ and $A,C\neq E(W-r)$. Then $A$ and $C$ are inclusion incomparable and $A,C\nsubseteq B_S$. By \autoref{lem:no_three_weak}, there is exactly one edge $e_b$ in $E(W-r)$ weakly coloured by $b$ and it is also coloured by $a$ and $c$. There are also edges $e_a\in B_S\cap (A\setminus C)$ and $e_c\in B_S\cap (C\setminus A)$. We proceed analogously to in Subcase 1a of the `Intransitive' section.

    \textbf{Case 2.} We have $A,C=E(W-r)$ and $B\neq E(W-r)$. Thus we have $A_S\neq C_S$ and $B\nsubseteq A_S,C_S$. Let $e_a$ and $e_c$ be the unique distinct edges in $E(W-r)$ coloured weakly by $a$ and by $c$, respectively. We have $e_a,e_c\in B_S$. There is also some other edge $e_{ac}\in (A_S\cap C_S)\setminus B$. We replace the branch at $r$ by a \autoref{lem:crude_two-torsoing}.

    \textbf{Case 3.} We have $A,B,C\neq E(W-r)$. Then $A,B,C$ are pairwise inclusion-incomparable. If there are $e_a,e_b,e_c$ that are $a,b,c$, respectively, but not the other two, then we are done by applying \autoref{lem:crude_three-torsoing} to the branch at $r$. So assume this is not the case. Without loss of generality. $B\subseteq A\cup C$. So pick $e_a,e_{ab},e_{bc},e_c$ so $e_a\in A\setminus (B\cup C)$, $e_{ab}\in (A\cap B)\setminus C$, $e_{bc}\in (B\cap C)\setminus A$ and $e_a\in C\setminus (A\cup B)$. In this case, we use \autoref{lem:crude_three-torsoing} to replace the branch at $r$ by something almost diatomic and not spanned by~$b$.
\end{pf}

Finally, we prove that disagreeable obstructions end in radiolarian obstructions.

\begin{lem}\label{lem:disagree_go_away}
    Let $(M,c)$ be a colour-complete connected constraint matroid that is disagreeable, but stable and balanced. Then $(M,c)$ contains an unrealisable diatom as a constraint minor.
\end{lem}
\begin{pf}
    Let $(A,B)$ be a disagreeable $2$-separation in $(M,c)$. Then there are two colours, say $f_1$ and $f_2$, that witness the disagreement. That is, they are grouped differently by the left and right quasicolourings. Since $(M,c)$ is stable and balanced, we can use \autoref{lem:careful_two-torsoing}, once at each side of $(A,B)$, to replace each side by diatomic $A'$ and $B'$. The resulting constraint matroid is either a $4$-circuit or an unrealisable diatom. In the $4$-circuit case, we observe that colour-completion guarantees a special colour for each of $A',B'$, which proves that the circuit is unrealisable.
\end{pf}

\section{A detailed look at circular torsoids}
\label{sec:circuits_detailed}

\begin{dfn}[Torsoid]
    Let $(M,c)$ be a connected constraint matroid which is balanced, stable and agreeable. Then \autoref{thm:matroids_tutte_analogue} gives a mixed-sum tree-decomposition of $(M,c)$ into bags which are constraint circuits, constraint diatoms and $3$-connected constraint matroids. We call these bags the \emph{torsoids} of $(M,c)$.
\end{dfn}

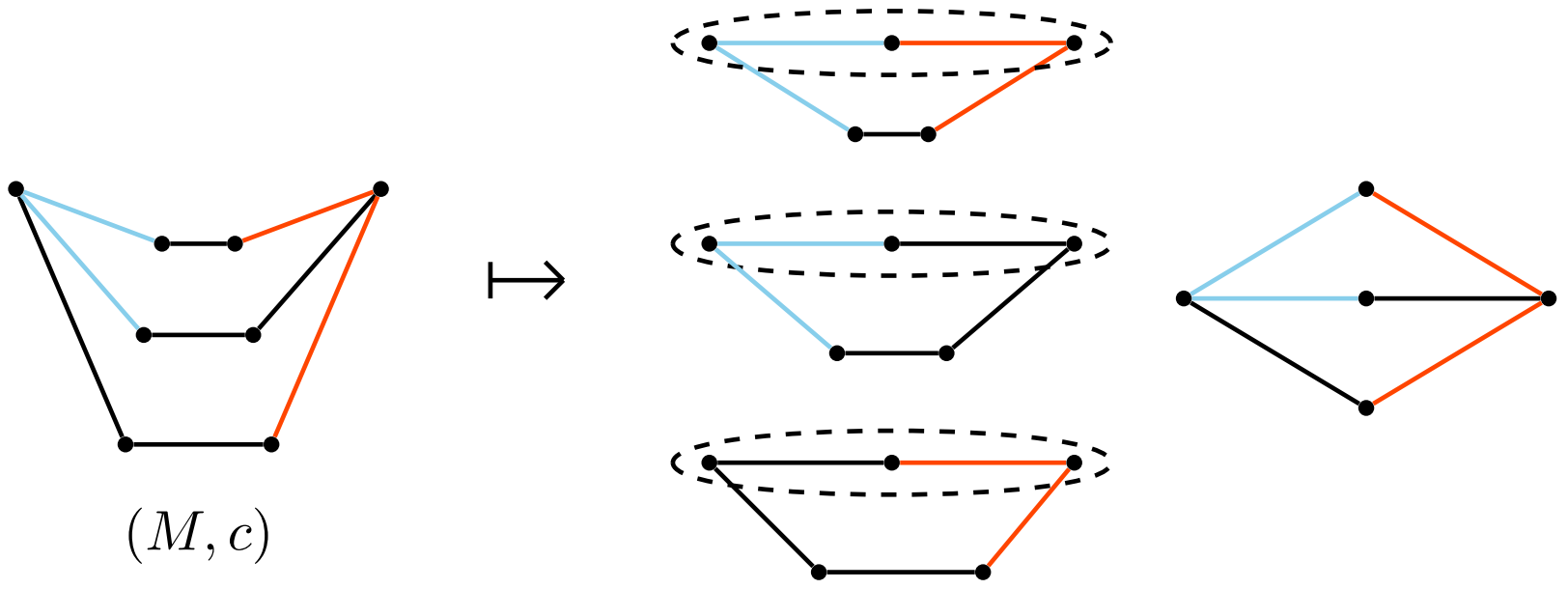
\begin{figure}
    \centering
    \begin{tikzpicture}
	\begin{pgfonlayer}{nodelayer}
		\node [style=none] (281) at (-1.5, 0.5) {};
		\node [style=none] (283) at (4.5, 0.5) {};
		\node [style=none] (284) at (-4, -0.25) {};
		\node [style=none] (285) at (-4, 0) {};
		\node [style=none] (286) at (-4, 0.25) {};
		\node [style=none] (287) at (-3, 0) {};
		\node [style=none] (288) at (-3.25, 0.25) {};
		\node [style=none] (289) at (-3.25, -0.25) {};
		\node [style=black dot] (331) at (-10.5, 1.25) {};
		\node [style=black dot] (332) at (-5.5, 1.25) {};
		\node [style=black dot] (333) at (-8.5, 0.5) {};
		\node [style=black dot] (334) at (-9, -2.25) {};
		\node [style=black dot] (335) at (-8.75, -0.75) {};
		\node [style=black dot] (336) at (-1, 3.25) {};
		\node [style=black dot] (337) at (4, 3.25) {};
		\node [style=black dot] (339) at (1.5, 3.25) {};
		\node [style=black dot] (340) at (-1, 0.5) {};
		\node [style=black dot] (341) at (4, 0.5) {};
		\node [style=black dot] (342) at (1.5, 0.5) {};
		\node [style=black dot] (344) at (-1, -2.5) {};
		\node [style=black dot] (345) at (4, -2.5) {};
		\node [style=black dot] (346) at (1.5, -2.5) {};
		\node [style=black dot] (348) at (5.5, -0.25) {};
		\node [style=black dot] (349) at (10.5, -0.25) {};
		\node [style=black dot] (350) at (8, 1.25) {};
		\node [style=black dot] (351) at (8, -1.75) {};
		\node [style=black dot] (352) at (8, -0.25) {};
		\node [style=black dot] (353) at (-7.5, 0.5) {};
		\node [style=black dot] (354) at (-7.25, -0.75) {};
		\node [style=black dot] (355) at (-7, -2.25) {};
		\node [style=black dot] (356) at (1, 2) {};
		\node [style=black dot] (357) at (2, 2) {};
		\node [style=black dot] (358) at (0.75, -1) {};
		\node [style=black dot] (359) at (2.25, -1) {};
		\node [style=black dot] (360) at (0.5, -4) {};
		\node [style=black dot] (361) at (2.75, -4) {};
		\node [style=none] (362) at (-1.5, -2.5) {};
		\node [style=none] (363) at (4.5, -2.5) {};
		\node [style=none] (364) at (-1.5, 3.25) {};
		\node [style=none] (365) at (4.5, 3.25) {};
		\node [style=none] (366) at (-8, -3.5) {$(M,c)$};
	\end{pgfonlayer}
	\begin{pgfonlayer}{edgelayer}
		\draw [style=d1] (281.center)
			 to [bend right=90, looseness=0.25] (283.center)
			 to [bend left=270, looseness=0.25] cycle;
		\draw [style=c0] (284.center) to (286.center);
		\draw [style=c0] (285.center) to (287.center);
		\draw [style=c0] (287.center) to (288.center);
		\draw [style=c0] (287.center) to (289.center);
		\draw [style=c0] (331) to (334);
		\draw [style=c1] (331) to (333);
		\draw [style=c1] (331) to (335);
		\draw [style=c1] (336) to (339);
		\draw [style=c2] (339) to (337);
		\draw [style=c1] (340) to (342);
		\draw [style=c0] (342) to (341);
		\draw [style=c0] (344) to (346);
		\draw [style=c2] (346) to (345);
		\draw [style=c0] (352) to (349);
		\draw [style=c0] (348) to (351);
		\draw [style=c1] (348) to (350);
		\draw [style=c2] (350) to (349);
		\draw [style=c2] (351) to (349);
		\draw [style=c1] (348) to (352);
		\draw [style=c0] (333) to (353);
		\draw [style=c0] (335) to (354);
		\draw [style=c0] (334) to (355);
		\draw [style=c0] (354) to (332);
		\draw [style=c2] (353) to (332);
		\draw [style=c2] (332) to (355);
		\draw [style=c0] (356) to (357);
		\draw [style=c0] (358) to (359);
		\draw [style=c0] (360) to (361);
		\draw [style=c0] (360) to (344);
		\draw [style=c0] (359) to (341);
		\draw [style=c1] (336) to (356);
		\draw [style=c2] (357) to (337);
		\draw [style=c2] (345) to (361);
		\draw [style=c1] (340) to (358);
		\draw [style=d1] (362.center)
			 to [bend right=90, looseness=0.25] (363.center)
			 to [bend left=270, looseness=0.25] cycle;
		\draw [style=d1, bend right=90, looseness=0.25] (364.center) to (365.center);
		\draw [style=d1, bend left=270, looseness=0.25] (365.center) to (364.center);
	\end{pgfonlayer}
\end{tikzpicture}
    \caption{Using \autoref{thm:matroids_tutte_analogue}, the constraint matroid $(M,c)$ (left) splits into torsoids (right). In the three circular torsoids, the colours arising from the totally-nested $2$-separations have been indicated using a dotted circle, whereas in the central torsoid, they have been omitted. Observe that the highest-positioned circular torsoid in the illustration is not a constraint minor of $(M,c)$, nor even of a colour-completion thereof.}
    \label{fig:not_minor}
\end{figure}

\begin{eg}
    A torsoid is not necessarily a constraint minor of its host. See \autoref{fig:not_minor}.
\end{eg}

In this section, we prove another characterisations of realisablity of constraint circuits (\autoref{lem:technical_scary}) and then use them to prove \autoref{lem:as_above_so_below}, which lets us translate assumptions of minimality of an obstruction constraint matroid into assumptions of minimality for one of its circular torsoids.

\subsection{Well-joinedness and gaudy colouring}
A \emph{constraint circuit with join colours} is a triple $(C,c,S)$ where $(C,c)$ is a constraint circuit and $S$ is a distinguished subset of the colours used by $c$, called \emph{join colours}, so that each colour class for a colour in $S$ is of size two and all of these colour classes are pairwise disjoint. Later, this will model the circuit having branches attached to it.

\begin{obs}\label{one_non_join}
    Let $(C,c,S)$ be a constraint circuit with join colours. If $(C,c)$ is not realisable, then there is at least one colour used by $c$ that is not in $S$. \qed
\end{obs}

Let $(C,c,S)$ be a constraint circuit with join colours. For a join colour $s\in S$, we say that a colour $f$ is \emph{total for $s$} if either both edges coloured $s$ are also coloured $f$ or both edges coloured $s$ are not coloured $f$. Otherwise we say $f$ is \emph{non-total for $s$}.

A constraint circuit with join colours $(C,c,S)$ is \emph{well-joined} if for each $s\in S$ there are at most two colours which are non-total for $s$.

\begin{dfn}
    Let $(O_Y,c_Y)$ be the constraint circuit on $6$ elements $r_1,b_1,r_2,b_2,r_3,b_3$ with colour classes $\{r_1,b_1\},\{r_2,b_2\},\{r_3,b_3\}$ and $\{r_1,r_2,r_3\}$.
\end{dfn}

\begin{obs}
    The constraint circuit $(O_Y,c_Y)$ is not realisable. \qed
\end{obs}

Given a constraint circuit $(C,c)$ and a colour $f$ used by $c$, the \emph{constraint circuit obtained from~$(C,c)$ by complementing $f$} is the constrant circuit $(C,c')$ where $c'$ is obtained from $c-f$ by adding a new colour $f^c$ which colours any edge that was not coloured by $f$ under $c$. A constraint circuit is obtained from another \emph{by colour complementation} if it is obtained by picking a subset of the colours and complementing them.

\begin{lem}\label{lem:technical_scary}
    Let $(C,c,S)$ be a constraint circuit with join colours. Then $(C,c)$ admits a realisation if and only if it does not contain constraint circuit with join colours $(C',c',S')$ so that either:\vspace{-.25cm}
    \begin{itemize}
        \item $(C',c')$ is isomorphic to $(O_Y,c_Y)$, possibly after complementing some colours; or
        \vspace{-.25cm}
        \item $(C',c',S')$ is not realisable and it is well-joined.
        \vspace{-.25cm}
    \end{itemize}
    More specifically, we obtain such a constraint minor using the following operations:\vspace{-.25cm}
    \begin{enumerate}
        \item edge-contraction outside of the colour classes of $S$; and\vspace{-.25cm}
        \item colour-forgetting outside of $S$; and\vspace{-.25cm}
        \item contracting both edges of a colour class in $S$ at once; and\vspace{-.25cm}
        \item if a colour class in $S$ is such that both edges receive identical colours, replacing it by a single edge with those colours; and\vspace{-.25cm}
        \item edge-contraction in a colour class of $S$ for which only one colour is non-total.\vspace{-.25cm}
    \end{enumerate}
\end{lem}
To prove the above, we will require some preparation.

\begin{lem}\label{PC-spectre}
    Let $(C,c)$ be a constraint circuit and $M=\lozenge(C,c)$ its overlay and suppose that no colours in $c$ colour identical sets of edges. Let $W$ be a torso of~$M$. Then union of the torso edges of $W$ with the edges $C\cap E(W)$ form a circuit in $W$. Furthermore, every other edge in $M$ is a chord of this circuit.
\end{lem}
\begin{pf}
    We will prove a stronger statment: Suppose $M$ has $k$ totally-nested $2$-separations. We will prove by induction that, for each $i\le k$, splitting $M$ into a $2$-sum $M_1\oplus_2\ldots\oplus_2M_i$ along its totally nested separations yields graphs $M_j$ with the required property. The base case is $i=0$ and $M$ has no torso edges, so we are done. Now suppose each summand in $M=M_1\oplus_2\ldots\oplus_2M_i$ has the required propery and that there is some $j$ so that $M_j$ has a totally-nested $2$-separation $(A,B)$ (otherwise $i=k$ and we are done). The outer edges from $C$ along with the torso edges in $M_j$ form a circuit $C'$. Suppose for a contradiction that $C'$ is contained entirely in $A$. Then $B$ is a single chord of $A$, but this is not a proper $2$-separation, a contradiction. We similarly get that $C'$ is not contained entirely in $B$. We conclude that each factor $M_A$ and $M_B$ has the required property.
\end{pf}
We call the circuit given by the above lemma the \emph{outer circuit} of the torso $W$.

Perhaps confusingly, in the remainder of this subsection, we will make use of a matroid with a more specialised sort of colouring. In what follows, we will assume that no colouring of a constraint matroid contains the special colours \textsc{red} and \textsc{blue}, which are reserved for the new kind of object.
\begin{dfn}[\textsc{red/blue}-circuit matroid]
    A \textsc{red/blue}-\emph{circuit matroid} is a pair $(M,C)$ where $M$ is a connected graphic matroid and $C$ is a circuit in $M$ whose elements are coloured \textsc{red} or \textsc{blue} (or both). We say that a graphic realisation $G$ of $M$ is \emph{good} if the cyclic order on the edges of~$C$ partitions the \textsc{red} and \textsc{blue} edges into segments.

    We say that $(M,C)$ is \emph{gaudy} if either one of the two following conditions holds:
    
    \noindent (\emph{$3$-overcoloured}) the circuit $C$ contains three edges, each of which recieves both \textsc{red} and \textsc{blue}; or
    
    \noindent (\emph{$4$-overcoloured}) $M$ is $3$-connected, so admits a unique graphic realisation and, in this realisation, the cyclic order of the circuit $C$ contains a suborder $(e_1,e_2,e_3,e_4)$ so that both $e_1$ and $e_3$ are coloured \textsc{red} and both $e_2$ and $e_4$ are coloured \textsc{blue}.
\end{dfn}

\begin{obs}
    A gaudy \textsc{red/blue}-circuit matroid admits no good realisations. \qed
\end{obs}

\begin{lem}\label{lem:gaudy}
    Let $(M,C)$ be a \textsc{red/blue}-circuit matroid. For each node $t$ of the Tutte decomposition tree of $M$, suppose that there is a torso $W_t$ whose torso edges, along with the edges of $C$ restricted to $W_t$, form a circuit $C_t$. Colour the torso edges \textsc{red}, \textsc{blue} or both, if such colours are in the branch corresponding to this torso edge. Then $(M,C)$ admits a good realisation if and only if none of the \textsc{red/blue}-circuit matroids $(W_t,C_t)$ are gaudy.
\end{lem}
\begin{pf}
    The proof is a direct analogue of the proof of Lemma 4.31 in \cite{me2025}. One need only note that the coloured `bond' can be substituted for the circuit coloured \textsc{red} and \textsc{blue} by appropriate dualisation.
\end{pf}

Now we are equipped to prove the main lemma of this subsection.
\begin{pf}[Proof of \autoref{lem:technical_scary}]
    Let $(C,c,S)$ be a constraint circuit with join colours which is not realisable and minimal with this property under the operations (i)-(v). We are required to show that either $(C,c)$ is isomorphic to $(O_Y,c_Y)$, possibly after complementing some colours; or $(C,c,S)$ is unrealisable and well-joined.
    
    Let $f$ be a colour used by $c$ and $c-f$ the colouring obtained from $c$ by forgetting $f$. By \autoref{one_non_join}, we may assume that $f$ is not a join colour.
    Let $M=\lozenge (C,c-f)$. By minimality and \autoref{overlay}, the matroid $M$ is graphic. Define the \textsc{red/blue}-circuit matroid $(M,C)$ by colouring the $f$-coloured edges \textsc{red} and all the other edges of $C$ \textsc{blue}.

    \begin{claim}
        The \textsc{red/blue}-circuit matroid $(M,C)$ has no good realisation.
    \end{claim}
    \begin{cproof}
        Suppose for a contradiction that $G$ is a good realisation of $(M,C)$. Then by \autoref{overlay}, the graph $G$ realises $(C,c-f)$. Furthermore, since $G$ is good, the \textsc{red} edges in $C$ form a segment in its cyclic order in $G$. Hence $G$ is also a realisation of $(C,c)$, a contradiction.
    \end{cproof}
    For each node $t$ of the Tutte decomposition of $M$, let $W_t$ be its torso and $C_t$ the outer cycle as per \autoref{PC-spectre}. In particular, the torso edges inherit any colours in their corresponding branches. By the above claim and \autoref{lem:gaudy}, there is a node $t$ in the Tutte decomposition tree so that $(W_t,C_t)$ is gaudy. There are two cases.

    \textbf{Case 1.} The \textsc{red/blue}-matroid $(W_t,C_t)$ is $3$-overcoloured, witnessed by edges $e_1$, $e_2$ and $e_3$. Since no edge of $C$ is both \textsc{red} and \textsc{blue}, we get that $e_1$, $e_2$ and $e_3$ are torso edges. Let $Z$ be the set of edges of $M$ not in a branch at $e_1,e_2$ or $e_3$.
    \begin{claim}
        The set $Z$ is empty.
    \end{claim}
    \begin{cproof}
        Suppose $Z$ is nonempty and consider the \textsc{red/blue}-circuit matroid obtained from $(M,C)/Z$ by contracting $Z$. Then $(M,C)/Z$ has a triangle torso, around which the branches at $e_1,e_2,e_3$ are arranged. This triangle torso is gaudy, and hence by \autoref{lem:gaudy} $(M,C)/Z$ has no good realisation. It remains to show that we can `model' this contraction in the constraint circuit $(C,c)$. That is, contracting $Z\cap E(C)$ in $(C,c)$ can be achieved by the operations (i)-(v). Indeed, since each of the torso edges $e_1$, $e_2$ and $e_3$ correspond to a totally-nested $2$-separation, no colour from $c$ colours edges from $Z\cap E(C)$ as well as ones in the branches at $e_1,e_2$ and $e_3$, except possibly \textsc{red}. Hence we can use operations (i), (iii).
    \end{cproof}
    
    By the above claim, $W_t$ consists of only the edges $e_1$, $e_2$ and $e_3$. For each $i\in\{1,2,3\}$, let $C_i$ be the segment of $C_i$ in the branch at $e_i$.

    \begin{claim}
        For each $i\in \{1,2,3\}$, there is an edge $f_i$ of $M$ in the branch at $e_i$ so the subsegment of $C_i$ defined by $f_i$ contains a purely \textsc{red} edge $r_i$ and a purely \textsc{blue} edge $b_i$.
    \end{claim}
    \begin{cproof}
        Suppose for a contradiction that this is not the case. Then we can bipartition the chords of $C$ in the branch at $e_i$ into those whose subsegments of $C_i$ contain only \textsc{red} and those which contain only \textsc{blue}. Let these sets, along with their subsegments of $C$ be written $R$ and~$B$, respectively. Then there is a $2$-separation with one side $R$ and the branch at $e_{i+1}$ and the other side $B$ and the branch at $e_{i-1}$. This $2$-separation crosses the $2$-separation corresponding to $e_i$, a contradiction to the fact that it is totally-nested.
    \end{cproof}

    Let $(O',c')$ be the constraint circuit obtained from $(C,c)$ by contracting, for each $i\in\{1,2,3\}$, all of $C_i$ except for $r_i$ and $b_i$. Observe that the constraint matroid we obtain is isomorphic to $(O_Y,c_Y)$, up to colour-complementation.
    We claim that such a constraint minor can be found using only the operations (i)-(v). Let $i\in\{1,2,3\}$. First suppose that there is a join colour which colours both $r_i$ and $b_i$. In this case, we can contract everything in $C_i$ outside of these two edges using (i) and (iii). So assume that there is no join colour class in $C_i$, ones of whose edges is \textsc{red} and one of whose edges in \textsc{blue}. Since every edge is either \textsc{red} or \textsc{blue}, there is a join colour class in $C_i$, both whose edges are \textsc{red} and a join colour class in $C_i$, both whose edges are \textsc{blue}. By forgetting any other overlapping colours which are non total for one of these, we can then use (iv).

    \textbf{Case 2.} The \textsc{red/blue}-matroid $(W_t,C_t)$ is $4$-overcoloured, witnessed by a quadruple $(e_1,e_2,e_3,e_4)$ where both $e_1$ and $e_3$ are \textsc{red}, both $e_2$ and $e_4$ are coloured \textsc{blue} and the cyclic order in which the edges $(e_1,e_2,e_3,e_4)$ appear in the outer cycle of $W_t$ in its unique graphic realisation is as written.

    \begin{claim}
        The torso $W_t$ has no torso edges.
    \end{claim}
    \begin{cproof}
        Suppose for a contradiction that there is a torso edge $e$ in $W_t$. Then there are at least two edges in the branch at $e$ in $M$ that come from $C$. Suppose first that $e\notin\{e_1,e_2,e_3,e_4\}$. Then using similar reasoning in Case 1, we obtain a contradiction to minimality by contracting all the edges of $C$ in the branch at $e$ except for one. To do so, we might need to use operation (v).

        So suppose without loss of generality that either $e=e_1$ or $e=e_2$. We claim that there exists an edge in $C$ in the branch at $e$ in $M$ which we can contract and contradict minimality. For this to work, we notice that in the case $e=e_1$ at least one \textsc{red} edge must remain and in the case $e=e_2$ one \textsc{blue} edge must remain. To achieve this, we observe that we can apply colour-forgetting inside the branch until it is possible to apply (v), reminding the reader that \textsc{blue} is an auxiliary colour which is not used by $c$.
    \end{cproof}

    By the above claim, $W$ has no torso edges. That means, $W=M$ and it is $3$-connected. Let $G$ be the unique graphic realisation of $W$. Let $P$ be the subpath of the realisation of $C$ in $G$ that is obtained from $C$ by deleting the edge $e_4$. We assume that $e_1$ and $e_3$ are as close together as possible on $P$. That is, all of the edges in the maximal subpath $Q\subseteq P$ between (and excluding) $e_1$ and $e_3$ are~\textsc{blue}. Observe that each chord $f$ of $C$ in $M$ defines a subpath of $P$, which we denote by $P(f)$.
    \begin{claim}
        There exists a sequence of chords $(f_1,\ldots,f_\ell)$ of $C$ so that:\vspace{-.25cm}
        \begin{enumerate}
            \item we have $e_1\in P(f_1)$ and $e_3\notin P(f_1)$; and\vspace{-.25cm}
            \item we have $e_1\notin P(f_\ell)$ and $e_3\in P(f_\ell)$; and \vspace{-.25cm}
            \item for $i\in[\ell-1]$, $P(f_i)$ and $P(f_{i+1})$ have nonempty intersection; and \vspace{-.25cm}
            \item for distinct $i,j,k\in[\ell]$, no edge in $C$ is in all three of $P(f_i),P(f_j)$ and $P(f_k)$.\vspace{-.2cm}
        \end{enumerate}
    \end{claim}
    \begin{cproof}
        Let $f_1$ be chosen satisfying (i) and maximising the number of edges of $Q$ that are contained in~$P(f_1)$. Let $f_{-1}$ be chosen satisfying (ii) and maximising the number of edges of $Q$ that are contained in~$P(f_{-1})$.
        
        Since otherwise we would be done, suppose that $P(f_1)$ and $P(f_{-1})$ do not overlap.
        Let $b_1$ and $b_{-1}$ be the last and first edges in $P(f_1)$ and $P(f_{-1})$, respectively. We repeat the process, with $b_1$ and $b_{-1}$ in place of $e_1$ and $e_3$. We obtain $f_2$ and, if necessary, $f_{-2}$ Observe that $e_1\notin P(f_2)$ and $e_3\notin P(f_{-2})$. Continuing in this way completes the proof of our claim.
    \end{cproof}

    Since otherwise it would contradict minimality, we know that the only colours present in $(C,c,S)$ are the colours $f_1,\ldots,f_\ell$, \textsc{red} and perhaps some additional join colours not in that list. It remains to prove that $(C,c,S)$ is well-joined. Let $s\in S$ be an arbitrary join colour and consider its colour class $c_s$ in $C$, which is connected in the graphic realisation. Since $s$ is a join-colour, it $c_s$ overlaps no other join colour. By construction in the above claim, $c_s$ is not partially overlapped by more than two colours from $f_1,\ldots,f_\ell$. Since the only \textsc{red} edges in $P$ are outside of the subpath between $e_1$ and $e_3$, we also observe that it is not possible for an edge of $C$ to lie in two segments from $f_1,\ldots,f_n$ and be coloured \textsc{red}. This completes our proof.
\end{pf}

\subsection{Microcosm, macrocosm}
\begin{lem}\label{lem:as_above_so_below}
    Let $(M,c)$ be a colour-complete connected constraint matroid which is stable, balanced and agreeable. Let $W$ be a circular torso of $M$ and $(W_*,c_*)$ its corresponding torsoid. Let $S$ be the set of colours used by $c_*$ that complete totally-nested $2$-separations left-incident with $W$ in $M$. Let $(W_*',c_*')$ be a constraint minor of $(W_*,c_*)$ obtained via the operations (i)-(v) from \autoref{lem:technical_scary}, where $S$ is the set of join-colours. Then there is a constraint minor of $(M,c)$ with $(W_*',c_*')$ as a torsoid.
\end{lem}
\begin{pf}
    We will prove this lemma in the case that $(W_*',c_*')$ is obtained from $(W_*,c_*)$ by one use of one of the operations (i)-(v) from \autoref{lem:technical_scary}, and split into cases based on the type of operation.

    \textbf{Case (i).} $(W_*',c_*')$ is obtained from $(W_*,c_*)$ by contracting an edge outside the colour classes of $S$. In this case, the edge is a non-torso edge of $W$ and can be contracted in $M$.

    \textbf{Case (ii).} $(W_*',c_*')$ is obtained from $(W_*,c_*)$ by forgetting a colour outside of $s$. In this case, we simply forget the colour in $(M,c)$.

    \textbf{Case (iii).} $(W_*',c_*')$ is obtained from $(W_*,c_*)$ by contracting both edges of a colour class $s\in S$ at once. The colour class of $s$ corresponds to a torso edge $e$ of $W$. In this case, we simply contract everything in the branch at $e$ in $M$.

    \textbf{Case (iv).} $(W_*',c_*')$ is obtained from $(W_*,c_*)$ by replacing a colour class $s\in S$, both whose edges are identically coloured, by a single edge with those colours. The colour class of $s$ corresponds to a torso edge $e$ of $W$. We find a bond in the branch at $e$ in $M$ and contract everything outside the bond and then parallel merge everything from the bond. The resulting single edge contains all the colours from the two edges in $W_*$, since by the construction in \autoref{thm:matroids_tutte_analogue} all such colours span the $2$-separation in the branch side.

    \textbf{Case (v).} $(W_*',c_*')$ is obtained from $(W_*,c_*)$ by conracting one of the edge in a colour class of $s\in S$, for which only one other colour $f$ is non-total for it. The colour class of $s$ corresponds to a torso edge $e$ of $W$. Suppose first that the contracted edge in $W_*$ has colour $f$. Then we can use \autoref{lem:cut_control} to replace the branch at $e$ in $M$ by a single edge not coloured by $f$ but having all the other colours of the class of $s$. Otherwise, the edge $e$ is the one not coloured $f$. In this case, we use \autoref{lem:cut_grab}.
\end{pf}

\begin{lem}\label{lem:find_now_we_have_well}
    Let $(M,c)$ be a colour-complete connected constraint matroid which is stable, balanced and agreeable. Let $W$ be a circular torso of $M$ and $(W_*,c_*)$ its corresponding torsoid. Let $S$ be the set of colours used by $c_*$ that complete totally-nested $2$-separations left-incident with $W$ in $M$. Suppose that $c$ is well-joined with respect to $S$ as the set of join colours. Then $(M,c)$ contains a radiolarian $(W_{**},c_{**})$ as a constraint minor. Furthermore, if $(W_*,c_*)$ is not realisable, then $(W_{**},c_{**})$ is not realiable.
\end{lem}
\begin{pf}
    For each join colour $s\in S$, there is a corresponding torso edge $e_s$ in $W$. Since $(W_*,c_*)$ is well-joined with respect to $S$, each $e_s$ has at most two weak colours in the torso colouring $c_W$. Let $(W_{**},c_{**})$ denote the constraint matroid obtained from $W$ by using \autoref{lem:crude_one-torsoing} and \autoref{lem:crude_two-torsoing} to replace each $e_s$ by a diatomic side $B_s$. Since using \autoref{lem:crude_one-torsoing} and \autoref{lem:crude_two-torsoing}, we have no control over whether the two edges we get respect the quasicolouring at $e_s$, and hence whether $(W_{**},c_{**})=(W_*,c_*)$. However, we claim that, if $(W_*,c_*)$ is not realisable, then $(W_{**},c_{**})$ is not realisable.

    If each $B_s$ has more than one series class, we observe that $(W,c_W)$ is the central coloured torso of the constraint radiolarian $(W_{**},c_{**})$: in fact, in the decomposition of \autoref{thm:matroids_tutte_analogue} it corresponds to a torsoid which is exactly $(W_*,c_*)$. This finishes the proof. Otherwise, not every $B_s$ has more than one series class. But in that case, we can artificially add in a new edge that spans each $B_s$ which is a single series class and this does not impact the realisability of the constraint radiolarian, or the fact that one of its torsoids is exactly $(W_*,c_*)$.
\end{pf}

\begin{cor}\label{lem:radio_resolve}
    Let $(M,c)$ be a colour-complete connected constraint matroid which is stable, balanced and agreeable. Let $W$ be a circular torso of $M$ and $(W_*,c_*)$ its corresponding torsoid. If $(W_*,c_*)$ is not realisable, then $(M,c)$ contains a radiolarian obstruction as a constraint minor.
\end{cor}
\begin{pf}
    By \autoref{lem:as_above_so_below} with \autoref{lem:technical_scary}, we may assume $(W_*,c_*)$ is well-joined with respect to $S$. Now we can apply \autoref{lem:find_now_we_have_well} and this completes our proof.
\end{pf}

\section{A detailed look at diatom torsoids}\label{sec:diatoms_detailed}

\begin{lem}\label{lem:diatom_min}
    Let $(D,c)$ be a clean constraint diatom which is unrealisable. If $(D,c)$ is a minimal under colour-forgetting, then it has no series classes of size two containing more than two colours.
\end{lem}
\begin{pf}
    Since $(D,c)$ is not realisable, by \autoref{lem:quick_check_diatom} we get that its frustration graph has an unbalanced cycle. We require the following elementary claim.

    \begin{claim}
        Let $O$ be a minimum-length unbalanced cycle in a signed graph. Then there is no path $P$ which meets $O$ exactly in its end vertices so that the union $O\cup P$ contains a smaller cycle than $O$.
    \end{claim}
    \begin{cproof}
        Suppose for a contradiction that $P$ such a path. Then the endvertices of $P$ split $O$ into two paths, one with an even number of negative edges and one with an odd number. Nomatter how many negative edges $P$ contains, one of the cycles using $P$ and a path from $O$ is unbalanced, contradicting the minimality of $O$.
    \end{cproof}

    Let $O$ be an unbalanced cycle in the frustration graph of $(D,c)$ of minimal length. By the above claim, $O$ has no chords. Write $v_1:=\{a_1,b_1\},\ldots,v_\ell:=\{a_\ell,b_\ell\}$ for the series classes which are the vertices of $O$, in that order. Since $O$ has no chords, there is a unique colour $f_i$ for each $i\in[\ell]$ which is shared only between $\{a_i,b_i\}$ and $\{a_{i+1},b_{i+1}\}$. Since otherwise it would not be minimal under colour-forgetting, we may assume that $f_1,\ldots,f_\ell$ are the only colours used by $c$.

    Suppose for a contradiction that there is a series class $v=\{a_0,b_0\}$ not in $O$ that contains at least three distinct colours from $f_1,\ldots,f_\ell$, say  $f_i,f_j,f_k$. Then $v$ forms a triangle with the edge $e_q:=v_qv_{q+1}$ for each $q\in\{i,j,k\}$. 

    For each $q\in\{i,j,k\}$, let $O_{qq+1}$ be the maximal path in $O$ meeting $e_q$ and $e_{q+1}$ but not containing $e_{q}$, $e_{q+1}$ or $e_{q+2}$.
    For distinct $p,q\in\{i,j,k\}$, let $e_{pq}$ be the edge from $v$ to the vertex incident with both $e_i$ and $O_{pq}$ (or $O_{qp}$). Let $C_q$ be the cycle formed by $e_q$ and $e_{qq+1},e_{qq-1}$. Let $C_{pq}$ be the cycle formed by $O_{pq}$ and $e_{pq}$ and $e_{qp}$.

    Let $\sigma_1,\sigma_2,\ldots,\sigma_6$ be the signs of the segments $e_i$, $O_{ij}$, $e_j$, $O_{jk}$, $e_k$, $O_{ki}$, respectively. Let $\gamma_1\ldots\gamma_6$ be the signs of the edges $e_{ij}$, $e_{ji}$, $e_{jk}$, $e_{kj}$, $e_{ki}$, $e_{ik}$, respectively.

    \begin{claim}
        We have $\gamma_{i+1}=\sigma_{i+1}\cdot\gamma_i$.
    \end{claim}
    \begin{cproof}
        The claim follows by the minimality of $O$.
    \end{cproof}

    The above claim proves that the number of negative $\sigma$s must be even. But then $O$ is not unbalanced, a contradiction.
\end{pf}

\begin{lem}\label{lem:further_cleaning}
    Let $(M,c)$ be an unrealisable constraint matroid that is minimal under constraint minors and is stable and balanced. Suppose that $W$ is a bond torso of $M$ whose corresponding torsoid is an unrealisable constraint diatom. Then the totally-nested $2$-separations incident with $W$ have at most two crossing non-spanning colours each in $c$.
\end{lem}
\begin{pf}
    Let $(D,c')$ be the unrealisable constraint diatom torsoid.
    We know that by \autoref{lem:quick_check_diatom} that the frustration graph of the clean version of $(D,c)$ is unbalanced and by \autoref{lem:diatom_min} that each series class contains as most $2$ non-spanning colours. By colour forgetting, we are done.
\end{pf}

\begin{lem}\label{lem:deal_with_diatoms}
    Let $(M,c)$ be a connected constraint matroid that is stable and balanced. Suppose it has a bond torso $W$ whose corresponding diatom torsoid $(W^*,c^*)$ is not realisable. Suppose not more than two crossing colours for each torso edge. Then either it contains a radiolarian obstruction or a diatom obstruction as a constraint minor.
\end{lem}
\begin{pf}
    By \autoref{lem:further_cleaning}, we may assume that the series classes of the torsoid $(W^*,c^*)$ of length two contain at most two colours each. Now by \autoref{lem:careful_two-torsoing} applied at each series class which corresponds to a branch in $(M,c)$, we obtain either a radiolarian obstruction as a constraint minor, or build a diatom obstruction. This completes our proof.
\end{pf}

\section{Proof of \autoref{thm:forb}}\label{sec:done}

\begin{lem}\label{lem:strict_reality}
    Let $(M,c)$ be a connected constraint matroid that uses one colour and $W$ a $3$-connected torso of $M$ and $c'$ its strict torso colouring. Then the constraint matroid $(W,c')$ is a constraint minor of $(M,c)$.
\end{lem}
\begin{pf}
    Let $f$ be the single colour used in $c$. We use \autoref{lem:cut_grab} and \autoref{lem:cut_control} on the branches at each torso edge of $(W,c')$ to replace the branches by single edges coloured by $f$, in the case where the torso edge was strongly coloured by $f$, and not coloured at all, in the case where the torso edge was weakly coloured by $f$.
\end{pf}

\begin{cor}\label{cor:3-con_resolution}
    Let $(M,c)$ be a connected constraint matroid that is stable, balanced and agreeable. If the mixed-sum tree decomposition of \autoref{thm:matroids_tutte_analogue} has a $3$-connected bag which is not realisable, then $(M,c)$ contains as a constraint minor one of the $3$-connected obstructions.
\end{cor}
\begin{pf}
    Recall that the mixed-sum tree decomposition has as $3$-connected torsoids the $3$-connected torsos of the Tutte decomposition along with the strict torso colouring.
\end{pf}

Now we are finally ready to prove \autoref{thm:forb}.

\begin{pf}[Proof of \autoref{thm:forb}]
    Let $(M,c)$ be a connected constraint matroid which is colour-complete. We are required to show that $(M,c)$ is realisable if and only if it does not contain as a constraint minor either:\vspace{-.25cm}
    \begin{enumerate}
        \item one of the $3$-connected obstructions depicted in \autoref{fig:obst3con}; or\vspace{-.25cm}
        \item an unrealisable diatom (definined \autoref{sec:diatoms}); or\vspace{-.25cm}
        \item an unrealisable radiolarian (defined \autoref{sec:radiolarians}); or\vspace{-.25cm}
        \item one of the unbalanced hexagon obstructions, depicted in \autoref{fig:unbalanced_hexagon}.\vspace{-.25cm}
    \end{enumerate}
    Suppose that $(M,c)$ is not realisable. We will show that either $(M,c)$ contains of the above obstructions as a constraint minor, or it contains as a constraint minor a connected constraint matroid which is colour complete and also not realisable. By induction, this suffices.
    
    First suppose that $(M,c)$ contains one of the non-local obstructions. That is, it is unstable, unbalanced or disagreeable. If it is unstable, then we obtain our smaller colour-complete obstruction by \autoref{lem:unstable_never_happens}. So suppose $(M,c)$ is stable, but that it is unbalanced. Then by \autoref{lem:unbalanced_resolved}, we find either a radiolarian obstruction or an unbalanced hexagon as a constraint minor. So suppose that $(M,c)$ is balanced and stable. If $(M,c)$ is disagreeable, then we are done by \autoref{lem:disagree_go_away}.

    So we may assume that $(M,c)$ is balanced, stable and agreeable. In this case, we can apply \autoref{thm:matroids_tutte_analogue} to get that $(M,c)$ has an unrealisable torsoid $(W^*,c_W^*)$. There are three cases: this torsoid is $3$-connected, a diatom, or a circuit. In the first case, we find a $3$-connected obstruction as a constraint minor by \autoref{cor:3-con_resolution}. In the second case, we obtain a radiolarian obstruction as a constraint minor by \autoref{lem:radio_resolve}. In the last case, we obtain a diatom obstruction as a constraint minor by \autoref{lem:deal_with_diatoms}. We are done.
\end{pf}

\section{Disconnected constraint matroids}
\label{sec:disconnected}
\subsection{Overview}

The main result of this section is the following.

\begin{thm}\label{thm:new_red_to_con}
    Let $(M,c)$ be a constraint matroid. For each component $B$ of $M$, there is an extension $c_B^*$ of the restriction of $c$ to $B$ so that\vspace{-.25cm}
    \begin{enumerate}
        \item every realisation of $(M,c)$ induces for each component $B$ of $M$ a realisation for $(B,c_B^*)$; and\vspace{-.25cm}
        \item if $(M,c)$ admits a realisation, then given for each component $B$ of $M$ a realisation $G_B$ of $(B,c_B^*)$, there is a realisation of $(M,c)$ which induces the graphs $G_B$.\vspace{-.25cm}
    \end{enumerate}
\end{thm}

In \autoref{sec:spaces}, we will develop the theory of connected constraint matroids to understand their spaces of possible realisations. In \autoref{sec:pointy} and \autoref{subsec:quasibranches}, we will introduce new colouring tools and ways of decomposing disconnected constraint matroids, respectively. Finally, in \autoref{subsec:proof2}, we will prove \autoref{thm:new_red_to_con}.

\subsection{Spaces of realisations}
\label{sec:spaces}
In this subsection, we will develop the theory of realisations of connected constraint matroids slightly further.
Let $G$ be a $2$-connected graph. Suppose $(A,B)$ is a $2$-separation of $G$. There are exactly two ways to recover a graph by gluing the sides $G[[A]]$ and $G[[B]]$. One yields $G$ and the other yields a new graph $H$ which we call the \emph{twist of $G$ about $(A,B)$}. The graph $H$ is also a realisation of $M(G)$, the cycle matroid of $G$.
\begin{lem}[Whitney's $2$-isomorphism Theorem 1933 \cite{whitney19332}]\label{thm:2-iso}
    Let $M$ be a connected matroid. Then any realisation of $M$ can be obtained from any other via a sequence of twists.
\end{lem}
Given a connected constraint matroid, it is of course possible using the above lemma to move between realisations using twists. However, it is no a priori obvious that at each step we have a graph which actually realises our constraint matroid. However, that is indeed the case.
\begin{lem}\label{lem:homotwisty}
    Let $(M,c)$ be a connected constraint matroid. Then any realisation of $(M,c)$ can be obtained from any other via a sequence of twists, so that each intermediate graph realises $(M,c)$.
\end{lem}
For the proof of the above, we will need the following three lemmata.
\begin{lem}\label{lem:homotwisty_circuit}
    Let $(C,c)$ be a constraint circuit. Then any realisation of $(C,c)$ can be obtained from any other via a sequence of twists, so that each intermediate graph realises $(C,c)$.
\end{lem}
\begin{pf}
    Let $M=\lozenge(C,c)$ denote the overlay of $(C,c)$. Then $M$ is a connected matroid and by \autoref{overlay}, the realisations of $(C,c)$ are exactly the restrictions to $E(C)$ of the realisations of $M$. By \autoref{thm:2-iso}, each pair of realisations of $M$ are connected by a sequence of twists. Our result follows from the observation that a twist of a realisation of $M$ corresponds to a twist of a realisation of $(C,c)$.
\end{pf}
\begin{lem}\label{lem:homotwisty_diatom}
    Let $(D,c)$ be a constraint diatom. Then any realisation of $(D,c)$ can be obtained from any other via a single twist.
\end{lem}
\begin{pf}
    Let $\{e_1,f_1\},\ldots,\{e_n,f_n\}$ denote the series classes with two edges in $D$. Let $G$ and $G'$ be realisations of $(D,c)$ and without loss of generality assume that $e_i$ meets $x$ and $f_i$ meets $y$ in $G$. Let $I\subseteq [n]$ be the indices so that $e_i$ does not meet $x$ in $G'$. Then $\cup_{i\in I}\{e_i,f_i\}$ forms one side of a totally-nested $2$-separation of $D$. Twisting by this separation transforms $G$ into $G'$ and vice versa.
\end{pf}

\begin{lem}\label{lem:homotwisty_extended_3-block}
    Let $(M,c)$ be a realisable extended $3$-block. Then each realisation can be obtained from any other via a sequence of twists, so that each intermediate graph realises $(M,c)$.
\end{lem}
\begin{pf}
    For the proof of this lemma, we observe that every realisation of $(M,c)$ is obtained from the unique graphic realisation of its central torso by taking a $2$-sum with a triangle at each torso edge. There are exactly two choices for how to perform each $2$-sum, and clearly they are all equivalent up to twisting at the totally-nested $2$-separator corresponding to the torso edge.
\end{pf}

\begin{pf}[Proof of \autoref{lem:homotwisty}]
    Let $(M,c)$ be a connected constraint matroid. We are required to show that any realisation of $(M,c)$ can be obtained from any other via a sequence of twists, so that each intermediate graph realises $(M,c)$.

    Let $(T,\{(M_u,c_u)\},\{\Ecal_u^e\})$ be the mixed-sum tree decomposition obtained in \autoref{thm:matroids_tutte_analogue}. Given a representation $G_u$ for each $(M_u,c_u)$, we obtain a representation $G$ for $(M,c)$ which induces each~$G_u$. This $G$ is unique up to a twist at each totally-nested $2$-separation of $M$ which has no weakly-crossing colours, along with a choice of how to realise each extended $3$-block from a realisation of its strictly-coloured central torso.
    
    Our key observation is that a twist of a representation $G_u$ (or a choice of realisation of an extended $3$-block) corresponds in $G$ to a twist at a (not necessarily totally-nested) $2$-separation.
    
    By \autoref{lem:homotwisty_circuit} and \autoref{lem:homotwisty_diatom} there is a sequence of twists for each $(M_u,c_u)$ so that any rep can be reached from any other via a sequence of twists so that each intermediate graph is a rep. Combining this fact with \autoref{lem:homotwisty_extended_3-block} completes our proof.
\end{pf}

\subsection{Pointy constraint matroids}\label{sec:pointy}
We introduce the following, which slightly generalises constraint matroids and their realisations.
\begin{dfn}
    A \emph{pointy constraint matroid} $(M,c,\Dcal)$ is a constraint matroid $(M,c)$ along with a family $\Dcal=\{D_1,\ldots,D_k\}$ of subsets of colours from $c$. A \emph{realisation} of a pointy constraint matroid $(M,c,\Dcal)$ is a realisation $G$ of $(M,c)$ so that for each $i\in[k]$, there is a vertex of $G$ seeing the colours~$D_i$.
\end{dfn}

Given a matroid $M$, we say that a colouring $c'$ of $M$ \emph{extends} a colouring $c$ of $M$ if for every edge $e$ in $M$, the colours assigned to $e$ by $c$ is a subset of the colours assigned to $e$ by $c'$. In this case, we get that every realisation for $(M,c')$ is also a realisation for $(M,c)$.

\begin{dfn}\label{dfn:using_D}
Let $(M,c)$ be a constraint matroid and $\Dcal$ a family of subsets of the colours used by $c$. Let $c'$ be the extension of the colouring $c$ defined as follows: For each pair of distinct colours $p$ and $q$ in some shared $D\in\Dcal$, create a new colour $p*q$ which colours every edge receiving colour $p$ or receiving colour $q$. Call $c'$ the \emph{extension of $c$ using $\Dcal$}.
\end{dfn}

\begin{lem}\label{cor:D_works_MORE}
    Let $(M,c,\Dcal)$ be a pointy constraint matroid and $c'$ extension of $c$ using $\Dcal$. If $(M,c,\Dcal)$ admits a realisation, then the set of realisations of the connected constraint matroid $(M,c')$ is identical to the set of realisations of $(M,c,\Dcal)$.
\end{lem}

To simplify matters, we will prove the above for the case where $\Dcal$ is a singleton \{D\}. We will simplify notation by writing $D$ instead of $\{D\}$ in some places.
\begin{lem}\label{lem:using_D_works}
    Let $(M,c,D)$ be a pointy constraint matroid and $c'$ extension of $c$ using $D$. If $(M,c,D)$ admits a realisation, then the set of realisations of the connected constraint matroid $(M,c')$ is identical to the set of realisations of $(M,c,D)$.
\end{lem}
\begin{pf}[Proof of \autoref{cor:D_works_MORE} assuming \autoref{lem:using_D_works}]
    Given $D\in\Dcal$ observe that a realisation of $(M,c,\Dcal)$ is always also a realisation of $(M,c,\{D\})$ and that the extension of the $c$ using $\{D\}$ is extended by the extension of $c$ using $\Dcal$. Now we simply apply \autoref{lem:using_D_works} for each $D\in\Dcal$.
\end{pf}

The remainder of the subsection will be used to prove \autoref{lem:using_D_works}.
\begin{eg}
    The assumption in \autoref{lem:using_D_works} that $(M,c,D)$ must have a realisation is necessary. Indeed, consider the constraint circuit $(C,c)$ having three elements $e_1,e_2,e_3$ coloured $f_1,f_2,f_3$, respectively. Consider also the pointy constraint circuit $(C,c,D:=\{f_1,f_2,f_3\})$. Clearly $(C,c)$ has a realisation but $(C,c,D)$ does not. In spite of this, let $c'$ be the colouring obtained from $(C,c,D)$ as in \autoref{lem:using_D_works}. We observe that $(M,c')$ is realisable.
\end{eg}

We first prove the connected case of \autoref{lem:using_D_works}.
\begin{lem}\label{lem:red_to_con_from_pty_con}
    Let $(M,c,D)$ be a pointy connected constraint matroid and $c'$ be the extension of $c$ using $D$. If $(M,c,D)$ admits a realisation, then the set of realisations of the connected constraint matroid $(M,c')$ is identical to the set of realisation of $(M,c,D)$.
\end{lem}
\begin{pf}
    Firstly let $G$ be an arbitrary realisation of $(M,c,D)$. Let $u$ be a vertex in $G$ which sees all the colours $D$. Now for each pair $p,q\in D$, the connected subgraph of $G$ whose edges are coloured $p$ contains $u$ and the connected subgraph of $G$ whose edges are coloured $q$ meets $u$. Hence the subgraph induced on the edges coloured $p*q$ is connected. This proves that $G$ is a realisation for~$(M,c')$.

    Now let $H$ be an arbitrary realisation of $(M,c')$. By \autoref{lem:homotwisty}, $G$ can be obtained from $H$ by a series of twists, so that each intermediate graph is a realisation of $(M,c')$. We claim that $H$ is a realisation of $(M,c,D)$. By induction on the number of twists, we assume that $H$ is obtained from $G$ by a single twist about a $2$-separation $(A,B)$ of $G$. Let $\{x,y\}$ be its separator in $G$. Let $\{x,y\}$ also denote the vertices in $H$ that separate $A$ and $B$ so that the edges in $A$ that are incident with $x$ and $y$, respectively, are the same in $G$ and $H$. By assumption, there is a vertex $u$ in $G$ that that sees all of $D$. If $u$ is not in the separator $\{x,y\}$, then clearly $H$ also has a vertex seeing all of $D$, and we would be done. Hence we assume without loss of generality that $u=x$. Let $X_A$ and $X_B$ be the sets of colours in $D$ so that there is an edge of that colour incident with $x$ in $A$ and $B$, respectively. Likewise, define $Y_A$ and $Y_B$ for the vertex~$y$.

    Suppose for a contradiction that $H$ is not a realisation of the pointy constraint matroid $(M,c,D)$. We will use the following two claims.

\begin{claim}\label{clm:bothy}
  Let $f\in D$ so that $f$ colours both some edge in $A$ and some edge in $B$. Then $(A,B)$ is $f$-spanned.
\end{claim}
\begin{cproof}
  Either $X_A$ or $X_B$ contains $f$. Without loss of generality, $f\in X_A$. Since $f$ colours some edge in $B$ there is an $f$-coloured path between an edge in $A$ incident with $x$ and an edge of $B$. Hence there is an $f$-coloured edge in $B$ incident with $x$ or $y$. In the latter case, there is also an $f$-coloured edge in $A$ incident with $y$ or an $f$-coloured edge in $B$ incident with $x$, and we are done. So assume that $f\in X_A$. No observe that there is an $f$-coloured path between the separating vertices of $(A,B)$ in $H$. By \autoref{lem:f-spanned_equiv}, our claim follows.
\end{cproof}

\begin{claim}\label{clm:get_rid}
  Let $f\in D$ so that $(A,B)$ is $f$-spanned. Then the pointy constraint matroid $(M,c,D-f)$ is realised by $G$ but not realised by $H$.
\end{claim}
\begin{cproof}
  Suppose for a contradiction that $H$ is a realisation of $(M,c,D-f)$. Then there is a vertex $v$ of $H$ that sees $D-f$. Since otherwise $G$ would be a realisation of $(M,c,D)$, the vertex $v$ is one of the vertices of $\{x,y\}$. But by hypothesis, there are edges at both $x$ and $y$ in $H$ coloured by $f$. Hence $H$ realises $(M,c,D)$, a contradiction.
\end{cproof}

By \autoref{clm:get_rid}, we assume that no $f\in D$ is so that $(A,B)$ is $f$-spanned.

\begin{claim}\label{clm:exactly_one}
  Every colour $f\in D$ is in exactly one of the sets $X_A,X_B,Y_A,Y_B$.
\end{claim}
\begin{cproof}
  By \autoref{clm:bothy} each colour in $D$ either colours edges in $A$ or edges in $B$ but not both. With \autoref{lem:f-spanned_equiv}, this proves our claim.
\end{cproof}

By the assumption that $H$ is not a realisation for $(M,c,D)$, there must exist $p\in X_A$ and $q\in X_B$ so that $p$ and $q$ are distinct colours. By \autoref{clm:exactly_one}, the subgraph of $H$ induced on edges coloured $p$ is disjoint from the subgraph of $H$ induced on edges coloured $q$. Hence the edges that receive colour $p*q$ do not form a connected subgraph of $H$ and hence it is not a realisation of $(M,c')$, a contradiction.
\end{pf}

\begin{lem}\label{lem:no_cont}
    Let $(M,c,D)$ be a pointy constraint matroid and $c'$ the extension of $c$ using $D$. Suppose that no component of $M$ contains all the colours of $D$. Then the set of realisations of the connected constraint matroid $(M,c')$ is identical to the set of realisations of $(M,c,D)$.
\end{lem}
\begin{pf}
    Clearly every realisation of $(M,c,D)$ is a reasliation of $(M,c')$ so we are required to show the converse. Let $G$ be a realisation of $(M,c')$.
    Let $L$ denote the subgraph of $G$ induced on the edges having colours from $D$. Since for every pair of distinct colours $p,q\in D$, we have that the edges coloured $p*q$ form a connected subgraph of $G$, we conclude that $L$ is connected. By hypothesis, there is a pair of colours $p,q\in D$ and a cutvertex $x$ of $G$ that separates the edges receiving $p$ and the edges receiving $q$. Let $r$ be an arbitrary colour of $D$. Suppose for a contradiction that no edge of colour $r$ is incident with $x$. But then the subgraph of $G$ coloured $r$ is not incident with the subgraph of $G$ induced on edges coloured $p$ or the subgraph of $G$ induced on edges coloured $q$, and hence either the subgraph induced on edges coloured $p*r$ or $q*r$ is not connected, a contradiction. Hence $x$ sees $D$ and this proves that $G$ is a realisation of $(M,c,D)$.
\end{pf}

\begin{pf}[Proof of \autoref{lem:using_D_works}]
    Let $(M,c,D)$ be a pointy constraint matroid and $c'$ the extension of $c$ using $D$. We are required to show that if $(M,c,D)$ admits a realisation, then the set of realisations of the connected constraint matroid $(M,c')$ is identical to the set of realisations of $(M,c,D)$.
    Let $G$ be an arbitrary realisation of $(M,c')$ and let $H$ be an arbitrary realisation of $(M,c,D)$. Given a component $B$ of $M$, we write $c|_B$ and $c'|_B$ for the restrictions to $B$ of the colourings $c$ and $c'$, respectively.
    
    Firstly, if there is no component of $M$ whose edges are coloured using all of the edges of $D$ then we are done by \autoref{lem:no_cont}. Secondly, if there are two components $B_1$ and $B_2$ of $M$, each of which uses all of $D$, then any separating vertex between $G[[B_1]]$ and $G[[B_2]]$ sees all of $D$ and in this case we are done. So we assume that there exists a component $B$ of $M$ containing all of the colours of $D$ and that it is the only one.
    \begin{claim}
        The subgraph $H[[B]]$ contains a vertex which sees $D$ in $H[[B]]$.
    \end{claim}
    \begin{cproof}
        Let $x$ be a vertex in $H$ that sees $D$. Since there is only one component which uses all the colours of $D$, we have that $x$ must lie in $H[[B]]$. Suppose for a contradiction that $x$ does not see all of $D$ in $H[[B]]$. Then there is some edge $e\in E(H)\setminus B$ incident with $x$ in $H$, and it receives some colour $f\in D$ which is not received by any edge incident with $x$ in $H[[B]]$. But then $e$ lies in a component separated from $H[[B]]$ by $x$ and since $B$ contains some edge receiving colour $f$, the subgraph of $H$ induced on edges of colour $f$ is disconnected. But this contradicts the fact that $H$ realises $(M,c,D)$. 
    \end{cproof}

    By the above claim, $H[[B]]$ is a realisation of $(B,c|_B,D)$. By \autoref{lem:red_to_con_from_pty_con}, every realisation of $(B,c'|_B)$ is a realisation of $(B,c|_B,D)$.

    \begin{claim}
        The subgraph $G[[B]]$ is a realisation of $(B,c'|_B)$.
    \end{claim}
    \begin{cproof}
        Let $p$ and $q$ be distinct colours from $D$. We are required to show that the subgraph induced on the edges of colour $p*q$, that is receiving either $p$ or $q$, is connected. Suppose for a contradiction that the subgraph $P$ of $G[[B]]$ induced of edges coloured $p$ is disjoint from the subgraph $Q$ of $G[[B]]$ on edges coloured $q$. But since $G$ is a realisation of $(M,c')$, we get that there is a path in $G\setminus E(B)$ between a vertex of $P$ and a vertex of $Q$. However, this contradicts the fact that $B$ forms a connected component of $G$.
    \end{cproof}

    By the above claim, $G$ is a realisation of $(M,c,D)$, and this completes our proof.
\end{pf}

\subsection{Quasibranches}\label{subsec:quasibranches}
Let $G$ be a graph. A \emph{block} of $G$ is a maximal connected subgraph of $G$ without a cutvertex. A \emph{branch} of $G$ is a subgraph $H$ so that $H$ meets the subgraph of $G$ induced on $E(G)\setminus E(H)$ in a single vertex. We will refer to a block which itself is branch as a \emph{leaf block} of $G$.

Given a constraint matroid $(M,c)$, let $\Bcal(M,c)$ denote its set of components. Let $\Hcal(M,c)$ denote the hypergraph whose vertices are the components $\Bcal(M,c)$ of $M$ and whose edges are the colours used by $c$, incident with a component whenever that component contains an edge with that colour.

We would like to be able to characterise using the data $\Hcal(M,c)$ which components of $(M,c)$ are leaves in its realisations. The following shows that such a characterisation is impossible in general.
\begin{eg}
    Let $\Hcal$ be an arbitrary hypergraph. Let $M$ be a matroid with edges $e_x$ for each vertex $x$ of $\Hcal$, and which has no circuits. Let $(M,c)$ be the constraint matroid obtained from $M$ by colouring $e_x$ by colours $f_h$ for each hyperedge $h$ in $\Hcal$ incident with $x$. Then $\Hcal=\Hcal(M,c)$. However, $(M,c)$ can be realised as a star. That is, there is a realisation where all of the components of $M$ appear as leaves.
\end{eg}

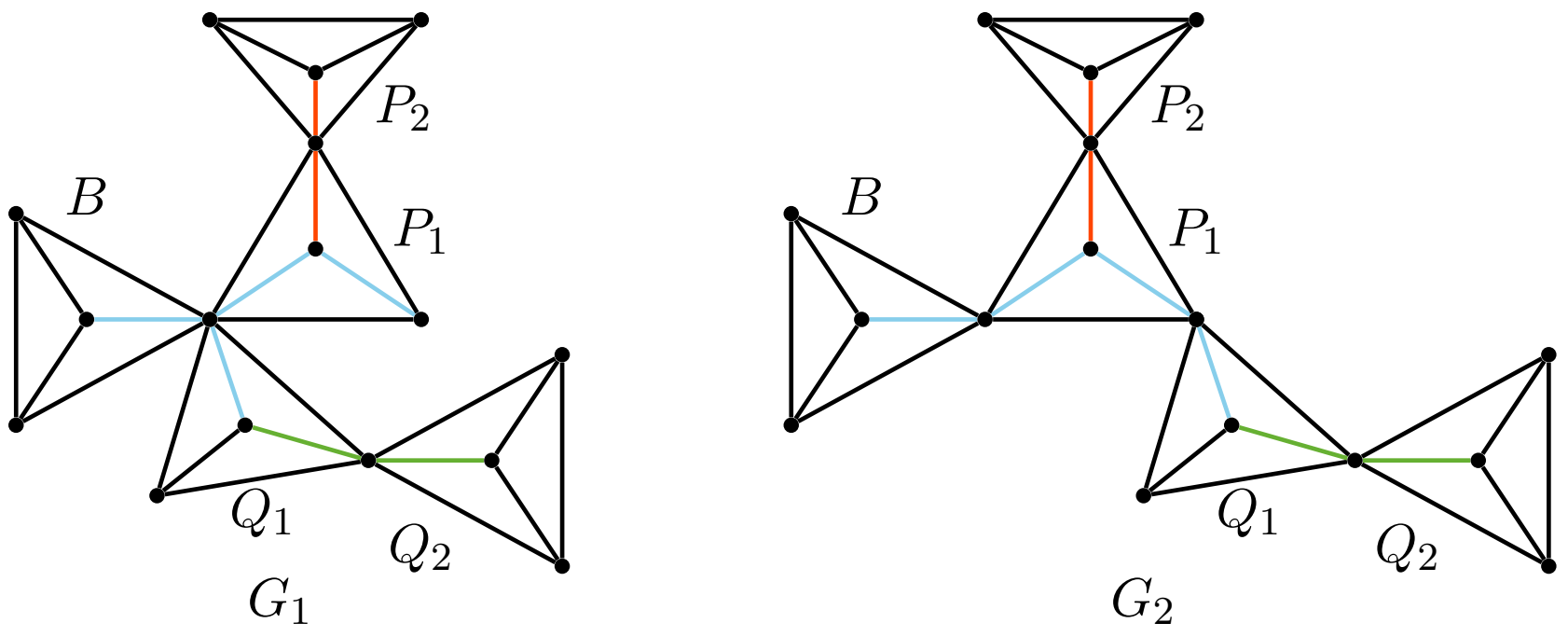
\begin{figure}
    \centering
    \begin{tikzpicture}
	\begin{pgfonlayer}{nodelayer}
		\node [style=black dot] (121) at (-8, 3) {};
		\node [style=black dot] (122) at (-8, 0) {};
		\node [style=black dot] (123) at (-5.25, 1.5) {};
		\node [style=black dot] (124) at (-7, 1.5) {};
		\node [style=black dot] (125) at (-3.75, 2.5) {};
		\node [style=black dot] (126) at (-6, -1) {};
		\node [style=black dot] (127) at (-2.25, 1.5) {};
		\node [style=black dot] (128) at (-4.75, 0) {};
		\node [style=black dot] (129) at (-3, -0.5) {};
		\node [style=black dot] (130) at (-3.75, 4) {};
		\node [style=black dot] (131) at (-5.25, 5.75) {};
		\node [style=black dot] (132) at (-2.25, 5.75) {};
		\node [style=black dot] (133) at (-3.75, 5) {};
		\node [style=black dot] (134) at (-0.25, -2) {};
		\node [style=black dot] (135) at (-0.25, 1) {};
		\node [style=black dot] (136) at (-1.25, -0.5) {};
		\node [style=none] (137) at (-7, 3.25) {$B$};
		\node [style=none] (138) at (-2.25, 2.75) {$P_1$};
		\node [style=none] (139) at (-2.5, 4.5) {$P_2$};
		\node [style=none] (140) at (-4.5, -1.25) {$Q_1$};
		\node [style=none] (141) at (-2.25, -1.75) {$Q_2$};
		\node [style=black dot] (142) at (3, 3) {};
		\node [style=black dot] (143) at (3, 0) {};
		\node [style=black dot] (144) at (5.75, 1.5) {};
		\node [style=black dot] (145) at (4, 1.5) {};
		\node [style=black dot] (146) at (7.25, 2.5) {};
		\node [style=black dot] (148) at (8.75, 1.5) {};
		\node [style=black dot] (151) at (7.25, 4) {};
		\node [style=black dot] (152) at (5.75, 5.75) {};
		\node [style=black dot] (153) at (8.75, 5.75) {};
		\node [style=black dot] (154) at (7.25, 5) {};
		\node [style=none] (158) at (4, 3.25) {$B$};
		\node [style=none] (159) at (8.75, 2.75) {$P_1$};
		\node [style=none] (160) at (8.5, 4.5) {$P_2$};
		\node [style=none] (163) at (-4.25, -2.5) {$G_1$};
		\node [style=none] (164) at (8, -2.5) {$G_2$};
		\node [style=black dot] (165) at (8.75, 1.5) {};
		\node [style=black dot] (166) at (8, -1) {};
		\node [style=black dot] (168) at (9.25, 0) {};
		\node [style=black dot] (169) at (11, -0.5) {};
		\node [style=black dot] (170) at (13.75, -2) {};
		\node [style=black dot] (171) at (13.75, 1) {};
		\node [style=black dot] (172) at (12.75, -0.5) {};
		\node [style=none] (173) at (9.5, -1.25) {$Q_1$};
		\node [style=none] (174) at (11.75, -1.75) {$Q_2$};
	\end{pgfonlayer}
	\begin{pgfonlayer}{edgelayer}
		\draw [style=c0] (121) to (122);
		\draw [style=c0] (122) to (124);
		\draw [style=c0] (124) to (121);
		\draw [style=c0] (121) to (123);
		\draw [style=c0] (123) to (122);
		\draw [style=c1] (124) to (123);
		\draw [style=c1] (123) to (125);
		\draw [style=c2] (130) to (125);
		\draw [style=c1] (125) to (127);
		\draw [style=c0] (127) to (123);
		\draw [style=c0] (123) to (130);
		\draw [style=c0] (130) to (127);
		\draw [style=c0] (123) to (129);
		\draw [style=c0] (129) to (126);
		\draw [style=c0] (126) to (123);
		\draw [style=c1] (123) to (128);
		\draw [style=c0] (128) to (126);
		\draw [style=c3] (128) to (129);
		\draw [style=c0] (130) to (131);
		\draw [style=c0] (131) to (132);
		\draw [style=c0] (132) to (133);
		\draw [style=c0] (133) to (131);
		\draw [style=c2] (130) to (133);
		\draw [style=c0] (132) to (130);
		\draw [style=c0] (129) to (135);
		\draw [style=c0] (135) to (134);
		\draw [style=c0] (134) to (129);
		\draw [style=c3] (129) to (136);
		\draw [style=c0] (136) to (134);
		\draw [style=c0] (136) to (135);
		\draw [style=c0] (142) to (143);
		\draw [style=c0] (143) to (145);
		\draw [style=c0] (145) to (142);
		\draw [style=c0] (142) to (144);
		\draw [style=c0] (144) to (143);
		\draw [style=c1] (145) to (144);
		\draw [style=c1] (144) to (146);
		\draw [style=c2] (151) to (146);
		\draw [style=c1] (146) to (148);
		\draw [style=c0] (148) to (144);
		\draw [style=c0] (144) to (151);
		\draw [style=c0] (151) to (148);
		\draw [style=c0] (151) to (152);
		\draw [style=c0] (152) to (153);
		\draw [style=c0] (153) to (154);
		\draw [style=c0] (154) to (152);
		\draw [style=c2] (151) to (154);
		\draw [style=c0] (153) to (151);
		\draw [style=c0] (165) to (169);
		\draw [style=c0] (169) to (166);
		\draw [style=c0] (166) to (165);
		\draw [style=c1] (165) to (168);
		\draw [style=c0] (168) to (166);
		\draw [style=c3] (168) to (169);
		\draw [style=c0] (169) to (171);
		\draw [style=c0] (171) to (170);
		\draw [style=c0] (170) to (169);
		\draw [style=c3] (169) to (172);
		\draw [style=c0] (172) to (170);
		\draw [style=c0] (172) to (171);
	\end{pgfonlayer}
\end{tikzpicture}
    \caption{Two realisations $G_1$ and $G_2$ of a constraint matroid with five components, namely $B,P_1,P_2,Q_1,Q_2$. The quasibranches with respect to $B$ are $\{P_1,P_2\}$ and $\{Q_1,Q_2\}$.}
    \label{fig:branches_everywhere}
\end{figure}

\begin{dfn}[Quasibranch]
    Let $\Hcal(M,c)$ be the hypegraph associated with a constraint matroid $(M,c)$. Let $B$ be a component of $M$. Then the \emph{quasibranches of $(M,c)$ with respect to $B$} are the vertex sets of connected components remaining after removing $B$ and all its incident hyperedges from $\Hcal(M,c)$. The \emph{quasibranches} of $(M,c)$ is the union of the sets of quasibranches with respect to each component of $(M,c)$.
\end{dfn}
See \autoref{fig:branches_everywhere} for a depiction of quasibranches.

\begin{lem}\label{lem:quasibranch_connected}
    Let $(M,c)$ be a constraint matroid and $\Bcal$ a quasibranch thereof. Then $\Bcal$ induces a connected subgraph in every realisation of $(M,c)$.
\end{lem}
\begin{pf}
    Let $B$ denote the component of $M$ with respect to which $\Bcal$ is a quasibranch. Let $G$ be an arbitrary realisation of $(M,c)$. Let $B'$ and $B''$ be components in $\Bcal$. Then there is a sequence of components $B'=B_1,B_2,\ldots,B_n=B''$ so that each consecutive pair $B_i,B_{i+1}$ shares some colour $f_i$ outside of $D$. Hence each consecutive pair is joined by a path in the colour $f_i$ in $G$. Any component containing the $f_i$-coloured path is in $\Bcal$. This proves that $\Bcal$ is connected in $G$.
\end{pf}

The following lemma describes in what sense `quasibranches' can be considered `branches' and the subsequent example shows the limits of this correspondence.
\begin{lem}\label{lem:quasibranch_main}
    Let $(M,c)$ be a constraint matroid and $\Bcal_1\ldots,\Bcal_k$ a set of quasibranches thereof, with respect to the same block. If $(M,c)$ admits a realisation, then there is a realisation of $(M,c)$ in which each $\Bcal_i$ induces a branch.
\end{lem}
\begin{pf}
    Let $B$ be a component of $M$ with respect to which each $\Bcal_i$ is a quasibranch. Denote the collection of quasibranches with respect to $B$ by $\Bcal_1,\ldots,\Bcal_n$. Let $G$ be an arbitrary realisation of $(M,c)$. By \autoref{lem:quasibranch_connected}, every quasibranch $\Bcal_i$ induces a connected subgraph $H_i$ of $G$. For each $i\in[n]$, let $x_i$ be the unique vertex of $G$ contained in $H_i$ that separates $\Bcal_i$ from $B$. That is, $x_i$ is the closest separator to $H_i$ of all separators between $\Bcal_i$ and $B$. Let $H_0$ be the subgraph of $G$ induced on $B$. We require the following three claims.
    \begin{claim}\label{clm:quasi1}
        For each $i\in[n]$, the vertex $x_i$ sees in $H_i$ the colours appearing in both $B$ and $\Bcal_i$.
    \end{claim}
    \begin{cproof}
        The vertex $x_i$ is the unique vertex on all paths between $H_i$ and $H_0$ in $G$. Since $G$ realises $(M,c)$, the vertex $x_i$ sees in $H_i$ all colours appearing in both $B$ and $\Bcal_i$, as required.
    \end{cproof}

    \begin{claim}\label{clm:quasi2}
        For each $i\in[n]$, the graph $H_i$ is a representation for the restriction $(\Bcal_i,c|_{\Bcal_i})$.
    \end{claim}
    \begin{cproof}
        Let $f$ be a colour used by $c|_{\Bcal_i}$.  If $f$ only appears in $\Bcal_i$, then the edges coloured $f$ form a connected subset of $H_i$. Otherwise, $f$ appears in both $\Bcal_i$ and $B$. Let $e$ be an edge in $H_i$ recieving colour $f$. Then there is an $f$-coloured path in $G$ from $e$ to $H_0$. This path must use $x_i$ and hence all such edges $e$ form a connected subgraph of $H_i$.
    \end{cproof}

    \begin{claim}\label{clm:quasi3}
        For each $i\in[n]$, there is a vertex $y_i$ in $H_0$ that sees in $H_0$ the colours appearing in both $B$ and $\Bcal_i$.
    \end{claim}
    \begin{cproof}
        Let $y_i$ be the unique vertex in $H_0$ on a path between $H_0$ and $H_i$ that is also in the graph induced on $E(G)\sm E(H_0)$ in $G$. Similarly to in the proof of \autoref{clm:quasi1}, the vertex $y_i$ sees in $H_0$ all of the colours appearing in both $B$ and $\Bcal_i$.
    \end{cproof}

    Let $G'$ be the graph obtained from the disjoint union of the graphs $H$ and $H_1,\ldots,H_n$ by identifying the pairs of vertices $(x_i,y_i)$ for each $i\in[n]$. By Claims~\ref{clm:quasi1}, \ref{clm:quasi2} and \ref{clm:quasi3}, we have that $G'$ is a realisation of $(M,c)$. Furthermore, the quasibranch $\Bcal$ induces a branch in $G'$.
\end{pf}

\begin{eg}
    Not in all realisations are quasibranches branches. See \autoref{fig:branches_everywhere}.
\end{eg}

\subsection{The proof of \autoref{thm:new_red_to_con}}\label{subsec:proof2}
Let $(M,c)$ be a constraint matroid. Let $B$ be a component of $M$ and $\Bcal_1,\ldots,\Bcal_k$ the quasibranches of $(M,c)$ with respect to $B$. For $i\in[k]$, let $D_i$ denote the set of colours used by $c$ in both $B$ and~$\Bcal_i$. The \emph{quasibranch colouring} of $B$ is the extension of the restriction $c|_B$ using the sets $D_1,\ldots,D_k$ and we denote it by $c_B^*$.

The \emph{fragments of $(M,c)$} with respect to $B$ are the constraint matroids $(B_i^+,c_i^+)$ where $B_i^+$ is obtained as the disjoint union of the quasibranch $\Bcal_i$ with a new loop $\ell$, and the colouring $c_i^+$ restricts to $c$ on $\Bcal_i$ and the loop $\ell$ receives exactly the colours $D_i$. That is, we add a loop to the quasibranch having the colours that cross into $B$.

\begin{lem}\label{lem:fragments}
    Let $(M,c)$ be a constraint matroid and $B$ a component of $M$. Let $(B_i^+,c_i^+)$ denote the fragments of $(M,c)$ with respect to $B$. Then\vspace{-.25cm}
    \begin{enumerate}
        \item every realisation of $(M,c)$ induces a realisation for $(B,c_B^*)$ and each $(B_i^+,c_i^+)$; and\vspace{-.25cm}
        \item if $(M,c)$ admits a realisation, then given a realisation $G_0$ of $(B,c_B^*)$ and for each $i\in[k]$ a realisation $G_i$ of $(B_i^+,c_i^+)$, there is a realisation of $(M,c)$ which induces the graphs $G_0,G_1,\ldots,G_k$.\vspace{-.25cm}
    \end{enumerate}
\end{lem}
\begin{pf}
    Let $D_i$ be the set of colours used by $c$ appearing in both $B$ and $B_i$.
    For (i), let a realisation $G$ of $(M,c)$ be given. By \autoref{lem:quasibranch_connected}, the quasibranches of $(M,c)$ with respect to $B$ are connected in $G$. Hence for each $i\in[k]$, there is a vertex $x_i$ in $G[[B]]$ so that $x_i$ sees $D_i$ in $G[[B]]$. We conclude that $G[[B]]$ is a realisation for the pointy constraint matroid $(B,c|_B,\{D_1,\ldots,D_k\})$. By \autoref{cor:D_works_MORE}, $G[[B]]$ is a realisation for $(B,c_B^*)$, as required.

    For (ii), let the graphs $G_0,G_1,\ldots,G_k$ be given. By \autoref{cor:D_works_MORE}, $G_0$ is a realisation for the pointy coloured matroid $(B,c|_B,\{D_1,\ldots,D_k\})$. For each $i\in[k]$, let $x_i$ denote a vertex in $G_0$ that sees the colours $D_i$. For each $i\in[k]$, let $y_i$ be the vertex in $G_i$ incident with the additional loop~$\ell_i$. Then $y_i$ sees the colours $D_i$ in $G_i-\ell_i$. Let $G$ be the graph obtained from the disjoint union of $G_0$ and $G_i-\ell_i$ for each $i\in[k]$ by identifying $x_i$ with $y_i$ for each $i\in[k]$. We claim that $G$ is a realisation of $(M,c)$. Indeed, let $p$ be a colour used by $c$ that is used in more than one of the parts $B,B_1,\ldots,B_k$. Then by the definition of quasibranches, we get that $p$ appears in $B$. Hence the edges coloured $p$ form a connected subgraph of $G$, as required. This completes our proof.
\end{pf}

\begin{eg}\label{eg:silly}
    Let $(M,c)$ be a constraint matroid, all whose elements are loops. Then $(M,c)$ admits a realisation. Namely, we glue all loops to a single shared vertex.
\end{eg}

\begin{pf}[Proof of \autoref{thm:new_red_to_con}]
    Let $(M,c)$ be a constraint matroid. We are required to find for each component $B$ of $M$, an extension $c_B^*$ of the restriction of $c$ to $B$ so that\vspace{-.25cm}
    \begin{enumerate}
        \item every realisation of $(M,c)$ induces for each component $B$ of $M$ a realisation for $(B,c_B^*)$; and\vspace{-.25cm}
        \item if $(M,c)$ admits a realisation, then given for each component $B$ of $M$ a realisation $G_B$ of $(B,c_B^*)$, there is a realisation of $(M,c)$ which induces the graphs $G_B$.\vspace{-.25cm}
    \end{enumerate}
    For a component $B$ of $M$, we take $c_B^*$ to be the quasibranch colouring for $B$.
    
    For (i), let $G$ be a realisation of $(M,c)$ and $p$ and $q$ colours used in $c$ so that $p*q$ is a colour used by $c_B^*$. Then there exists a block $B$ of $M$ and a quasibranch $\Bcal$ with respect to $B$ so that both $p$ and $q$ appear in both $B$ and $\Bcal$. By \autoref{lem:quasibranch_connected}, $\Bcal$ forms a connected subgraph of $(M,c)$. Hence we conclude that there is a vertex $x$ in $G[[B]]$ that is in every paths between $G[[B]]$ and $G[[\Bcal]]$ in $G$. Hence $x$ sees $p$ and $q$ in $G[[G]]$. Since $G$ is a realisation of $(M,c)$ the edges coloured $p$ (resp.\ $q$) form a connected subgraph of $G$ and since these connected subgraphs share $x$, the edges coloured $p*q$ in $G$ also form a connected subgraph of $G$.

    For (ii), we proceed by induction on the number of components of $M$ with more than one edge. Observe that if there are no components of $M$ with more than one edge, then all of the components of $M$ are single edges and the single edges have unique realisations, which in particular respect any colouring. Furthermore, $(M,c)$ is realisable by gluing all of these single edges at some shared vertex. So we are done in this case.

    So suppose that $M$ has at least one component with more than one edge and suppose for induction that our hypothesis holds for all matroids with fewer than this many components with more than one edge. Let $B$ be a component of $M$ having at least two edges. Let $(B_i^+,c_i^+)$ denote the fragments of $(M,c)$ with respect to $B$. By \autoref{lem:fragments}, if $(M,c)$ admits a realisation, then given a realisation $G_0$ of $(B,c_B^*)$ and for each $i\in[k]$ a realisation $G_i$ of $(B_i^+,c_i^+)$, there is a realisation of $(M,c)$ which induces the graphs $G_0,G_1,\ldots,G_k$.
    \begin{claim}
        For some $i\in[k]$, let $B'$ be a component of $B_i^+$ that is not a loop. Then the $B'$ is a component of $M$ and the quasibranch colouring of $B'$ in the context of $(B_i^+,c_i^+)$ is extended by the quasibranch colouring of $B'$ in the context of $(M,c)$.
    \end{claim}
    \begin{cproof}
        Clearly $B'$ is a component of $M$. Our claim thus follows from the fact that the quasibranches of $B'$ in $(B_i^+,c_i^+)$ are either quasibranches of $B'$ in $(M,c)$ or contain the loop $\ell$ added to $B_i$ to obtain $B_i^+$.
    \end{cproof}
    
    By induction, the realisations of each fragment $(B_i^+,c_i^+)$ can be obtained from a realisation of each of their components under the quasibranch colouring in the context of $(B_i^+,c_i^+)$. Since by the above claim, this colouring extends that of the quasibranch colouring of the blocks with respect to $(M,c)$, we are done.
\end{pf}

\begin{figure}
    \centering
    \begin{tikzpicture}
	\begin{pgfonlayer}{nodelayer}
		\node [style=black dot] (121) at (-8, 3) {};
		\node [style=black dot] (122) at (-8, 0) {};
		\node [style=black dot] (123) at (-5.25, 1.5) {};
		\node [style=black dot] (124) at (-7, 1.5) {};
		\node [style=black dot] (125) at (-4.5, 2.75) {};
		\node [style=black dot] (127) at (-2.5, 2.5) {};
		\node [style=black dot] (130) at (-5.25, 4.25) {};
		\node [style=none] (137) at (-6.75, 3.25) {$B$};
		\node [style=none] (138) at (-3, 3.75) {$P$};
		\node [style=none] (164) at (-5, -1.75) {$A$};
		\node [style=black dot] (170) at (3.5, 1.5) {};
		\node [style=black dot] (172) at (4.25, 2.75) {};
		\node [style=black dot] (173) at (6.25, 2.5) {};
		\node [style=black dot] (174) at (3.5, 4.25) {};
		\node [style=none] (176) at (5.5, 3.75) {$P$};
		\node [style=none] (177) at (5.25, -0.5) {$Q$};
		\node [style=black dot] (178) at (3.75, 0) {};
		\node [style=black dot] (179) at (7.25, 0) {};
		\node [style=black dot] (180) at (5.75, 1) {};
		\node [style=black dot] (181) at (-0.5, 3) {};
		\node [style=black dot] (182) at (-0.5, 0) {};
		\node [style=black dot] (183) at (2.25, 1.5) {};
		\node [style=black dot] (184) at (0.5, 1.5) {};
		\node [style=none] (185) at (0.75, 3) {$B$};
		\node [style=black dot] (187) at (-5.25, 1.5) {};
		\node [style=black dot] (188) at (-7.75, -1) {};
		\node [style=black dot] (189) at (-4.25, -1) {};
		\node [style=black dot] (190) at (-5.75, 0) {};
		\node [style=none] (191) at (-4.25, 0.5) {$Q$};
		\node [style=none] (192) at (0.75, -1.75) {$B$};
		\node [style=none] (193) at (5.25, -1.75) {$C$};
	\end{pgfonlayer}
	\begin{pgfonlayer}{edgelayer}
		\draw [style=c0] (121) to (122);
		\draw [style=c0] (122) to (124);
		\draw [style=c0] (124) to (121);
		\draw [style=c0] (121) to (123);
		\draw [style=c2] (123) to (122);
		\draw [style=c1] (124) to (123);
		\draw [style=c0] (123) to (125);
		\draw [style=c0] (130) to (125);
		\draw [style=c0] (125) to (127);
		\draw [style=c3] (127) to (123);
		\draw [style=c1] (123) to (130);
		\draw [style=c0] (130) to (127);
		\draw [style=c0] (170) to (172);
		\draw [style=c0] (174) to (172);
		\draw [style=c0] (172) to (173);
		\draw [style=c3] (173) to (170);
		\draw [style=c1] (170) to (174);
		\draw [style=c0] (174) to (173);
		\draw [style=c0] (173) to (179);
		\draw [style=c0] (179) to (178);
		\draw [style=c2] (178) to (173);
		\draw [style=c3] (173) to (180);
		\draw [style=c0] (180) to (178);
		\draw [style=c0] (180) to (179);
		\draw [style=c0] (181) to (182);
		\draw [style=c0] (182) to (184);
		\draw [style=c0] (184) to (181);
		\draw [style=c0] (181) to (183);
		\draw [style=c2] (183) to (182);
		\draw [style=c1] (184) to (183);
		\draw [style=c0] (187) to (189);
		\draw [style=c0] (189) to (188);
		\draw [style=c2] (188) to (187);
		\draw [style=c3] (187) to (190);
		\draw [style=c0] (190) to (188);
		\draw [style=c0] (190) to (189);
	\end{pgfonlayer}
\end{tikzpicture}
    \caption{Three graphs $A$, $B$ and $C$. Let $(M,c)$ be the coloured matroid given by $A$. Then $(M,c)$ is realised uniquely by $A$ and has blocks $B,P$ and $Q$. The set $\{P,Q\}$ is a quasibranch with respect to $B$. The restrictions of $(M,c)$ to $B$ and $\{P,Q\}$ are realised by $B$ and $C$, respectively. However, $C$ is not the restriction of a realisation of $(M,c)$.}
    \label{fig:quasi_necessity}
\end{figure}
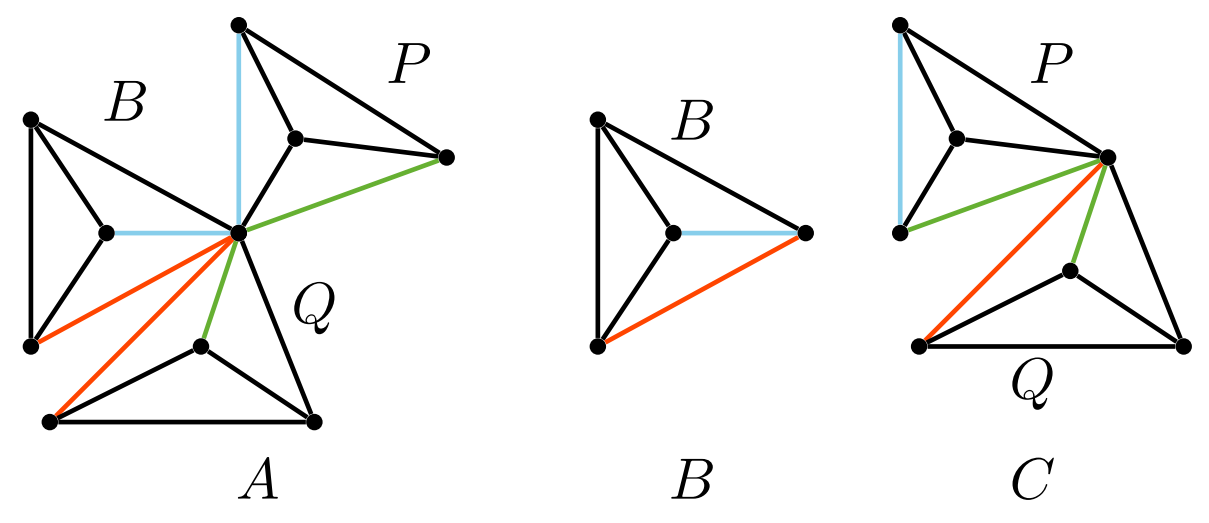

\begin{eg}
    One might imagine there is a simpler colouring than the quasibranch colouring. Indeed, what happens if one just adds a colour $p*q$ for every pair $p$ and $q$ that both appear in two distinct blocks? This is insufficient, as \autoref{fig:quasi_necessity} shows.
\end{eg}

\section{Algorithmics}
\label{sec:algorithmics}
In this section, we prove \autoref{cor:alg}, which is recalled below for the convenience of the reader.
\begingroup
\def\thethm{\ref{cor:alg}}
\begin{thm}
    The realisability of a constraint matroid $(M,c)$ can be determined in polynomial time with respect to the number of edges of $M$ and the number of colours used by $c$.
\end{thm}
\addtocounter{lem}{-1}
\endgroup

For the above, we require the following specialisation to connected constraint matroids.
\begin{lem}\label{lem:alg_con}
    The realisability of a connected constraint matroid $(M,c)$ can be determined in polynomial time with respect to the number of edges of $M$ and the number of colours used by $c$.
\end{lem}
\begin{pf}
    Let $(M,c)$ be a connected constraint matroid and let $n$ and $m$ denote the number of elements of $M$ and colours used by $c$, respectively.

    We begin by computing the Tutte decomposition of $M$, which runs in polynomial~\cite{cunningham1980combinatorial}. The Tutte decomposition has at most $O(n)$ torsos. For each torso of the Tutte decomposition, we check to see if it is unbalanced or unstable, which are both clearly tractable, and if any of them are unbalanced or unstable, we output that $(M,c)$ admits no realisation (as per \autoref{cor:not_unbalanced} and \autoref{lem:not_unstable}). Otherwise, we compute the quasicolourings of the directed totally nested $2$-separations of $(M,c)$, of which there are at most $O(n)$. If some totally nested $2$-separation is disagreeable, we output that $(M,c)$ admits no realisation (as per \autoref{cor:not_disagreeable}). Otherwise, we compute the mixed-sum decomposition into $3$-connected constraint matroids, constraint circuits and constraint diatoms given by \autoref{thm:matroids_tutte_analogue}. The original constraint matroid $(M,c)$ admits a realisation if and only if all of its pieces admit realisations. Checking whether these $O(n)$ pieces are realisable is a tractable process which is outlined in \autoref{sec:circuits}, \autoref{sec:diatoms} and, in the $3$-connected case, by inspecting the unique realisation. We note that each piece is coloured using at most $m+O(n)$ colours, since we add at most one colour for every totally-nestd $2$-separator $(A,B)$ left-incident with each torso, namely the new colour named `$f_{(A,B)}$'.
\end{pf}

\begin{pf}[Sketch proof of \autoref{cor:alg}]
    Let $(M,c)$ be a constraint matroid and let $n$ and $m$ denote the number of elements of $M$ and colours used by $c$, respectively.

    Our first aim is to find a realisation $G_B$ for each component $B$ of $M$ a realisation of the connected constraint matroid $(B,c_B^*)$ (of which there are at most $n$). This can be done in polynomial time by \autoref{lem:alg_con}, since each block has at most $n$ edges and $c_B^*$ uses at most $m+m^2$ colours.

    Now it remains to assemble these graphs $G_B$ into a realisation for $(M,c)$. This is possible by \autoref{thm:new_red_to_con}, but it is not clear from the hypothesis of this theorem how one should efficiently compute such an assembly.
    Following the proof of \autoref{thm:new_red_to_con}, we use \autoref{lem:fragments} recursively. Let $\#_{>1}(M)$ denote the number of components of $M$ with more than one edge and $\#_{=1}(M)$ denote the number of components of $M$ with exactly one edge. Then \autoref{lem:fragments} reduces our problem to a block $(B,c_B^*)$ (to which we assign the realisation $G_B$) and $k$ fragments $(B_i^+,c_i^+)$ so that $k\le n$,
    \[
        \sum_{i\in[k]} \#_{>1}(B_i) = \#_{>1}(M)-1
    \]
    and
    \[
        \sum_{i\in[k]} |E(B_i)| \le n + k,
    \]
    which means our approach works in polynomial time, as required.
\end{pf}

\bibliographystyle{plain}
\bibliography{wliterature}

@book{oxley2006,
  title={Matroid Theory},
  author={Oxley, J.G.},
  isbn={9780199202508},
  series={Oxford graduate texts in mathematics},
  year={2006},
  publisher={Oxford University Press}
}

@article{diestel2006,
  title={Graph Theory, Graduate Texts in Mathematics. Springer},
  author={Diestel, R},
  year={2006}
}

@article{carmesin2017embeddingIII,
  title={Embedding simply connected 2-complexes in 3-space--III. Constraint minors},
  author={Carmesin, Johannes},
  journal={arXiv preprint arXiv:1709.04645},
  year={2017}
}

@article{cunningham1980combinatorial,
  title={A combinatorial decomposition theory},
  author={Cunningham, William H and Edmonds, Jack},
  journal={Canadian Journal of Mathematics},
  volume={32},
  number={3},
  pages={734--765},
  year={1980},
  publisher={Cambridge University Press}
}

@article{carmesin2017embeddingIV,
  title={Embedding simply connected 2-complexes in 3-space--IV. Dual matroids},
  author={Carmesin, Johannes},
  journal={arXiv preprint arXiv:1709.04652},
  year={2017}
}

@article{geelen2018quasi,
  title={Quasi-graphic matroids},
  author={Geelen, Jim and Gerards, Bert and Whittle, Geoff},
  journal={Journal of Graph Theory},
  volume={87},
  number={2},
  pages={253--264},
  year={2018},
  publisher={Wiley Online Library}
}

@article{maimani2023structure,
  title={On the structure of matroids arising from the gain graphs},
  author={Parsaei-Majd, L},
  journal={Discrete Mathematics},
  volume={346},
  number={12},
  pages={113637},
  year={2023},
  publisher={Elsevier}
}

@article{me2025,
  title={A Graph Minors Approach to Temporal Sequences},
  author={Carmesin, Johannes and Turner, Will J.},
  publisher = {arXiv},
  year = {2025},
  copyright = {arXiv.org perpetual, non-exclusive license}
}

@article{matroid_overlay,
  title={Planarity for Contraction Trees of Graphs},
  author={Turner, Will J.},
  year = {in preparation}
}

@article{whitney19332,
  title={2-isomorphic graphs},
  author={Whitney, Hassler},
  journal={American Journal of Mathematics},
  volume={55},
  number={1},
  pages={245--254},
  year={1933},
  publisher={JSTOR}
}

@article{tutte1959matroids,
  title={Matroids and graphs},
  author={Tutte, William Thomas},
  journal={Transactions of the American Mathematical Society},
  volume={90},
  number={3},
  pages={527--552},
  year={1959},
  publisher={JSTOR}
}

@article{tutte1960algorithm,
  title={An algorithm for determining whether a given binary matroid is graphic},
  author={Tutte, WT},
  journal={Proceedings of the American Mathematical Society},
  volume={11},
  number={6},
  pages={905--917},
  year={1960},
  publisher={JSTOR}
}

\appendix
\begin{center}
{\bf \LARGE Appendix}
\end{center}

\section{Bunches}\label{appendix:bunches}
\begin{figure}
    \centering
    \begin{tikzpicture}
	\begin{pgfonlayer}{nodelayer}
		\node [style=black dot] (152) at (-1.25, 0) {};
		\node [style=black dot] (153) at (1.5, 1.5) {};
		\node [style=black dot] (154) at (1.5, -1.5) {};
		\node [style=none] (155) at (-2.75, 1) {};
		\node [style=none] (156) at (-2.75, -1) {};
		\node [style=none] (159) at (1.75, -2.25) {};
		\node [style=none] (160) at (1.75, 2.25) {};
		\node [style=none] (161) at (1.5, 2) {};
		\node [style=none] (162) at (1.5, -2) {};
		\node [style=none] (163) at (2, 1.75) {};
		\node [style=none] (164) at (-1.75, -0.25) {};
		\node [style=none] (165) at (2, -1.75) {};
		\node [style=none] (166) at (-1.75, 0.25) {};
		\node [style=black dot] (167) at (-4.25, 0) {};
		\node [style=black dot] (168) at (-7, 1.5) {};
		\node [style=black dot] (169) at (-7, -1.5) {};
		\node [style=none] (170) at (-2.75, 1) {};
		\node [style=none] (171) at (-2.75, -1) {};
		\node [style=none] (172) at (-7.25, -2.25) {};
		\node [style=none] (173) at (-7.25, 2.25) {};
		\node [style=none] (174) at (-7, 2) {};
		\node [style=none] (175) at (-7, -2) {};
		\node [style=none] (176) at (-7.5, 1.75) {};
		\node [style=none] (177) at (-3.75, -0.25) {};
		\node [style=none] (178) at (-7.5, -1.75) {};
		\node [style=none] (179) at (-3.75, 0.25) {};
		\node [style=none] (180) at (-2.75, -2.5) {$\Hcal$};
	\end{pgfonlayer}
	\begin{pgfonlayer}{edgelayer}
		\draw [style=d1, in=0, out=-180] (160.center) to (155.center);
		\draw [style=d1, in=-180, out=0] (156.center) to (159.center);
		\draw [style=d1, bend right=90, looseness=0.50] (159.center) to (160.center);
		\draw [style=d1, bend left=90, looseness=0.25] (161.center) to (162.center);
		\draw [style=d1, bend left=90, looseness=0.25] (162.center) to (161.center);
		\draw [style=d1, bend right=105, looseness=0.25] (163.center) to (164.center);
		\draw [style=d1, bend right=75, looseness=0.25] (164.center) to (163.center);
		\draw [style=d1, bend left=90, looseness=0.25] (166.center) to (165.center);
		\draw [style=d1, bend left=105, looseness=0.25] (165.center) to (166.center);
		\draw [style=d1, in=180, out=0] (173.center) to (170.center);
		\draw [style=d1, in=0, out=180] (171.center) to (172.center);
		\draw [style=d1, bend left=90, looseness=0.50] (172.center) to (173.center);
		\draw [style=d1, bend right=90, looseness=0.25] (174.center) to (175.center);
		\draw [style=d1, bend right=90, looseness=0.25] (175.center) to (174.center);
		\draw [style=d1, bend left=105, looseness=0.25] (176.center) to (177.center);
		\draw [style=d1, bend left=75, looseness=0.25] (177.center) to (176.center);
		\draw [style=d1, bend right=90, looseness=0.25] (179.center) to (178.center);
		\draw [style=d1, bend right=105, looseness=0.25] (178.center) to (179.center);
	\end{pgfonlayer}
\end{tikzpicture}
    \caption{A connected hypergraph $\Hcal$, consisting of 6 vertices and 7 hyperedges, which has neither a singleton quasibranch, nor component with respect to which there is at most one quasibranch.}
    \label{fig:noquasibranch}
\end{figure}
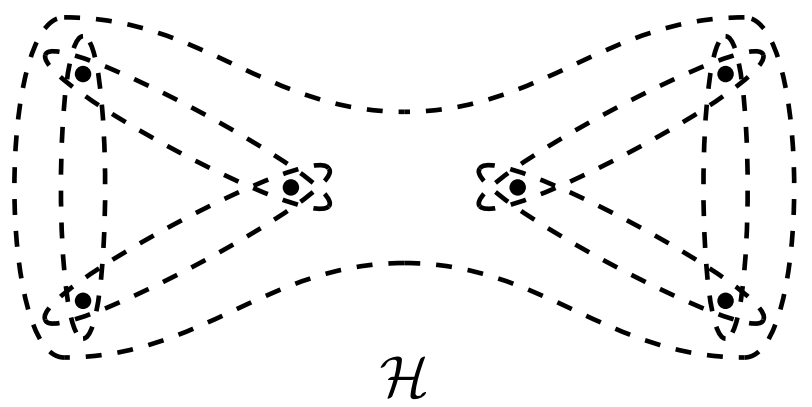

If a given realisable constraint matroid $(M,c)$ contains a singleton quasibranch, or a component which has at most one quasibranch with respect to it, then by \autoref{lem:quasibranch_main} there is some realisation in which it is a leaf. However, there are some constraint matroids where neither object exists (see~\autoref{fig:noquasibranch}).
We would like to determine which blocks of $M$ can appear as leaves in some realisation of $(M,c)$.

Let $G$ be a graph. A \emph{bunch} of $G$ is a subset of the leaf blocks of $G$ that are all incident with a shared vertex in $G$. A bunch is \emph{terminal} if it is incident with at most one other block in $G$.
\begin{lem}
    Every graph contains a terminal bunch.
\end{lem}
\begin{pf}
    Let $G$ be a graph and $H$ the subgraph of $G$ induced on the edges of $G$ that are not in leaf blocks. If $H$ is empty, then $G$ itself is a terminal bunch. So assume $H$ is nonempty. Let $B_0$ be a leaf block of $H$. Now observe that there is a bunch of $G$ which meets only $B_0$, as desired.
\end{pf}

\begin{lem}\label{lem:bunch_branch}
    Given a realisation $G$ of a constraint matroid $(M,c)$, every terminal bunch of $G$ contains a quasibranch of $(M,c)$.
\end{lem}
\begin{pf}
    Let $B$ be a terminal bunch of $G$. If $B$ is the whole of $G$, then we are done. So assume that this is not the case and let $B_0$ be the unique block of $G$ incident with $B$ but not contained in it. Since $B$ is terminal, it is an inclusion-maximal bunch outside of $B_0$. Hence we observe that at least one of the quasibranches of $(M,c)$ with respect to $B_0$ is contained in $B$. Indeed, let $\Bcal$ be a quasibranch of $M$ with respect to $B_0$. Then in $G$ every block in $\Bcal$ meeting $B_0$ meets a shared vertex of $B_0$. If $\Bcal$ has nonempty intersection with $B$, then $B$ must coincide with $\Bcal$.
\end{pf}
A \emph{quasibunch} of a constraint matroid $(M,c)$ is a minimal nonempty quasibranch under inclusion (viewed as sets of components of $M$).

\begin{lem}\label{lem:branch_brunch}
    Let $(M,c)$ be a constraint matroid and $\Bcal$ a quasibunch thereof. Suppose that $G$ is a realisation of $(M,c)$ where $\Bcal$ is a branch. Then $\Bcal$ induces a bunch in $G$.
\end{lem}
\begin{pf}
    Let $B$ denote the branch induced by $\Bcal$ in $G$. Suppose for a contradiction that $\Bcal$ does not induce a bunch in $G$. Then there is an terminal bunch $B'$ of $G$ which is strictly contained in $B$. By \autoref{lem:bunch_branch}, the bunch $B'$ contains a quasibranch of $(M,c)$, contradicting the minimality of~$\Bcal$.
\end{pf}

In conclusion, the members of a quasibunch, of which there is at least one, can be thought of as a sort of `quasileaves'.

\section{The case of $3$-connected constraint matroids}\label{appendix:translation}
In \cite{carmesin2017embeddingIII}, Carmesin proves the following.
\begin{thm}[Carmesin 2017 \cite{carmesin2017embeddingIII}]\label{og}
    Let $(M,c)$ be a $3$-connected constraint matroid that uses at most one colour. Then it is realisable if and only if we cannot obtain a constraint matroid isomorphic to one of the ones depicted in \autoref{fig:obst3con} by successive use of the following operations:\vspace{-.25cm}
    \begin{itemize}
        \item contracting an edge; or\vspace{-.25cm}
        \item deleting an edge that does not recieve a colour.\vspace{-.25cm}
    \end{itemize}
\end{thm}

In this appendix, we show how the characterisation does not change if we replace the operations in the above theorem by our own operations for constraint-minors.

\begin{lem}
    Let $(M,c)$ be a $3$-connected constraint matroid. Then it is realisable if and only if we cannot obtain a constraint matroid isomorphic to one of the ones depicted in \autoref{fig:obst3con} by successive use of the operations:\vspace{-.25cm}
    \begin{itemize}
        \item forgetting a colour; and\vspace{-.25cm}
        \item contracting an edge; and\vspace{-.25cm}
        \item merging parallel edges.\vspace{-.25cm}
    \end{itemize}
\end{lem}
\begin{pf}
    First observe that we may forget all but one colour used by $c$ and still have an unrealisable constraint matroid. Moving forwards with the proof, we therefore assume that $c$ uses only one colour.

    By \autoref{og}, there is a constraint matroid $(N,c')$, isomorphic to one from the list in \autoref{fig:obst3con}, obtained from $(M,c)$ by contracting edges and deleting colourless edges. Since contracting edges is in our allowed list of operations, assume that $(N,c')$ is obtained from $(M,c)$ just by deleting a set $F\subseteq E(M)$ of colourless edges.

    Let $B$ be a base for $N$ and $B'$ an extension of $B$ to a base for $M$. Consider the matroid $N' := M/(B'\setminus B)$.

    If $N'=N$ then we are done. Otherwise, a graphic realisation $H'$ for $N'$ can be obtained from the graphic realisation $H$ of $N$ by adding colourless edges (including possibly loops) spanning edges which are already there. Since we are allowed to merge parallel edges and contract loops, we may assume that all added edges span a pair of vertices not already spanned. Inspecting \autoref{fig:obst3con}, we see that this always results in a constraint matroid which contains a constraint minor whose underlying matroid is $M(K_4)$, which is itself in the list from \autoref{fig:obst3con}.
\end{pf}

\section{Irreducible obstructions}\label{appendix:evil_class}
There is a clear way to use \autoref{thm:new_red_to_con} to figure out if a constraint matroid is representable. In particular, we check if each block (with an appropriate extension of its colouring) is realisable, and then we try to assemble a realisation for each block into a realisation for the whole constraint matroid. However, \autoref{thm:new_red_to_con} does not tell us beforehand, how to detect whether this will be possible.

The barricade to doing so comes from the assumption in part (ii) of the theorem. In turn, this comes from our methods on pointy constraint matroids. In particular \autoref{lem:using_D_works}, which requires knowing beforehand whether a pointy constraint matroid has a realisation.

Is it possible to characterise when such a realisation exists, using a forbidden substructure?

In Section 5 of \cite{carmesin2017embeddingIV}, Carmesin exhibits an infinite family of unrealisable constraint matroids which do not seem to easily reduce to simpler examples. We generalise that construction here in the following defintion.

Let $(M,c)$ be a connected constraint matroid. Let $F_1,\ldots,F_k$ denote the sets of edges receiving each colour $f_1,\ldots,f_k$ in $c$. We say that $(M,c)$ is a \emph{rainbow arrangement of tightness $t<k$} when there exists a realisation $G$ of $(M,c)$ so that for each subset $\Fcal$ of $\{f_1,\ldots,f_k\}$ of size $t$, there exists a vertex $x$ seeing $\Fcal$ in $G$, but there is not vertex in $G$ seeing the whole of $\{f_1,\ldots,f_k\}$.

Clearly $G$ is not a realisation of the pointy coloured matrid $(M,c,D)$ where $D=\{f_1,\ldots,f_k\}$. One may show however that $G$ realises $(M,c')$ where $c'$ is the extension of $c$ using $D$. Therefore, by \autoref{lem:using_D_works}, $(M,c,D)$ is not realisable.

We can't even fix the problem by introducing more colours in our extension. For some $t<k$ let $c^{(t)}$ be the extension of $c$ be including a new colour for each size-$t$ subset of $\{f_1,\ldots,f_k\}$ and call it the \emph{order $t$ extension of $c$ using $D$}. We can still find numerous rainbow arrangements $(M,c)$ (in particular, of tightness $t$), that exhibit realisations for $(M,c^{(t)})$ but so that $(M,c,D)$ is still not realisable for $D=\{f_1,\ldots,f_k\}$.

\end{document}